\documentclass[10pt]{amsart}
\pdfoutput=1
\usepackage{ mathrsfs }
\usepackage{amsmath, amsthm, amssymb,slashed,stmaryrd}
\usepackage{ifpdf}
\usepackage[pdftex]{graphicx}
\usepackage{tikz}
\usetikzlibrary{matrix,arrows,calc}
\usepackage[all]{xy}

\usepackage[pdftex,plainpages=false,hypertexnames=false,pdfpagelabels]{hyperref}
\usepackage{amsthm}

\theoremstyle{definition}
 \newtheorem{thm}{Theorem}[section]
    
 \newtheorem{cor}[thm]{Corollary}
 \newtheorem{lem}[thm]{Lemma}
 \newtheorem{prop}[thm]{Proposition}
 
 \newtheorem{defn}[thm]{Definition}
 
 \newtheorem{ex}[thm]{Example}
  
 \newtheorem*{thm*}{Theorem}
 \theoremstyle{remark}
 \newtheorem{rmk}[thm]{Remark}

\def\beq{\begin{eqnarray}}
\def\eeq{\end{eqnarray}}
 \newcommand{\bp}{\begin{proof}[Proof]}
 \newcommand{\ep}{\end{proof}}

\DeclareSymbolFont{bbold}{U}{bbold}{m}{n}
\DeclareSymbolFontAlphabet{\mathbbold}{bbold}

\def\image{{\rm image}}

\def\ord{{\rm ord}}
\def\Inf{{\rm Inf}}
\def\M{{\mathcal{M}}}
\def\MM{{\mathfrak{M}}}
\def\sTr{{\sf sTr}}

\def\KO{{\rm KO}}
\def\Ind{{\rm Ind}}
\def\K{{\rm K}}
\def\Mod{{\sf Mod}}

\def\Ch{{\rm Ch}}
\def\CS{{\rm CS}}

\def\dKO{\widehat{\rm KO}{}}
\def\dK{\widehat{\rm K}{}}
\def\Cl{{\rm Cl}}
\def\cCl{\mathbb{C}{\rm l}}
\def\Fred{{\rm Fred}}
\def\cl{{\rm cl}}

\def\H{{\rm H}}

\def\Spin{{\rm Spin}}

\def\U{{\rm U}}

\def\SO{{\rm SO}}

\def\Det{{\rm Det}}

\def\pt{{\rm pt}}
\def\ev{{\rm ev}}
\def\odd{{\rm odd}}

\def\bS{{\mathbb{S}}}
\def\A{{\mathbb{A}}}
\def\Aff{{\mathbb{A}{\rm ff}}}
\def\B{{\mathbb{B}}}
\def\Bc{{\rm{B}}}
\def\R{{\mathbb{R}}}
\def\F{{\mathbb{F}}}

\def\N{{\mathbb{N}}}
\def\id{{{\rm id}}}
\def\C{{\mathbb{C}}}

\def\Z{{\mathbb{Z}}}

\def\End{{\rm End}}

\def\Hom{{\sf Hom}}

\newcommand{\op}{{\sf{op}}}

\newcommand\nc{\newcommand}
\nc\mf\mathfrak
\nc\mc\mathcal
\nc\mb\mathbb

\begin{document}

\title{The families Clifford index and differential KO-theory}
\author{Daniel Berwick-Evans}
\begin{abstract}
Extending ideas of Atiyah--Bott--Shapiro and Quillen, we construct a model for differential $\KO$-theory whose cocycles are families of Clifford modules with superconnection. 
The model is built to accommodate an analytic pushforward for bundles of spin manifolds, affording a differential refinement of Atiyah and Singer's families index.
\end{abstract}

%\date{\today}

\maketitle 
\setcounter{tocdepth}{1}
\tableofcontents

\section{Introduction}

This paper gives a model for differential $\KO$-theory in degree $n\in \Z$ that generalizes a generators and relations description of $\dK^0$, differential complex K-theory in degree~0. 

For a smooth manifold $M$, generators for $\dK^0(M)$ are triples $(V,\A,\phi)$ where~$V\to M$ is a $\C$-vector bundle with hermitian metric, $\A$ is a hermitian superconnection, and $\phi$ is an odd differential form on~$M$. 
These data are subject to an equivalence relation involving a Chern--Simons form; see~\S\ref{sec:FL} for details. 
At first pass, generators for $\dKO^n(M)$ replace~$\C$ by the real Clifford algebra $\Cl_n$: consider $(V,\A,\phi)$ where $V$ is a real, metrized vector bundle with self-adjoint $\Cl_n$-action, $\A$ is a $\Cl_n$-linear metric-preserving connection on~$V$, and $\phi$ is a differential form in degrees $n-1$ mod~4. The appropriate equivalence relation on these data generalizes the Atiyah--Bott--Shapiro equivalence relation on Clifford modules~\cite{ABS}, see also Donovan--Karoubi~\cite{DonovanKaroubi}. Roughly, triples $(V,\A,\phi)$ and $(V',\A',\phi')$ are equivalent when there is an inclusion $g\colon V\hookrightarrow V'$ of metrized vector bundles together with the data of a $\Cl_n\otimes \Cl_{-1}$-action on the orthogonal complement~$V^\perp\subset V'$, plus a condition on a Chern--Simons form. The main subtlety is that this equivalence relation needs to be implemented locally on $M$, meaning an element of $\dKO^n(M)$ is given by an open cover $\{U_\alpha\}$ with $(V_\alpha,\A_\alpha,\phi_\alpha)$ over each $U_\alpha$ and equivalences $g_{\alpha\beta}$ for each $U_\alpha\bigcap U_\beta$ satisfying a cocycle condition that refines the one in the \v{C}ech--de~Rham complex; \S\ref{sec:superconn} and~\S\ref{sec:dKOsec} for details.

\begin{defn}[Sketch of Definition~\ref{defn:diffKO}]\label{defn:differentialKO} Let $M$ be a smooth manifold. Consider the free abelian group on equivalence classes of tuples $\{V_\alpha\to U_\alpha,\A_\alpha,\phi_\alpha,g_{\alpha\beta}\}$, modulo concordances with exact Chern--Simons form. 
Define $\dKO^n(M)$ as the quotient gotten from identifying the direct sum of bundles with addition in~$\dKO^n(M)$; see~\eqref{eq:FLequiv4}. 
\end{defn}

\begin{thm}\label{thm1} The graded ring $\dKO^\bullet(M)$ is naturally isomorphic to the differential KO-theory of~$M$. When $M=\pt$, this description refines the Atiyah--Bott--Shapiro isomorphism that sends a Clifford module to a class in $\KO^\bullet(\pt)$. 
\end{thm}

We emphasize that the real vector bundles $V_\alpha\to U_\alpha$ in Definition~\ref{defn:differentialKO} are finite rank. By contrast, Atiyah and Singer's construction of (non-differential) KO-theory comes from families of $\Cl_n$-linear Fredholm operators acting on infinite-rank bundles~\cite{AtiyahSingerskew}. As we explain below, an index bundle construction permits a comparison between Atiyah and Singer's description of~$\KO^n(M)$ and the one from Theorem~\ref{thm1}; see~\S\ref{sec:unbounded}-\S\ref{subsec:ABSindex} and Proposition~\ref{prop:ClnKO}. As outlined in~\S\ref{sec:motivate}, Definition~\ref{defn:differentialKO} emerges naturally when considering a differential refinement of the index bundle for a family of Clifford linear Dirac operators.
The kernel of such a family is usually not a vector bundle (it need not have locally constant rank), which leads to the subtle gluing data for the index bundle and its differential refinement. 
Definition~\ref{defn:differentialKO} is purpose-built to accommodate this information, providing a natural home for the analytic index in differential $\KO$-theory. 

\begin{thm}\label{thm2}
Let $\pi\colon X\to B$ be a proper family of $n$-dimensional Riemannian spin manifolds and $V\to X$ be a metrized real vector bundle with compatible connection $\nabla$. The families Clifford index determines a degree~$-n$ differential cocycle $\widehat{\pi}_!(V,\nabla)$ on $B$ whose associated class $[\widehat{\pi}_!(V,\nabla)]\in \dKO^{-n}(B)$ refines the analytic pushforward in $\KO$-theory. 
\end{thm}

\subsection{Comparison with other work}
In contrast to differential K-theory for complex bundles, the literature on differential KO-theory is quite sparse. The most thorough treatment to date is by Grady and Sati~\cite{GradySati}. They work with Hopkins and Singer's model for differential cohomology~\cite{HopSing} phrased in terms of sheaves of spectra~\cite{BNV}. When applied to KO, this general construction elucidates structural aspects like differential twists and communicates with powerful homotopical tools, e.g., the differential Atiyah--Hirzebruch spectral sequence. 

The Hopkins--Singer differential refinement of complex K-theory is frequently augmented by geometrically-motivated cocycle models. These more tailored presentations of~$\dK$ are crucial in certain constructions such as the differential analytic index in~\cite{LottGorokhovsky} and differential index theorem in~\cite{FreedLott}. In turn, these constructions allow one to refine and enhance classical results, e.g.,~\cite{LottRZ,FreedQFTK,bunkeindex,BunkeString}.
The analog of these geometric models and their ensuing applications remain under-developed in differential KO-theory. Freed sketches one candidate construction~\cite[page~34]{FreedDirac}, and further suggestions are given by Bunke and Schick~\cite[\S4.9]{BunkeSchick2}; however, these ideas were never fleshed out. More recently, Gomi and Yamashita~\cite{GomiYamashita} developed a model as an enhancement of Karoubi's KO-theory~\cite{KaroubiBook}, with the goal of studying differential refinements of invariants of condensed matter systems~\cite{FreedHopkins}. 

The cocycle model for differential~KO constructed in this paper is motivated by the geometry of index bundles with intended applications in real differential index theory. At the heart of this model is the definition of a $\Cl_n$-bundle (see Definition~\ref{defn:Clnbundle}), which arises as a finite-dimensional families index of $\Cl_n$-linear Dirac operators, see~\S\ref{sec:KOindexsketch}. This geometry resonates with Furuta's vectorial bundles~\cite{Furuta,Gomi} and Henrique's quasi-bundles \cite{AndreBott}, though below we emphasize the connection with the Atiyah--Bott--Shapiro construction~\cite{ABS}. In future work we plan to construct a differential topological index via a realization of the Atiyah--Bott--Shapiro Thom class in the framework of Theorem~\ref{thm1}.

Replacing the real Clifford algebras $\Cl_n$ with complex Clifford algebras $\cCl_n$ leads to analogs of Theorems~\ref{thm1} and~\ref{thm2} for differential complex K-theory. This makes contact with previous results. For example, when $n=0$ and $M$ is compact, Theorem~\ref{thm1} generalizes the generators and relations presentation of $\widehat{\K}^0(M)$ in terms of finite-rank vector bundles with connection and an odd degree differential form; see Remark~\ref{ex:standard}. One also finds overlap with models of a more analytic flavor: applying the index bundle construction (see~\S\ref{sec:motivate}) to the Hilbert bundle model for differential K-theory from~\cite{LottGorokhovsky} leads to a version of Definition~\ref{defn:differentialKO} with $\C$-coefficients, where the odd degree groups in~\cite[\S6]{LottGorokhovsky} can be understood in terms of $\cCl_1$-linear superconnections. 
Finally, Theorem~\ref{thm2} (both in the statement and proof) generalizes Freed and Lott's construction of the differential analytic index for families of spin$^c$-manifolds~\cite{FreedLott}.

Definition~\ref{defn:differentialKO} is also designed to receive a cocycle map from a class of $1|1$-dimensional Euclidean field theories in the style of Stolz and Teichner~\cite{ST04,ST11}. A fundamental assertion of the Stolz--Teichner program is that quantization of such field theories should correspond to the analytic pushforward in KO-theory. It is conjectured that a $2|1$-dimensional generalization supplies an analytic pushforward in the cohomology theory of topological modular forms. Theorems~\ref{thm1} and~\ref{thm2} give a context in which one can start to make some of these ideas precise, see~\cite{DBEChern,DBEEFT,DBEtorsion}.
Owing to its independent interest in geometry, the present paper is written to be independent from this discussion of field theories.

\subsection{Acknowledgements} 
I thank Theo Johnson-Freyd for patiently (and repeatedly) entertaining my questions about superalgebra signs rules. I also thank Yigal Kamel his careful reading and feedback on previous drafts. This work was supported by the National Science Foundation under grant DMS-2205835.

%\setcounter{tocdepth}{1}
%\tableofcontents

\section{The main idea from the families index}\label{sec:motivate}

Definition~\ref{defn:differentialKO} emerges naturally when generalizing Freed and Lott's families differential index~\cite[\S3]{FreedLott} to families of Clifford linear Dirac operators. In this section, we review versions of the families index and describe the relationship with our definitions and results.

\subsection{The families index in $\K$-theory}
Let $X$ be a compact even-dimensional spin${}^c$ manifold. The index of the Dirac operator $\slashed{D}$ on $X$ is 
\beq
\Ind(\slashed{D})=[{\rm ker}(\slashed{D})_\ev]-[{\rm ker}(\slashed{D})_\odd]\in \K(\pt)\simeq \Z,\label{eq:indexC}
\eeq
where we consider ${\rm ker}(\slashed{D})$ as a $\Z/2$-graded $\C$-vector space with even and odd subspaces ${\rm ker}(\slashed{D})_\ev$ and ${\rm ker}(\slashed{D})_\odd$, respectively. 
A family of Riemannian spin${}^c$ manifolds $\pi\colon X\to B$ with compact even-dimensional fibers determines a family of operators $\slashed{D}$ on~$B$. The fiber at $b \in B$ is a Dirac operator $\slashed{D}_b$ on the spin${}^c$ manifold $\pi^{-1}(\{b\})=X_b$. One would like to construct a class in $\K(B)$ from this data. 
Taking the fiberwise index yields a $\Z/2$-graded vector space ${\rm ker}(\slashed{D}_b)$ for each $b\in B$. However, these typically fail to assemble into a $\Z/2$-graded vector bundle, e.g., the dimension of the kernel need not be locally constant. Hence, this naive attempt at a families index does not determine a class in $\K(B)$ generalizing~\eqref{eq:indexC}. 

There are two well-known ways to remedy this. The first is analytic: the space of Fredholm operators represents the functor~$\K(-)$, e.g., see~\cite[\S{III.8}]{LM}. A family of spin${}^c$ manifolds $\pi\colon X\to B$ determines a $B$-family of Dirac operators, which in turn provides a map from $B$ into the space of Fredholm operators. The homotopy class of this map defines the \emph{families analytic index} $\Ind(\slashed{D}) \in \K(B)$. The second approach is geometric: consider the open cover~$\{U_\lambda\}_{\lambda\in \R_{>0}}$ of~$B$ defined as
\beq
U_\lambda:=\{b\in B\mid \lambda\notin {\rm Spec}(\slashed{D}_b^2)\},\label{eq:Ulambda}
\eeq
where ${\rm Spec}(\slashed{D}_b^2)$ denotes the spectrum of the self-adjoint, nonnegative Dirac Laplacian $\slashed{D}_b^2$. For $b\in U_\lambda$, let $\mathcal{H}^{<\lambda}_b$ denote the sum of eigenspaces of $\slashed{D}_b^2$ with eigenvalues less than~$\lambda$. Ellipticity of $\slashed{D}_b^2$ implies that the $\mathcal{H}^{<\lambda}_b$ fit together into a finite rank, smooth, $\Z/2$-graded vector bundle $\mathcal{H}^{<\lambda}\to U_{\lambda}$. Hence we obtain a $\K$-theory class $[\mathcal{H}^{<\lambda}]=[\mathcal{H}^{<\lambda}_\ev]-[\mathcal{H}^{<\lambda}_\odd]\in {\rm K}(U_\lambda)$. On overlaps $U_\lambda\bigcap U_\mu$ with $\lambda<\mu$, there is an evident inclusion of $\Z/2$-graded vector bundles,
\beq
g_{\lambda\mu}\colon \mathcal{H}^{<\lambda}\hookrightarrow \mathcal{H}^{<\mu}. \label{eq:Hilbcocycle}
\eeq
The restriction of $\slashed{D}$ to the orthogonal complement of the inclusion~\eqref{eq:Hilbcocycle} determines an invertible odd endomorphism $e_{\lambda\mu}$ that witnesses an equality of $\K$-theory classes $[\mathcal{H}^{<\lambda}]=[\mathcal{H}^{<\mu}]\in \K(U_\lambda\bigcap U_\mu)$, see~\eqref{eq:CLminus1action}. 
In fact, the inclusions~\eqref{eq:Hilbcocycle} and odd isomorphisms $e_{\lambda\mu}$ satisfy a cocycle condition on $U_\lambda\bigcap U_\mu\bigcap U_\mu$. This constructs a well-defined class 
\beq
{\rm Ind}(\slashed{D}):=[\{\mathcal{H}^{<\lambda}\to U_\lambda, g_{\lambda\mu},e_{\lambda\mu}\}] \in  {\rm K}(B)\label{eq:KUindexbundle}
\eeq
equal to the families analytic index of $\pi\colon X\to B$; see~\cite[\S7.2]{FreedGeoDirac}, \cite[\S9.5]{BGV}.

The geometric description~\eqref{eq:KUindexbundle} of the index also leads to a natural refinement in the differential $\K$-theory of~$B$ that uses the \emph{Bismut superconnection}, denoted~$\B$ and recalled in~\eqref{eq:Bismutsuper} below. The restriction of the degree~1 part of~$\B$ to $\mathcal{H}^{<\lambda}$ determines ordinary connections $\nabla^{<\lambda}$ on each finite-rank bundle~$\mathcal{H}^{<\lambda}\subset \mathcal{H}|_{U_\lambda}$ with Chern forms
\beq\label{eq:initialChernform}
\Ch(\nabla^{<\lambda})\in \Omega^{\rm even}(U_\lambda;\R[u^{\pm 1}]),\qquad |u|=-2. 
\eeq
Properties of the Bismut superconnection further determine a unique extension of the data~\eqref{eq:initialChernform} to a cocycle in the \v{C}ech--de~Rham complex of $M$ relative to the open cover~\eqref{eq:Ulambda} (see Proposition~\ref{prop:CechdeRham}),
\beq\label{eq:initialChernformprime}
&&\Ch(\B,\{U_\lambda\}_{\lambda\in \R_{>0}})\in \Omega^{\bullet}(U_*;\R[u^{\pm 1}]),\qquad (d+\delta)\Ch(\B,\{U_\lambda\}_{\lambda\in \R_{>0}})=0.
\eeq 
By construction, this \v{C}ech--de~Rham cocycle represents the Chern character of the index bundle. By the Riemann--Roch theorem, another representative is the integral along the fibers
\beq\label{eq:indexdensity1}
\int_{X/B}\Ch(\nabla^V)\wedge {\rm Todd}(\nabla^{X/B}) \in \Omega^\bullet_\cl(B;\R[\alpha,\alpha^{-1}])
\eeq
of the Chern form of $V$ modified by the ${\rm Todd}(\nabla^{X/B})$, the Todd form of the spin${}^c$-family. 

The families differential index comes from specifying a coboundary in the \v{C}ech--de~Rham complex that mediates between the cocycles~\eqref{eq:initialChernformprime} and~\eqref{eq:indexdensity1}. 
%We recall (e.g.,~\cite[\S8]{BottTu}) that the \v{C}ech differential $\delta$ requires an ordering on the indexing set for the cover, which in this case is inherited from the standard ordering on $\R_{>0}$.
Adapting methods of Bismut~\cite{Bismutindex}, Bismut--Cheeger~\cite{BismutCheeger}  and Freed--Lott~\cite[\S3]{FreedLott}, there are forms $\eta_\lambda\in \Omega^{\rm odd}(U_\lambda,\R[u,u^{-1}])$ satisfying
\beq\label{Eq:basicetaform}
&&\int_{X/B}\Ch(\nabla^V)\wedge {\rm Todd}(X/B)-\Ch(\B,\{U_\lambda\})=(d+\delta)\{\eta_\lambda\}\in \Omega^\bullet(U_*;\R[u,u^{-1}]),
\eeq
where $d$ is the de~Rham differential and $\delta$ is the \v{C}ech differential, i.e., the $\eta_\lambda$ determine the desired coboundary in the \v{C}ech--de~Rham complex. This determines a differential $\K$-cocycle 
\beq\label{eq:diffindex}
\widehat{\Ind}(\slashed{D})=[\{(\mathcal{H}^{<\lambda}\to U_\lambda,\nabla^{<\lambda},\eta_\lambda),(g_{\lambda\mu},e_{\lambda\mu})\}]\in \widehat{\K}(B)
\eeq
representing the differential index. 

\subsection{The families index in $\KO$-theory}\label{sec:KOindexsketch}
Now consider the Clifford linear Dirac operator~$\slashed{D}$ on an $n$-dimensional Riemannian spin manifold. The kernel of $\slashed{D}$ is a finite-rank Clifford module~\cite[II.7]{LM}. The equivalence class of this module under the Atiyah--Bott--Shapiro map gives a class in~$\KO^{-n}(\pt)$~\cite{ABS}
\beq
&&\Ind(\slashed{D}):=[{\rm ker}(\slashed{D})] \in  \KO^{-n}(\pt).\label{eq:KOindex}
\eeq
This is the \emph{Clifford index} of $\slashed{D}$. We refer to~\S\ref{sec:Cliffordindex} for further details on the Clifford linear Dirac operator and~\S\ref{sec:ABSreview} for a review of the Atiyah--Bott--Shapiro construction. 

When generalizing~\eqref{eq:KOindex} to families of $n$-dimensional spin manifolds $\pi\colon X\to B$, the kernels of the fiberwise Clifford linear Dirac operators fail to determine a smooth bundle of Clifford-modules for the same reasons as in the spin$^c$ case. One again has an analytic remedy: the space $\Fred_{n}$ of Clifford linear Fredholm operators represents the functor~$\KO^{-n}(-)$~\cite{AtiyahSingerskew},~\cite[III.10]{LM}. A family of Clifford linear Dirac operators associated to a bundle of spin manifolds $\pi\colon X\to B$ determines a map $B\to \Fred_{n}$ whose homotopy class defines the \emph{families analytic index} 
\beq
\Ind(\slashed{D})\in [B,\Fred_{n}]\simeq \KO^{-n}(B).\label{eq:familiesanalytic}
\eeq
The main goal of this paper is to develop the complementary geometric approach to the families Clifford index, in parallel to the index bundle for the spin$^c$-Dirac operator. In the remainder of this section we sketch the main ideas. 

\begin{rmk} 
There are several sign conventions to fix when working with Clifford linear Dirac operators. Our conventions (explained in Remarks~\ref{rmk:sign1} and~\ref{rmk:signs}) lead to the kernel of the Clifford linear Dirac operator being a (left) $\Cl_{-n}$-module. 
\end{rmk}

For a family $\slashed{D}$ of $\Cl_{-n}$-linear Dirac operators, define an open cover~\eqref{eq:Ulambda} as before: consider the subspace where $\lambda\in \R_{>0}$ is not an eigenvalue of $\slashed{D}^2$. Taking the direct sum of eigenspaces with eigenvalues $<\lambda$ yields finite-rank bundles of Clifford modules $\mathcal{H}^{<\lambda}\to U_\lambda$ over each component. Smoothness of this bundle follows from ellipticity of the Dirac Laplacian~$\slashed{D}^2$. On an overlap $U_\lambda\bigcap U_\mu$, the restriction of $\slashed{D}$ to the orthogonal complement to the inclusion of bundles~\eqref{eq:Hilbcocycle} determines an invertible odd endomorphism that commutes with the $\Cl_{-n}$-action, see~\eqref{eq:CLminus1action}. This odd endomorphism generates a $\Cl_{-1}$-action, leading to a $\Cl_{-n}\otimes \Cl_{-1}\simeq \Cl_{-n-1}$-action on the orthogonal complement extending the prior $\Cl_{-n}$-action. The existence of such an extension implies that the bundles of Clifford modules over $U_\lambda\bigcap U_\mu$ represent the same class in $\KO$-theory. In fact, the analogous information~\eqref{eq:KUindexbundle} in the Clifford linear case determines a class in $\KO^{-n}(B)$. 

To summarize, the index bundle for a family of $\Cl_{-n}$-linear Dirac operators is data:
\begin{enumerate}
\item[(i)] an open cover with finite-rank $\Cl_{-n}$-module bundles $\mathcal{H}^{<\lambda}\to U_\lambda$; 
\item[(ii)] for each intersection $U_\lambda\bigcap U_\mu$, inclusions~\eqref{eq:Hilbcocycle} of $\Cl_{-n}$-modules together with the data of a $\Cl_{-n-1}$-action  on the orthogonal complement extending the $\Cl_{-n}$-action.
\end{enumerate} 
These data satisfy a cocycle condition. This essentially recovers Definition~\ref{defn:Clnbundle} of a $\Cl_{-n}$-bundle. In turn, $\Cl_{-n}$-bundles can be used to construct $\KO^{-n}(M)$ (Proposition~\ref{prop:ClnKO}), yielding a version of Theorem~\ref{thm1} in (non-differential) $\KO$-theory. 

A differential refinement comes from generalizing the differential index~\eqref{eq:diffindex} in complex K-theory. To start, the Bismut superconnection $\B$ defined in~\eqref{eq:Bismutsuper} for a family of Clifford linear Dirac operators is automatically $\Cl_{-n}$-linear. It restricts to an ordinary $\Cl_{-n}$-linear connection $\nabla^{<\lambda}$ on each~$\mathcal{H}^{<\lambda}$ with Chern forms analogous to~\eqref{eq:initialChernform} comprising part of the data of a \v{C}ech--de~Rham cocycle~\eqref{eq:initialChernformprime} representing the Pontryagin character of the index bundle. Replacing the Todd form in~\eqref{eq:indexdensity1} by the $\hat{A}$-form, the remaining data of the differential index comes from $\eta_\lambda$-forms determining a \v{C}ech--de~Rham coboundary as in~\eqref{Eq:basicetaform}. These are constructed from Clifford linear generalizations of the standard constructions. With these definitions in hand, the differential $\KO$-index is defined essentially the same as in~\eqref{eq:diffindex},
\beq\label{eq:diffindex2}
\widehat{\rm Ind}(\slashed{D})=[\{(\mathcal{H}^{<\lambda}\to U_\lambda,\nabla^{<\lambda},\eta_\lambda),(g_{\lambda\mu},e_{\lambda\mu})\}]\in \dKO^{-n}(B)
\eeq
where $\mathcal{H}^{<\lambda}\to U_\lambda$ is a bundle of $\Cl_{-n}$-modules, $g_{\lambda\mu}$ is an inclusion of Clifford module bundles, $e_{\lambda\mu}$ extends a $\Cl_{-n}$-action on the orthogonal complement of the inclusion to a $\Cl_{-n-1}$-action, $\nabla^{<\lambda}$ is a Clifford linear connection on $\mathcal{H}^{<\lambda}$, and the forms $\{\eta_\lambda\}$ determine a \v{C}ech--de~Rham coboundary. 

Definition~\ref{defn:superconngeneral} of a differential cocycle is designed so that the differential index~\eqref{eq:diffindex2} is an immediate example. Constructing a map from data like~\eqref{eq:diffindex2} to differential $\KO$-theory basically follows the same route as in complex K-theory. The main work in Theorem~\ref{thm1} is to show that this map is surjective; this comes down to a computation (Proposition~\ref{prop:diffsusp}). Theorem~\ref{thm2} then follows from the differential index construction sketched above. 

\section{Clifford modules}

In this section we review the theory of Clifford modules. Our focus is on real Clifford algebras and real Clifford modules, with a few pieces of the theory over~$\C$. For definitions and constructions that apply to either real or complex coefficients, we use $\F$ to denote either~$\R$ or~$\C$. 

\subsection{$\Z/2$-graded $*$-algebras} 
The category of \emph{real}, respectively \emph{complex}, \emph{super vector spaces} is the category whose objects are $\Z/2$-graded vector spaces over~$\R$, respectively~$\C$, and morphisms are linear maps that preserve the grading. We promote this to a symmetric monoidal category via the graded tensor product, where the adjective ``super" refers to the sign in the braiding isomorphism
$$
V\otimes_\F W\stackrel{\sim}{\to} W\otimes_\F V \qquad v\otimes w\mapsto (-1)^{|v||w|} w\otimes v \quad v\in V, \ w\in W
$$
where $v$ and $w$ are homogeneous elements of degree $|v|,|w|\in \Z/2$. For a super vector space~$V$, let $V_\ev\oplus V_\odd$ denote the direct sum decomposition of $V$ into its even and odd subspaces. The \emph{grading involution} is the linear map $(-1)^{\sf F}\colon V \to V$ that acts by $+1$ on $V_\ev$ and~$-1$ on $V_\odd$.\footnote{The notation $(-1)^{\sf F}$ is from physics where ${\sf F}$ stands for ``fermion", e.g., see~\cite[\S1]{susymorse}.}  Let $\Pi V$ denote the parity reversal of~$V$, i.e., $(\Pi V)_\pm=V_\mp$. Throughout, $\F^{n|m}$ will denote the super vector space with~$(\F^{n|m})_\ev=\F^n$ and~$(\F^{n|m})_\odd=\F^m$.

A \emph{superalgebra} over $\F$ is an algebra object in super vector spaces over $\F$. The \emph{super commutator} of elements in a superalgebra is defined as
\beq
[a,b]=ab-(-1)^{|a||b|}ba\qquad a,b\in A,\label{eq:supercommutator}
\eeq
for homogenous elements $a$ and $b$. The super commutator satisfies a graded version of the Jacobi identity. 

\begin{ex}\label{ex:End}
For a super vector space $V$, let $\End(V)$ denote the set of all linear maps $V\to V$ (not necessarily grading-preserving), where grading-preserving maps are regarded as even and grading-reversing maps are odd. Composition of linear maps endows $\End(V)$ with the structure of a superalgebra.
\end{ex}

\begin{rmk}
Generalizing the previous example, for super vector spaces $V$ and $W$ one can consider the super vector space~$\underline{\Hom}(V,W)$ of all linear maps $V\to W$, graded according to whether a linear map is grading preserving versus grading reversing. This endows the category of super vector spaces with a closed monoidal structure, i.e., $-\otimes V$ is left adjoint to $\underline\Hom(V,-)$. 
\end{rmk}

A \emph{super trace} on a superalgebra is a map of super vector spaces 
\beq
\sTr_A\colon A\to \F\qquad \sTr_A([a,b])=0\label{eq:superalgsupertrace}
\eeq
that vanishes on super commutators, and hence is determined by a linear map~$A/[A,A]\to \F$. The \emph{opposite} of a superalgebra $A$ is a superalgebra $A^\op$ with the same underlying super vector space and multiplication 
$$
a\cdot_{\op} b:=(-1)^{|a||b|}b\cdot a
$$
where $b\cdot a$ is the multiplication in $A$. For superalgebras $A$ and $B$, let $A\otimes B$ denote the superalgebra with multiplication
$$
(a\otimes b)\cdot (a'\otimes b')=(-1)^{|b||a'|} aa'\otimes bb'. 
$$

\begin{ex}
The \emph{super trace} of a linear endomorphism $T\colon V\to V$ is defined as
\beq
\sTr(T):= {\sf Tr}((-1)^{\sf F}\circ T)\label{eq:supetrace}
\eeq
where ${\sf Tr}$ is the ordinary trace. The super trace determines a trace on the superalgebra~$\End(V)$ from Example~\ref{ex:End}. 
\end{ex}
A \emph{$*$-structure} on a superalgebra is a $*$-structure on its underlying (ungraded) algebra, i.e., an anti-involution of the underlying algebra (that is conjugate-linear when $\F=\C$). 

\begin{ex}
The usual transpose determines a $*$-structure on $\End(V)$. 
\end{ex}

\begin{rmk}
A closely related notion is that of a \emph{$*$-superalgebra}, which is a superalgebra together with an involutive \emph{graded} anti-involution, meaning a homomorphism to the opposite algebra in the graded sense. There is a dictionary relating superalgebras with $*$-structures to $*$-superalgebras~\cite[pages~89-92]{strings1}. The signs associated with $*$-superalgebras (while perhaps more categorically pleasing) compare awkwardly with the standard conventions in index theory. This is why we have chosen to work with ordinary $*$-structures. 
\end{rmk} 

\subsection{Modules over superalgebras}\label{sec:supermod}
Given a superalgebra $A$, let ${}_A\Mod$ and $\Mod_A$ denote the categories of left and right $A$-modules, respectively, whose objects are super vector spaces with a (graded) $A$-action. These categories have symmetric monoidal structures from the direct sum of $A$-modules. Similarly, let ${}_A\Mod_B$ denote the category of $A-B$-bimodules. There are canonical equivalences of categories 
\beq
{}_{A\otimes B}\Mod\simeq {}_A \Mod_{B^\op} \simeq \Mod_{ A^\op\otimes B^\op}\label{eq:leftright}
\eeq 
where a left $A\otimes B$-bimodule $W$ is sent to the $A-B^\op$-bimodule with the same underlying super vector space and action
$$
a\cdot w\cdot b:=(-1)^{|w||b|} (a\otimes b) \cdot w
$$
with similar formulas defining a right $A^\op\otimes B^\op$-module. We often use notation like ${}_AW$ for an object of ${}_A\Mod$, and ${}_AW_B$ for an object of ${}_A\Mod_B$. 

A super trace on a superalgebra $A$ induces a super trace on the category of (finite-dimensional) $A$-modules. For right $A$-modules, this super trace comes from the composition
\beq
\End_A(W)\stackrel{\sim}{\leftarrow} W\otimes_A W^\vee \to A/[A,A]\stackrel{\sTr_A}{\to}  \F \label{eq:HattoriStallings}
\eeq
where $W^\vee:={\rm Hom}_A(W,A)$ is the dual (left) $A$-module, and the map to $A/[A,A]$ is the Hattori--Stallings trace~\cite{Hattori,Stallings} induced by the evaluation pairing,
$$
{\rm eval}\colon W\otimes_\F W^\vee \to A. 
$$
One obtains traces for bimodules (and left modules) from traces on right modules using~\eqref{eq:leftright}. 

Super algebras $A$ and $B$ are \emph{Morita equivalent}  if there exist bimodules ${}_AW_B$ and ${}_BV_A$ such that ${}_AW\otimes_B V_A\simeq {}_AA_A$ and ${}_BV\otimes_A W_B\simeq {}_BB_B$, where $A$ and $B$ are regarded as bimodules over themselves. For Morita equivalent algebras, there are equivalences of categories
$$
({}_B V\otimes_A -) \colon {}_A\Mod\leftrightarrows {}_B \Mod \colon ({}_AW\otimes_B-).
$$

\begin{rmk} Super algebras, bimodules, and bimodule maps are the objects, 1-morphisms and 2-morphisms of the (super) Morita bicategory, e.g., see~\cite[\S1.27]{Freedalg}. Objects in this bicategory are 1-isomorphic when they are Morita equivalent. We do not explicitly use this bicategorical structure below, but we do employ some of the terminology. In particular, a Morita equivalence will be referred to as an \emph{invertible} bimodule in reference to its role as an invertible 1-morphism in the Morita bicategory. 
\end{rmk}

\subsection{Clifford algebras}

For a real (ungraded) vector space $V$ with quadratic form $q$, the \emph{Clifford algebra} is the quotient of the tensor algebra of $V$ 
\beq
\Cl(V,q):=T(V)/\langle v\otimes v+q(v)\rangle.\label{eq:ClVdef}
\eeq
This is a superalgebra with grading involution determined by $v\mapsto -v$ for $v\in V\subset \Cl(V,q)$. 

\begin{rmk} \label{rmk:sign1}
Sign conventions in~\eqref{eq:ClVdef} vary; ours are consistent with~\cite{ABS,LM,ST04} but differ from~\cite[IX.9]{Bourbakialgebra}. 
\end{rmk} 

Let $\Cl_{n,m}$ denote the Clifford algebra associated to the real vector space $V=\R^{n+m}= \R^n\oplus \R^m$ for quadratic form given by the direct sum of the standard positive definite inner product on $\R^n$ and the standard negative definite inner product on $\R^m$. Hence, we have the generators and relations presentation,
\beq
&&\Cl_{n,m}=\langle f_1,\dots, f_n, e_1,\dots, e_m \mid  [f_i,f_j]=-2\delta_{ij},\ [e_i,e_j]=2\delta_{ij}, \ [f_i,e_j]=0\rangle \label{eq:gennotation}
\eeq
where the $f_i\in \R^n\subset \Cl_{n,m}$ and $e_i\in \R^m\subset \Cl_{n,m}$ are odd, and the relations all involve the super commutator~\eqref{eq:supercommutator}. Since quadratic forms over $\R$ are classified by their signature, there exists an isomorphism $\Cl(V,q)\simeq \Cl_{n,m}$ for any $(V,q)$ where $q$ has signature $(n,m)$. 

We equip $\Cl_{n,m}$ with a $*$-structure determined by the values on generators,
\beq
e_i^*=e_i, \qquad f_i^*=-f_i. \label{eq:starstructure}
\eeq
We also adopt the notation
\beq
&&\Cl_n:=\left\{\begin{array}{ll} \Cl_{n,0} &\phantom{BB} n\ge 0 \\ \Cl_{0,n} &\phantom{BB} n<0\end{array}\right.\nonumber
\eeq
and
\beq
\Gamma_{n,m}=2^{-(n+m)/2}f_1f_2\cdots f_ne_1e_2\cdots e_m,\qquad \Gamma_n:=\left\{\begin{array}{ll} \Gamma_{n,0} & n\ge 0\\ \Gamma_{0,n} & n<0. \end{array}\right.\label{eq:Gammadef}
\eeq
The element $\Gamma_{n,m}$ depends on a choice of orientation of $\R^{n}\oplus \R^m$, and changing the orientation introduces a sign. There are canonical isomorphisms
\beq
&&\Cl_n^\op\simeq \Cl_{-n},\quad \Cl_n\simeq \Cl_1^{\otimes n},\quad  \Cl_{n}\otimes \Cl_{m}\simeq \Cl_{n+m},\quad  \Cl_{n}\otimes \Cl_{-m}\simeq \Cl_{n,m},\label{eq:isosofCliff}
\eeq
for $n,m$ both nonnegative or both nonpositive. The isomorphism $\Cl_n^\op\simeq \Cl_{-n}$ is determined by the identity map on underlying super vector spaces, so below we will write~$\Cl_n^\op= \Cl_{-n}$. The other isomorphisms critically use that $\otimes$ is the graded tensor product of superalgebras. 

Define the complex Clifford algebras by $\cCl_{n,m}:=\Cl_{n,m}\otimes \C$. These have $*$-structures by taking the $\C$-antilinear extension of the assignments~\eqref{eq:starstructure}. There are isomorphisms $\cCl_{n,m}\simeq \cCl_{n+m}\simeq \cCl_{-n-m}$ implemented by multiplying the appropriate generators by~$i=\sqrt{-1}$. These isomorphisms do not preserve the elements $\Gamma_{n,m}\in \Cl_{n,m}\subset \cCl_{n,m}$ in~\eqref{eq:Gammadef}. 

\subsection{Clifford modules}
A (real) \emph{Clifford module} is a graded left module over $\Cl(V,q)$. Right Clifford modules and Clifford bimodules are defined analogously. The first isomorphism in~\eqref{eq:isosofCliff} determines equivalences of categories 
\beq
{}_{\Cl_n\otimes \Cl_{-m}}\Mod\simeq {}_{\Cl_n}\Mod_{\Cl_m}\simeq \Mod_{\Cl_{-n}\otimes \Cl_m}\label{eq:Cliffmods}
\eeq
which will allow us to translate any Clifford (bi)module into a left Clifford module. 

An \emph{inner product} on a real super vector space $W$ is an inner product on the ungraded vector space underlying $W$ such that the even and odd subspaces of $W$ are orthogonal. We take a similar definition for a hermitian inner product on a complex super vector space. 

\begin{rmk} 
Sign conventions in differential geometry often use a mix of the graded and ungraded tensor product for $\Z/2$-graded Hilbert spaces, Clifford algebras, and differential forms valued in sections of a super vector bundle, e.g., conventions in~\cite{BGV}. Making these conventions consistent with each other uses an equivalence between inner products on super vector spaces defined via the graded and ungraded tensor product. This equivalence involves signs and factors of $\sqrt{-1}$ (see~\cite[pages~90-91]{strings1}). 
Such factors appear when the two conventions collide. For example, self-adjoint superconnections inherit signs from this translation, see Definition~\ref{defn:sasuperconn}. 
\end{rmk}

A left Clifford module $\rho\colon \Cl(V,q)\to \End(W)$ is \emph{self-adjoint}~\cite[Definition~3.3]{BGV} if~$W$ is equipped with an inner product and
\beq
\rho(v^*)=\rho(v)^\dagger\label{eq:superadjoint}
\eeq
for all $v\in V\subset \Cl(V,q)$, using the $*$-structure~\eqref{eq:starstructure} and where $(-)^\dagger$ is the adjoint on $\End(W)$. Self-adjoint right modules and self-adjoint bimodules are defined using the equivalences~\eqref{eq:Cliffmods}. Given a self-adjoint left $\Cl(V,q)$-module $W$, a \emph{Clifford linear endomorphism} is a linear map $T\colon W\to W$ (not necessarily grading preserving) that graded commutes with the Clifford action, i.e., for all $v\in V\subset \Cl(V,q)$,
$$
T(\rho(v)\cdot w)=(-1)^{|T|}\rho(v)\cdot T(w), \quad \iff \quad [T,\rho(v)]=0
$$
where $[-,-]$ denotes the super commutator~\eqref{eq:supercommutator} on $\End(W)$. A Clifford linear odd endomorphism is \emph{self-adjoint} when it is self adjoint with respect to the inner product on $W$. Note that the grading involution on~$W$ is always self-adjoint. 

Complex Clifford modules are defined analogously, where self-adjoint modules are defined relative to a hermitian inner product. 

\begin{ex}\label{ex:Morita1}
Define the left $\Cl_{1,1}$-module whose underlying super vector space is $\R^{1|1}$ with action determined by
\beq
f\mapsto\left[\begin{array}{cc} 0 & 1 \\ -1 & 0\end{array}\right],\qquad e \mapsto\left[\begin{array}{cc} 0 & 1 \\ 1 & 0\end{array}\right]\label{eq:Cl11bimod}
\eeq
using the notation~\eqref{eq:gennotation} for generators. For the standard inner product this $\Cl_{1,1}$-action is self-adjoint. The assignment~\eqref{eq:Cl11bimod} induces an isomorphism of superalgebras, $\Cl_{1,1}\simeq \End(\R^{1|1})$. This implies that $\Cl_{1,1}$ is Morita equivalent to $\R$ with Morita bimodule $\R^{1|1}$. Using~\eqref{eq:isosofCliff}, we find that $\Cl_n\otimes \Cl_m$ is Morita equivalent to $\Cl_{n+m}$ for any $n,m\in \Z$. From~\eqref{eq:Cliffmods}, we also see that~\eqref{eq:Cl11bimod} determines an invertible $\Cl_1-\Cl_1$-bimodule. The isomorphism $\Cl_{1,1}\simeq \End(\R^{1|1})$ shows that there are two isomorphism classes of invertible $\Cl_1-\Cl_1$-bimodules, corresponding under Morita equivalence to the line $\R$ with its two possible gradings. The left $\Cl_{1,1}$-module from~\eqref{eq:Cl11bimod} corresponds to $\Cl_1$ as a bimodule over itself. Finally, the isometric automorphisms of $\Cl_1$ as a bimodule is the group $\Z/2=\Spin(1)$, acting by $\{\pm 1\}$ on $\R^{1|1}\simeq \Cl_1$. 
\end{ex}

\begin{ex}\label{ex:Moritann}
Taking tensor powers of the previous example yields a Morita equivalence between $\Cl_{n,n}$ and $\R$, determined by $\Cl_n$ as an invertible $\Cl_n$-bimodule for $n\in \N$. The group $\SO(n)$ acts on $\Cl_{\pm n}$ through algebra automorphisms via its standard action on generators $\R^n\subset \Cl_{\pm n}$. This gives a categorical action of $\SO(n)$ on the groupoid of invertible $\Cl_{\pm n}$-bimodules, given by twisting the left $\Cl_{\pm n}$-action by an algebra automorphism of $\Cl_{\pm n}$. Categorical actions are additional data, which in this case amounts to the $\Spin(n)$-extension of~$\SO(n)$. To see this, we first observe that by Schur's lemma there exists an isomorphism between $\Cl_{\pm n}$ as a bimodule over itself for the standard action, and the bimodule that twists the left action by an element of $\SO(n)$. This isomorphism is given by a (metric-preserving) scalar, i.e., $\{\pm 1\}$. Hence, the group of automorphisms of $\Cl_{\pm n}$ as a bimodule is a $\Z/2$-extension of $\SO(n)$, which turns out to be the spin double cover, e.g., see~\cite[Definition 2.3.1]{ST04}. 
\end{ex} 

\begin{ex}\label{ex:CL4}
Let $\mathbb{H}$ denote the quaternions, and fix the isomorphism of vector spaces $\mathbb{H}\simeq \R^4$ sending $\{1,i,j,k\}$ to the standard basis of $\R^4$. Define a $\Cl_4$-module whose underlying super vector space is $\R^{4|4}$ with action determined by the action of generators $\R^4\subset \Cl_4$
\beq\label{eq:CL4}
&&\R^4 \to \End(\mathbb{R}^{4|4})\qquad q\mapsto \left[\begin{array}{cc} 0 & q \\ -\bar q & 0\end{array}\right],\quad q\in \mathbb{H}
\eeq
where $q\in \mathbb{H}$ is an element of the quaternions acting on $\R^4$ via the left action of the quaternions on themselves, and the block diagonal matrix is for the decomposition $\R^{4|4}\simeq \R^{4|0}\oplus \R^{0|4}$, i.e., into even and odd subspaces. 
Checking the Clifford relations, the assignment~\eqref{eq:CL4} extends to an injective homomorphism of superalgebras $\Cl_4\hookrightarrow \End(\R^{4|4})$ with image the $\mathbb{H}$-linear endomorphisms of $\R^{4|4}\simeq \mathbb{H}\oplus\mathbb{H}$. 
\end{ex}

\begin{ex}\label{ex:Morita2} Let $\mathbb{O}$ denote the octonians, and fix an isomorphism of vector spaces $\mathbb{O}\simeq \R^8$ by sending the generators of the octonians to the standard basis of $\R^8$. Define a $\Cl_8$-module whose underlying super vector space is $\R^{8|8}$ with an action determined by the action of generators $\R^8\subset \Cl_8$
\beq\label{eq:Morita2}
\R^8 \to \End(\mathbb{R}^{8|8})\qquad x\mapsto \left[\begin{array}{cc} 0 & x \\ -\bar x & 0\end{array}\right],\quad x\in \mathbb{O}
\eeq
where $x\in \mathbb{O}$ acts on $\R^8$ (on the left) via the identification with left multiplication in the octonions, and the block diagonal matrix is for the decomposition $\R^{8|8}\simeq \R^{8|0}\oplus \R^{0|8}$.
The assignment~\eqref{eq:Morita2} leads to an injective map $\Cl_8\hookrightarrow \End(\R^{8|8})$; since $\dim(\Cl_8)=2^8=16\cdot 16=\dim(\End(\R^{8|8}))$, this map is an isomorphism. 
Together with Example~\ref{ex:Morita1}, this induces the famous Morita equivalences between the Clifford algebras $\Cl_{n+8}$ and $\Cl_n$ for all~$n$. We also note the isomorphism of $\Cl_8\simeq \Cl_4\otimes \Cl_4$-modules~\cite[page 70]{LM}
\beq
({}_{\Cl_4}\R^{4|4}) \otimes_\R ({}_{\Cl_4}\R^{4|4}) \simeq ({}_{\Cl_8} \R^{8|8})^{\oplus 4} \label{eq:KOrelation}
\eeq
where ${}_{\Cl_4}\R^{4|4}$ is from the previous example. 
\end{ex} 
%https://mathoverflow.net/questions/194279/explicit-isomorphism-between-cl8-and-mathbbr16

\begin{ex} \label{ex:cMorita2} The complexification of the module from Example~\ref{ex:Morita1} gives a $\Cl_{1,1}\otimes \C\simeq \cCl_2$-module that is an invertible $\C-\cCl_2$-bimodule. Explicitly, the action of $\cCl_2$ on $\C^{1|1}$ is determined by the action of the generators
\beq
f_1\mapsto\left[\begin{array}{cc} 0 & 1 \\ -1 & 0\end{array}\right],\qquad f_2 \mapsto\left[\begin{array}{cc} 0 & i \\ i & 0\end{array}\right], \label{eq:Cl2bimod}
\eeq
inducing an isomorphism $\cCl_2\stackrel{\sim}{\to} \End(\C^{1|1})$ and hence a Morita equivalence between $\cCl_2$ and $\C$. 
\end{ex}

\subsection{Super traces for Clifford modules}\label{sec:supertraces}
There is a canonical isomorphism~\cite[Proposition~3.21]{BGV}
$$
\Cl_n/[\Cl_n,\Cl_n]\simeq \Cl_n/\Cl_{n-1}\simeq  \Lambda^n\R^n= \Det(\R^n)
$$ 
with the determinant line bundle of $\R^n$. 
Hence, a (nontrivial) super trace~\eqref{eq:superalgsupertrace} on $\Cl_n$ is equivalent to the data of a trivialization of the determinant line on $\R^n$, i.e., a choice of volume form. Any such volume form is equivalent to an element $c\cdot \Gamma_n$ for $c\in \R$ and $\Gamma_n$ from~\eqref{eq:Gammadef}. This gives a formula for the associated super trace,
\beq
\sTr_{\Cl_n}(a):= c\cdot \sTr(\Gamma_n \circ a)\qquad c\in \R\label{eq:ClsTr}
\eeq
where the super trace on the right is taken for the action of $\Gamma_n \circ a$ on the vector space $\Lambda^\bullet\R^n\simeq \Cl_n$ coming from the left action of the Clifford algebra on itself. 

After choosing a scalar~$c$, equations~\eqref{eq:HattoriStallings} and~\eqref{eq:ClsTr} determine a super trace on free $\Cl_n$-modules. We extend this to a super trace on arbitrary Clifford modules by writing a finite-rank module as a summand of a free module. The resulting super trace
$$
\End_{\Cl_n}(W)\to \Det(\R^n)\stackrel{\sim}{\to} \R
$$ 
inherits all the properties of the Hattori--Stallings trace~\eqref{eq:HattoriStallings}, e.g., it vanishes on super commutators of Clifford linear endomorphisms. When comparing with the Pontryagin character in KO-theory, it turns out to be convenient to modify the above to include a formal parameter; compare~\cite[\S2.5]{FreedLott}. 

\begin{defn}\label{defn:Cliffordsupertrace}
Given a left $\Cl_{n,m}$-module $W$, the \emph{Clifford super trace} of a~$\Cl_{n,m}$-linear map $T$ is
\beq
{\rm sTr}_{\Cl_{n,m}}(T)={\rm sTr}(u^{(m-n)/2}\Gamma_{n,m}\circ T)\in \R[u^{1/2},u^{-1/2}],\label{eq:Cliffordtrace}
\eeq
where $u^{1/2}$ is an invertible formal variable of degree $-1$ with inverse $u^{-1/2}$. 
\end{defn}

From properties of the ordinary super trace and the definition of~$\Gamma_{n,m}$ in~\eqref{eq:Gammadef}, we observe that the Clifford super trace satisfies
\beq
{\rm sTr}_{\Cl_{n,m}}(T_1\oplus T_2)= {\rm sTr}_{\Cl_{n,m}}(T_1)+{\rm sTr}_{\Cl_{n,m}}(T_2)\label{eq:Cliffordtraceplus}
\eeq
for $T_i$ an endomorphism of a $\Cl_{n,m}$-module $W_i$ for $i=1,2$, and
\beq
{\rm sTr}_{\Cl_{n,m}\otimes \Cl_{n',m'}}(T_1\otimes T_2)= {\rm sTr}_{\Cl_{n,m}}(T_1)\cdot {\rm sTr}_{\Cl_{n',m'}}(T_2)\label{eq:Cliffordtracetimes}
\eeq
for $T_1$ an endomorphism of a $\Cl_{n,m}$-module $W_1$ and for $T_2$ an endomorphism of a $\Cl_{n',m'}$-module $W_2$. 

There is a completely analogous trace theory for modules over complex Clifford algebras, enjoying the same formal properties as the real Clifford super trace. The choice of normalization in~\eqref{eq:Cliffordtrace}  determines a normalization of the complex Clifford supertrace, using that $\Det(\C^n)\simeq \Det(\R^n)\otimes \C$, i.e., that a real trivialization of the determinant line also provides a trivialization of its complexification.

\begin{ex}\label{eq:protoChern}
For the $\Cl_{1,1}$-module in Example~\ref{ex:Morita1}, 
$$
{\rm sTr}_{\Cl_{1,1}}(\id_{\R^{1|1}})={\rm sTr}(\Gamma_{1,1})=2/2=1.
$$ 
From this we see that the Clifford super trace of a $\Cl_{n,m}$-linear operator $T$ is compatible with the Morita equivalences from Example~\ref{ex:Morita1}
$$
{\rm sTr}_{\Cl_{n,m}\otimes \Cl_{1,1}}(T\otimes \id_{\R^{1|1}})={\rm sTr}_{\Cl_{n,m}}(T)\cdot {\rm sTr}_{\Cl_{1,1}}(\id_{\R^{1|1}})={\rm sTr}_{\Cl_{n,m}}(T). 
$$
For the $\Cl_4$-module in Example~\ref{ex:CL4} we find 
\beq
{\rm sTr}_{\Cl_{4}}(\id_{\R^{4|4}})={\rm sTr}(u^{-4/2}\Gamma_4)=2^{-4/2}\cdot 8u^{-2}=2u^{-2},\label{eq:strCl4}
\eeq
and for the $\Cl_8$-module in Example~\ref{ex:Morita2}, 
\beq
{\rm sTr}_{\Cl_{8}}(\id_{\R^{8|8}})={\rm sTr}(u^{-8/2}\Gamma_8)=u^{-4}2^{-8/2}\cdot 16=u^{-4}.\label{eq:strCl8}
\eeq
We verify that~\eqref{eq:strCl4} and~\eqref{eq:strCl8} are compatible with the isomorphism of $\Cl_8$-modules in~\eqref{eq:KOrelation},
$$
2u^{-2}\cdot 2u^{-2}={\rm sTr}_{\Cl_{4}\otimes \Cl_4}(\id_{\R^{4|4}}\otimes \id_{\R^{4|4}})={\rm sTr}_{\Cl_8}((\id_{\R^{8|8}})^{\oplus 4})=4u^{-4}.
$$
From~\eqref{eq:strCl8} we also find that for a $\Cl_n$-linear operator $T$, 
$$
{\rm sTr}_{\Cl_{n}\otimes \Cl_{8k}}(T\otimes (\id_{\R^{8|8}})^{\otimes k})={\rm sTr}_{\Cl_{n}}(T)\cdot ({\rm sTr}_{\Cl_{8}}(\id_{\R^{8|8}}))^k={\rm sTr}_{\Cl_{n}}(T)\cdot u^{-4k},
$$
so that the $8$-periodicity of Clifford modules is reflected by the Clifford super trace, where the periodicity generator is $u^{\pm 4}\in \R[u^{1/2},u^{-1/2}]$.
\end{ex}

\begin{ex} The complex Clifford module from Example~\ref{ex:cMorita2} has
$$
{\rm sTr}_{\cCl_{2}}(\id_{\C^{1|1}})={\rm sTr}(\Gamma_{2})=2iu^{-1} /2=iu^{-1}.
$$ 
Hence, analogous to the 8-periodicity in the real case
\beq
&&{\rm sTr}_{\cCl_{n}\otimes \cCl_{2k}}(T\otimes (\id_{\C^{1|1}})^{\otimes k})={\rm sTr}_{\cCl_{n}}(T)\cdot ({\rm sTr}_{\cCl_{2}}(\id_{\C^{1|1}}))^k={\rm sTr}_{\cCl_{n}}(T)\cdot i^ku^{-k}.\label{eq:Botttrace}
\eeq
\end{ex}

\begin{rmk}
When working with complex Clifford modules, some authors incorporate a factor of $i^{n/2}$ in the definition~\eqref{eq:Gammadef} of $\Gamma_{n}$ so that ${\rm sTr}_{\cCl_{2}}(\id_{\C^{1|1}})=u^{-1}$, e.g., see~\cite[Proposition~3.21]{BGV} and~\cite[Definition~3.2.16]{ST04}. However, for real Clifford modules this is a bit less natural as the resulting $\Gamma_n$ is not an element of the real Clifford algebra~$\Cl_n$. 
One reason to incorporate the formal variable $u$ in~\eqref{eq:Cliffordtrace} is that any two choices of normalization can be related by rescaling~$u$ by the appropriate scalar. 
\end{rmk}

\section{The Atiyah--Bott--Shapiro construction and the spectrum $\KO$}\label{sec:ABS}

In this section we review a Clifford module construction of the spaces representing $\KO^n(-)$ from~\cite[\S7]{HST} and~\cite[\S1.3-1.4]{Cheung}. Our presentation emphasizes the connection with the Atiyah--Bott--Shapiro construction and the families index.

\subsection{The Atiyah--Bott--Shapiro construction}\label{sec:ABSreview}
Let $M$ be a manifold, not necessarily compact. A class in the compactly supported real K-theory group $\KO^0_c(M)$ can be specified by a pair $(V,Q)$ for a $\Z/2$-graded real vector bundle $V\to M$ and $Q\in \End(V)_\odd$ an odd endomorphism that is invertible outside of a compact subset of~$M$. Roughly,
\beq
[V]=[V_\ev]-[V_\odd]\in \KO^0(M) \label{eq:formaldiff}
\eeq
determines a class, and on the subspace of $M$ where $Q$ is invertible we obtain an isomorphism $V_\ev\stackrel{\sim}{\to} V_\odd$, so that the formal difference~\eqref{eq:formaldiff} vanishes on this subspace; we explain this in greater detail in Remark~\ref{rmk:ABScs}. 

Given a right $\Cl_n$-module $W$, consider the trivial bundle $\underline{W}=W\times \R^n\to \R^n$. Viewing this as a $\Z/2$-graded vector bundle, equip it with the odd operator $\cl_x\in \End(W)_\odd$ that at a point $x\in \R^n$ is Clifford multiplication by $x\in \R^n\subset \Cl_n$ on $W$. This constructs a map
\beq
\Mod_{\Cl_n}\to \KO^0_c(\R^n)\simeq \KO^{-n}(\pt),\qquad W\mapsto (\underline{W},\cl_x) \label{eq:ABS0}
\eeq
where the class has compact support because the operator $\cl_x$ is invertible away from $0\in \R^n$ as $(\cl_x)^2=-|x|^2$ is multiplication by a scalar. 
The isomorphism in~\eqref{eq:ABS0} is the suspension isomorphism, identifying the 1-point compactification of $\R^n$ with the $n$-sphere. Direct sums of $\Cl_n$-modules are sent to addition of classes in $\KO^{-n}(\pt)$ under~\eqref{eq:ABS0}.

Let $\MM_n$ denote the Grothendieck group of $\Mod_{\Cl_n}$. The standard inclusion $\R^n\subset \R^{n+1}$ yields a homomorphism of superalgebras $i\colon \Cl_n\hookrightarrow \Cl_{n+1}$ determined by the map on generators~\eqref{eq:gennotation}. This induces restrictions 
$$
\Mod_{\Cl_{n+1}}\to \Mod_{\Cl_{n}}, \quad i^*\colon \MM_{n+1}\to \MM_n, \qquad n\ge 0.
$$
Atiyah--Bott--Shapiro~\cite[Theorem~11.5]{ABS} show that~\eqref{eq:ABS0} factors through the quotient $\MM_n/i^*\MM_{n+1}$, and in fact determines isomorphisms of abelian groups
\beq
\MM_n/i^*\MM_{n+1}\stackrel{\sim}{\to} \KO^{-n}(\pt)\qquad n\ge 0.\label{eq:oldABS}
\eeq
Furthermore, the tensor product of Clifford modules is compatible with the graded ring structure on $\bigoplus_{n\in \N} \KO^{-n}(\pt)$. Together with the 8-fold periodicity of Clifford algebras and~$\KO$, the maps~\eqref{eq:oldABS} give a completely algebraic description of the graded ring~$\KO^\bullet(\pt)$.

%There is an automorphism of the monoidal category of super vector spaces that permits a version of the isomorphism~\eqref{eq:oldABS}, but with a different sign. This version turns out to be more natural for later geometric applications. Consider the lax monoidal functor $F\colon ({\sf SVect},\otimes) \to ({\sf SVect},\otimes) $ from super vector spaces to itself whose underlying functor is the identity functor and whose lax monoidal structure is 
%\beq
%&&V\otimes W=F(V)\otimes F(W)\to F(V\otimes W)=V\otimes W,\qquad v\otimes w\mapsto (-1)^{|v||w|}v\otimes w. \label{eq:Theo}
%\eeq
%Lax monoidal functors send algebra objects to algebra objects, and we observe that $F(\Cl_n)\simeq \Cl_{-n}$ for all $n\in \Z$. Hence, $F$ together with~\eqref{eq:oldABS} gives isomorphisms 
%\beq
%\MM_{-n}/F(i)^*\MM_{-n-1} \simeq \KO^{-n}(\pt)\qquad n\ge 0. \label{eq:ABS}
%\eeq

\begin{rmk}\label{rmk:signs}
We compare sign conventions in the Atiyah--Bott--Shapiro isomorphism~\eqref{eq:oldABS}. First, we recall the equivalence of categories $\Mod_{\Cl_n}\simeq {}_{\Cl_{-n}}\Mod$ using the canonical isomorphism $\Cl_{\pm n}^\op=\Cl_{\mp n}$. This allows one to rephrase~\eqref{eq:oldABS} in terms of left Clifford modules, 
\beq
{}_{-n}\MM/i^*{}_{-n-1}\MM \simeq \KO^{-n}(\pt)\qquad n\ge 0, \quad i\colon \Cl_{-n}\hookrightarrow \Cl_{-n-1} \label{eq:ABS}
\eeq
but with the indicated sign change. One can instead use the ungraded anti-automorphism of~$\Cl_n$ (sometimes called the \emph{transpose}) that gives a left-module version of~\eqref{eq:oldABS} but without a sign change. This second option is a bit awkward when considering tensor products of modules and the corresponding products in $\KO^\bullet(\pt)$: the ungraded tensor product of $\Cl_m$ and $\Cl_n$ is not a Clifford algebra. 

Our choice of sign convention below is fixed by our preference for left modules and the example of the Clifford linear Dirac operator on an $n$-dimensional Riemannian spin manifold. 
If one demands that the Dirac Laplacian is a nonnegative operator, then it must be defined on the associated bundle for the $\Spin(n)$ action on $\Cl_{n}$ as in~\eqref{eq:Riemspinbundle} (the associated bundle for~$\Cl_{-n}$ leads to a nonpositive Dirac Laplacian). This convention agrees with~\cite[II.7]{LM}. Using the canonical isomorphism $\Cl_n^\op=\Cl_{-n}$, the kernel of the Clifford linear Dirac operator is then a left $\Cl_{-n}$-module, which determines an element of $\KO^{-n}(\pt)$ using the convention~\eqref{eq:ABS}. 
\end{rmk}

\begin{rmk} The Atiyah--Bott--Shapiro isomorphism~\eqref{eq:ABS} can be verified completely explicitly using the generators and relations description~\cite[I.9]{LM}
\beq\label{eq:gensrelation}
&&\pi_*\KO=\Z[\eta,\alpha,\beta^{\pm1}]/\langle 2\eta,\eta^3,\eta\alpha,\alpha^2-4\beta\rangle,\quad |\eta|=-1, \ |\alpha|=-4, \ |\beta|=-8.
\eeq
We sketch surjectivity of~\eqref{eq:ABS}. The Clifford modules constructed in the previous section correspond to classes
$$
[_{\Cl_{-1}}\R^{1|1}]=\eta\in \KO^{-1}(\pt),\quad [_{\Cl_4}\R^{4|4}]=\beta^{-1}\alpha\in \KO^4(\pt),\quad [_{\Cl_8}\R^{8|8}]=\beta^{-1}\in \KO^8(\pt)
$$
where the $\Cl_{-1}$-action on $\R^{1|1}$ is defined in Example~\ref{ex:Morita1} (by forgetting the $\Cl_1$-action), the $\Cl_4$-action on $\R^{4|4}$ in Example~\ref{ex:CL4} and the $\Cl_8$-action on $\R^{8|8}$ in Example~\ref{ex:Morita2}. To see that $2\eta=\eta^3=0$, one must show that $\R^{1|1}\oplus \R^{1|1}\simeq \R^{2|2}$ extends to a $\Cl_{-1}\otimes \Cl_{-1}\simeq \Cl_{-2}$-module and the $\Cl_{-3}$-module $(\R^{1|1})^{\otimes 3}\simeq \R^{4|4}$ extends to a $\Cl_{-4}$-module; for the former, the second Clifford generator can be taken to be the isomorphism that exchanges the copies of $\R^{1|1}$ in the direct sum, and for the latter one uses a module structure closely related to~\eqref{eq:CL4}. To see $\eta\alpha=0$, it suffices to check $\eta\alpha\beta^{-1}=0$; starting with the $\Cl_{-1}\otimes \Cl_4$-module $\R^{1|1}\otimes \R^{4|4}$, a Morita equivalence provides a $\Cl_3$-module with underlying vector space $\R^{4|4}$, and this similarly has an extension to a $\Cl_4$-module. Invertibility of $\beta$ follows from the Morita equivalence implemented by~\eqref{eq:Morita2}, and the relation involving $\alpha$ and $\beta$ is equivalent to the isomorphism~\eqref{eq:KOrelation}. 
\end{rmk}

 \subsection{Unbounded operators that represent $\KO$}\label{sec:unbounded}

Next we review a space of unbounded operators that is known to represent $\KO$, following the presentation in~\cite[\S3]{HST}. Fix a real, infinite-dimensional, separable Hilbert space $\mathcal{H}$. An \emph{unbounded operator}~$D$ on~$\mathcal{H}$ is the data of a subspace ${\rm dom}(D)\subset \mathcal{H}$ (the domain of $D$) and a linear map $D\colon {\rm dom}(D)\to \mathcal{H}$. An unbounded operator is \emph{self-adjoint} if it is self-adjoint on the closure of its domain. An unbounded, self-adjoint operator $D$ has \emph{compact resolvent} if the spectrum of $D$ consists of eigenvalues of finite multiplicity that do not have an accumulation point in~$\R$. 

\begin{defn} 
Let $\Inf(\mathcal{H})$ denote the set of unbounded, self-adjoint operators on~$\mathcal{H}$ with compact resolvent. 
\end{defn} 

\begin{rmk}
The notation $\Inf$ comes from the terminology in~\cite{HST}: $D\in \Inf(\mathcal{H})$ is an \emph{infinitesimal generator} of a super semigroup representation. Although we do not use super semigroup representations, we keep the same notation for ease of comparison with~\cite{HST}. We caution that our sign conventions differ slightly from \cite{HST}; see Remark~\ref{eq:HSTsigns}. 
\end{rmk}

To simplify the exposition, hereafter we will typically omit the words ``compact resolvent" and ``self-adjoint," referring to an element of $\Inf(\mathcal{H})$ simply as an unbounded operator on~$\mathcal{H}$. The Cayley transform determines an injective map~\cite[Proposition 3.7]{HST} 
\beq
\Inf(\mathcal{H})\hookrightarrow O(\mathcal{H})\label{eq:orthogonalsubspace}
\eeq
from unbounded operators to the infinite orthogonal group, sending eigenspaces of an operator $D$ with (real) eigenvalue $\lambda$ to an orthogonal operator with (complex) eigenvalue $\frac{\lambda+i}{\lambda-i}$. 
%We recall that the operator norm topology on $\mathcal{U}(\mathcal{H})$ is determined by the norm
%$$
%\|u\|=\sup \{\|u(x)\| : \|x\|=1, \ x\in \mathcal{H}\}
%$$
%on bounded operators $u\in \mathcal{B}(\mathcal{H})$. 
Equip the set of unbounded operators with the subspace topology inherited from~$O(\mathcal{H})$, where $O(\mathcal{H})$ is given the norm topology. 

Next we consider variations on the space $\Inf(\mathcal{H})$ of unbounded operators that involve gradings and Clifford actions. Define a \emph{$\Z/2$-graded Hilbert space} as a separable Hilbert space whose underlying vector space has a $\Z/2$-grading such that the even and odd subspaces are infinite-dimensional and orthogonal to each other. A \emph{stable $\Cl_n$-module} is a $\Z/2$-graded Hilbert space with a self-adjoint left $\Cl_{n}$-module structure such that  every isomorphism class of irreducible $\Cl_{n}$-module appears with infinite multiplicity; for example, starting with a $\Z/2$-graded Hilbert space $\mathcal{H}$, then $\mathcal{H}_n:=\Cl_n\otimes \mathcal{H}$ is a stable $\Cl_n$-module. Any pair of stable $\Cl_n$-modules are isomorphic, and the tensor product of a stable $\Cl_n$-module with a stable $\Cl_m$-module determines a stable $\Cl_{n+m}$-module (using the Morita equivalence from Example~\ref{ex:Morita1} when $n$ and $m$ have different signs). We refer to~\cite[\S3]{Joachim} for further discussion of this tensor product of stable $\Cl_n$-modules. 

A self-adjoint unbounded operator $D$ on $\mathcal{H}_n$ is \emph{odd} if its domain is a graded subspace and $D\colon {\rm dom}(D)\to \mathcal{H}$ is an odd linear map. An odd self-adjoint unbounded operator $D$ on~$\mathcal{H}_n$ is \emph{Clifford linear} if ${\rm dom}(D)$ is a Clifford submodule and $D$ graded commutes with the $\Cl_{n}$-action. Let $\Inf_{-n}$ denote the set of odd, Clifford linear unbounded operators on~$\mathcal{H}_n$. We equip this set with the subspace topology for the evident inclusion
\beq
\Inf_{-n}\subset \Inf(\mathcal{H}_n).\label{eq:subspaceagain}
\eeq

\begin{rmk}\label{eq:HSTsigns}
In \cite{HST}, the authors use stable \emph{right} $\Cl_n$-modules rather than the stable left $\Cl_n$-modules above, compare Remark~\ref{rmk:signs}. The spaces $\Inf_{-n}$ in the present paper are homeomorphic to the ones in \cite{HST} (with the sign unchanged) by trading a left $\Cl_{n}$-action on $\mathcal{H}_n=\Cl_n\otimes \mathcal{H}$ for a right $\Cl_{-n}$-action. 
\end{rmk}

It is important in applications that the domain of an unbounded operator $D\in \Inf_{-n}$ is \emph{not} required to be dense in~$\mathcal{H}_n$; this permits the following example. 

\begin{ex} \label{ex:finiterank}
Consider a graded, $\Cl_{n}$-invariant subspace $W\hookrightarrow \mathcal{H}_n$. Define the unbounded operator $P_W$ whose domain is $W$ and $P_W\colon W\to \mathcal{H}_n$ is the inclusion. If $W$ is finite-dimensional, then $P_W$  has compact resolvent and~$P_W\in \Inf_{-n}$. 
\end{ex}

Perhaps the most well-known analytic model for the $\KO$-spectrum comes from the Atiyah--Singer spaces $\Fred_n$ of $\Cl_n$-linear Fredholm operators that for $n\ge 0$ represent~$\KO^{-n}$ \cite{AtiyahSingerskew}, \cite[III.10]{LM}. The spaces $\Inf_{n}$ give an unbounded version of these spaces in the following sense (we continue to use the sign convention from Remark~\ref{rmk:signs}). 

\begin{thm}[\cite{HST} Proposition~3.11 and Theorem~8.2]\label{thm:HST}
For $n\ge 0$, there are homotopy equivalences $\Inf_{n}\stackrel{\sim}{\to} \Fred_{n}$. Hence, $\Inf_n$ represents the functor $\KO^{-n}(-)$.
\end{thm}

The 8-periodicity of Clifford modules gives homeomorphisms $\Inf_{n}\simeq \Inf_{n-8}$, and so (by Bott periodicity) one concludes that $\Inf_{n}$ represents $\KO^{-n}(-)$ for all $n\in \Z$. Alternatively, one can extend the construction of the spaces $\Fred_n$ to arbitrary $n\in \Z$, and prove a version of Theorem~\ref{thm:HST} for all~$n$, e.g., see \cite[Remark 0.1]{HST}.

\subsection{Topological categories and their classifying spaces}

A \emph{topological category} is a category in which the objects and morphisms are topological spaces for which the structure maps (source, target, unit, composition) are continuous. For a topological category $\mathcal{C}$, let~$\mathcal{C}_0$ denote the space of objects and $\mathcal{C}_1$ the space of morphisms. 

Given an open cover of a topological space, the following two examples construct topological categories. 

\begin{ex} \label{ex:opencover0}
Let $Z$ be a topological space, and $\{U_i\}_{i\in I}$ an open cover of $Z$. For $\sigma\subset I$ a nonempty finite subset, let $U_\sigma=\bigcap_{i\in \sigma} U_i$ denote the intersection. Define the \emph{\v{C}ech category} of the cover, denoted $\check{{\rm C}}(U_i)$, as the topological category whose objects are $\check{{\rm C}}(U_i)_0=\coprod_{\sigma} U_\sigma$ and whose morphisms are the space $\check{{\rm C}}(U_i)_1=\coprod_{\sigma\subseteq \tau} U_\tau$. The source and target maps are determined by the projection and inclusion $U_\tau\subseteq U_\sigma$, the unit is induced by the identity map $U_\sigma \to U_\sigma$, and composition comes from nested inclusions.
\end{ex}

\begin{ex} \label{ex:opencover00}
Let $Z$ be a topological space, and $\{U_i\}_{i\in I}$ an open cover of $Z$ with the data of an ordering on the indexing set $I$. Then define a topological category $\check{{\rm C}}^\ord(U_i)$ with objects and morphisms
$$
\check{{\rm C}}^\ord(U_i)_0=\coprod_{i\in I} U_i,\qquad \check{{\rm C}}^\ord(U_i)_1=\coprod_{i\le j} U_{ij},\qquad U_{ij}=U_i\bigcap U_j.
$$
The source and target maps are the inclusions $U_i\bigcap U_j\subset U_i,U_j$, the unit map is the identity on~$U_i=U_i\bigcap U_i$, and composition comes from inclusions of triple intersections into double intersections. 
\end{ex}

\begin{defn} The \emph{nerve} of a topological category $\mathcal{C}$ is the simplicial space whose $k$-simplices are $N_k\mathcal{C}=\mathcal{C}_k$ for $k=0,1$, and for $k>1$ are length $k$ chains of composable morphisms
\beq
N_k\mathcal{C}&=&\{x_0\stackrel{f_1}{\to}x_1\stackrel{f_2}{\to}x_2\stackrel{f_3}{\to}\cdots \stackrel{f_k}{\to}x_k\mid x_i\in \mathcal{C}_0, f_i\in \mathcal{C}_1\}\nonumber\\&=&\mathcal{C}_1\times_{\mathcal{C}_0}\cdots \times_{\mathcal{C}_0} \mathcal{C}_1,\nonumber 
\eeq
topologized as the fibered product.
The face maps in  $N_\bullet \mathcal{C}$ are determined by composing morphisms, and degeneracies are determined by the unit map in $\mathcal{C}$. 
\end{defn}

The \emph{geometric realization} of a simplicial space $Z_\bullet$ is 
\beq
| Z_\bullet | :=(\coprod_n Z_n\times \Delta^n)/{\sim}\qquad (f^*x,t)\sim (x,f_*t) \quad \forall f\colon [l]\to [k]\label{eq:fatrealizae}
\eeq
where $\Delta^n$ is the standard $n$-simplex, and $f\colon [l]\to [k]$ is an order-preserving map. Geometric realization is natural: a map of simplicial spaces $F\colon Z_\bullet\to Y_\bullet$ determines a map between realizations, $|F|\colon |Z_\bullet|\to |Y_\bullet|$. 

\begin{defn} Given a topological category $\mathcal{C}$, define the \emph{classifying space of $\mathcal{C}$} as the geometric realization of its nerve, $
\Bc\mathcal{C}:=| N_\bullet \mathcal{C}|.$ 
\end{defn}

\begin{ex}\label{ex:opencover}
The nerve of the topological category $\check{C}^\ord(U_i)$ from Example~\ref{ex:opencover00} has as $k$-simplices the intersections
\beq\label{eq:intersectionnotion}
&&N_k\check{C}(U_i)=\coprod_{i_0\le i_1\le \cdots \le i_k} U_{i_0 i_1\cdots i_k},\qquad U_{i_0 i_1\cdots i_k}=U_{i_0}\bigcap U_{i_1}\bigcap \cdots \bigcap U_{i_k}.
\eeq
whereas the nerve of the topological category $\check{C}(U_i)$ from Example~\ref{ex:opencover0} has as $k$-simplices 
$$
N_k\check{C}(U_i)=\coprod_{\sigma_0\subseteq \sigma_1\subseteq \cdots \subseteq \sigma_k} U_{\sigma_k}. 
$$
Results of Segal~\cite[Proposition~4.12]{Segalclassifying}, and Dugger--Isaksen~\cite[Theorem~2.1 and Proposition~2.6]{DuggerIsaksen} show that the canonical maps
\beq
&&\Bc \check{C}^\ord(U_i)\stackrel{\sim}{\to} \Bc \check{C}(U_i),\qquad  \Bc \check{C}^\ord(U_i)\stackrel{\sim}{\to} Z, \qquad \Bc \check{C}(U_i)\stackrel{\sim}{\to} Z,\label{eq:Segalwe}
\eeq 
are homotopy equivalences.
 \end{ex}

There is a homotopy equivalence $\Bc\mathcal{C}^\op\simeq \Bc\mathcal{C}$ between the realization of a category and the realization of its opposite. Below we will often identify $\Bc \check{C}(U_i)$ with $\Bc \check{C}(U_i)^\op$.

\subsection{The Atiyah--Bott--Shapiro category} \label{sec:ABScat}

\begin{defn} 
Let $\M_n^\delta$ be the category whose objects are finite-dimensional $\Cl_{n}$-submodules $W\subset \mathcal{H}_n$ and whose morphisms are
\beq
&&{\rm Mor}(W_0,W_1):=\left\{ \begin{array}{cl} \begin{array}{c} \Cl_{n}\otimes \Cl_{-1}-{\rm action\ on\ } W_0^\perp \\ {\rm extending\ the }\ \Cl_{n}{\rm-action}\end{array}  & {\rm if} \ W_0\subseteq W_1\subset \mathcal{H}_n \\
\null\\
\emptyset & {\rm else} \end{array}\right.
\eeq
where $W_0^\perp\subseteq W_1$ denotes the orthogonal complement of the subspace $W_0\subseteq W_1$. Composition for a nested inclusion $W_0\subseteq W_1\subseteq W_2$ uses the direct sum of Clifford modules to obtain a $\Cl_{n}\otimes \Cl_{-1}$-module structure on the orthogonal complement of $W_0\subseteq W_2$. 
\end{defn} 

\begin{rmk}\label{rmk:esquared}
An extension to a $\Cl_{n}\otimes \Cl_{-1}$-action is the data of a self-adjoint odd endomorphism $e$ on $W_0^\perp$ that graded commutes with the $\Cl_n$-action and satisfies $e^2=1$. 
\end{rmk}
\begin{rmk}
For $n\ge 0$, a $\Cl_{-n}\otimes \Cl_{-1}$-module is the same data as a $\Cl_{-n-1}$-module. Hence, for $n\le 0$,~$\M_n^\delta$ is a category whose morphisms parameterize how the Atiyah--Bott--Shapiro map~\eqref{eq:ABS} identifies (left) Clifford modules in $\KO^{n}(\pt)$.
\end{rmk}

Let $\mathcal{B}(\mathcal{H}_n)$ denote the space of bounded $\Cl_{n}$-linear operators on $\mathcal{H}_n$. There are injections
\beq
{\rm Obj}(\M_n^\delta)\hookrightarrow \mathcal{B}(\mathcal{H}_n),\qquad {\rm Mor}(\M_n^\delta)\hookrightarrow \mathcal{B}(\mathcal{H}_n)^{\times 3}. \label{eq:bddops}
\eeq
The first map comes from identifying a submodule $W_0\subset \mathcal{H}_n$ with a finite rank projection operator. Similarly, $(W_0,W_1,e)\in {\rm Mor}(\M_n^\delta)$ determines a pair of finite rank projection operators and a bounded operator gotten by extending $e$ (in the notation of Remark~\ref{rmk:esquared}) to an operator on $\mathcal{H}_n$ via the zero operator on the orthogonal complement of $W_0^\perp\subset \mathcal{H}_n$. 

\begin{defn} Define the \emph{$n$th Atiyah--Bott--Shapiro category $\M_n$} as the topological category with underlying category $\M_n^\delta$ whose objects and morphisms are endowed with the subspace topology for the inclusions~\eqref{eq:bddops} using the norm topology on bounded operators. 
\end{defn}

\begin{rmk}\label{rmk:productsum}
Since a tensor product of $\Cl_n$ and $\Cl_m$-modules produces a $\Cl_{n+m}$-module and any pair of stable $\Cl_n$-modules are isomorphic,  
there are continuous functors
\beq
\M_n\times \M_m\stackrel{\otimes}{\to} \M_{n+m}, \label{eq:ABStensor}
\eeq
unique up to a contractible space of choices. The topological category $\M_n$ also has a partially-defined sum operation, defined for orthogonal subspaces $W,W'\subset \mathcal{H}_n$. 
\end{rmk}

\subsection{The families index and the Atiyah--Bott--Shapiro category} \label{subsec:ABSindex}

We will compare the classifying space of the $n$th Atiyah--Bott--Shapiro category with the space $\Inf_{-n}$ through an index map, similar to the families index described in~\S\ref{sec:motivate}. Define an open cover with components indexed by $\lambda\in \R_{>0}$
\beq
\Inf_{n,\lambda}:=\{D\in \Inf_{-n} \mid \lambda\notin {\rm Spec}(D^2)\}\subset \Inf_{-n}.\label{eq:infopencov}
\eeq
Let $\check{{\rm C}}(\Inf_{n,\lambda})$ denote the \v{C}ech category of this cover of $\Inf_{-n}$. Explicitly, the objects of this category are given by an unbounded operator $D$ and a $(k+1)$-tuple of positive real numbers,
$$
\check{{\rm C}}(\Inf_{n,\lambda})_0=\{D,\lambda_0\le \lambda_1\le \cdots \le \lambda_k \mid D\in \Inf_{-n}, \ \lambda_i\in \R_{>0}, \ \lambda_i \notin {\rm Spec}(D^2)\},
$$
satisfying the indicated property. Morphisms are given by the same data together with a specified subset $\{\lambda_{i_0},\dots,\lambda_{i_l}\}\subseteq \{\lambda_0,\dots,\lambda_k\}$ with~$\lambda_{i_0}\le \lambda_{i_1}\le\cdots \le\lambda_{i_l}$.

Define a continuous functor
\beq
\Ind \colon \check{{\rm C}}(\Inf_{n,\lambda})^\op\to \M_n\label{eq:contfunctor}
\eeq
as follows. 
On objects, $\Ind$ assigns to $(D,\lambda_0\le \lambda_1\le \cdots \le \lambda_k)$ the subspace $\mathcal{H}^{<\lambda_k}\subset\mathcal{H}_n$ given by the direct sum of eigenspaces of $D^2$ of eigenvalue less than $\lambda_k$,
$$
\mathcal{H}^{<\lambda_k}:=\bigoplus_{\mu<\lambda} \mathcal{H}_\mu. 
$$
The subspace $\mathcal{H}^{<\lambda_k}\subset\mathcal{H}_n$ is finite-dimensional because~$D$ has compact resolvent. Furthermore, $\mathcal{H}^{<\lambda_k}$ is a graded $\Cl_{n}$-submodule because~$D$ is a Clifford linear operator. On morphisms, the functor~\eqref{eq:contfunctor} assigns the inclusion of Clifford submodules $\mathcal{H}^{<\lambda_{i_l}}\subset \mathcal{H}^{<\lambda_k}\subset\mathcal{H}_n$, where on the orthogonal complement to the inclusion the operator $D$ is used to construct a $\Cl_{-1}$-action as follows. First decompose the orthogonal complement to the inclusion into $\mu$-eigenspaces for the operator~$D^2$,
$$
(\mathcal{H}^{<\lambda_{i_l}})^\perp\simeq \bigoplus_{\lambda_{i_l}<\mu<\lambda_{k}} \mathcal{H}_{\mu}.
$$ 
The sum is finite because $D^2$ has compact resolvent and is nonnegative, and the eigenspaces are orthogonal because $D^2$ is self-adjoint. Now define an odd operator 
$$
e\colon (\mathcal{H}^{<\lambda_{i_l}})^\perp\to (\mathcal{H}^{<\lambda_{i_l}})^\perp,\qquad  e:=\bigoplus_\mu (\frac{1}{\sqrt{\mu}}D|_{\mathcal{H}_{\mu}}). 
$$ 
By construction, $e^2=+1$, and hence determines a $\Cl_{-1}$-action (see Remark~\ref{rmk:esquared}). This action is self-adjoint and graded commutes with the $\Cl_n$-action because $D$ is self-adjoint and $\Cl_n$-linear. 

\begin{lem}\label{lem:Kthy} The continuous functor~\eqref{eq:contfunctor} induces a weak equivalence 
\beq
\Bc\Ind\colon \Bc\check{{\rm C}}(\Inf_{n,\lambda})\stackrel{\sim}{\to} \Bc\M_n,\label{eq:index}
\eeq
and hence $\Bc\M_n$ represents the functor $\KO^{n}(-)$. 
\end{lem}

\bp
In~\cite[page 53]{HST}, there is a topological category denoted~${\sf D}_n$ with the same objects as $\M_n$ and whose morphisms are denoted
$$
{\rm Mor}_{{\sf D}_n}(W_0,W_1)=\left\{\begin{array}{cl} \{ R\in O_{C_{n}}(W_2^\perp) \mid R \ {\rm odd} \ R^2=1\} & {\rm if } \ W_0\subset W_1 \\ \emptyset & {\rm else.} \end{array}\right.
$$
Since $R$ is orthogonal and $R^2=1$, it generates a self-adjoint action of $\Cl_{-1}$. From this we observe that ${\sf D}_n$ and~$\M_n$ are isomorphic as topological categories: their spaces of objects and morphisms are homeomorphic. Then~\cite[Proposition~3.11, Proposition~4.6, and Theorem~7.1]{HST} construct a homotopy equivalence $\Bc\M_n\stackrel{\sim}{\to} \Inf_{-n}$ that sits in the homotopy commutative triangle,
\beq
\begin{tikzpicture}[baseline=(basepoint)];
\node (A) at (0,0) {$\Bc\check{C}(\Inf_{n,\lambda})$};
\node (B) at (-1.5,-1.5) {$\Bc\M_n$};
\node (C) at (1.5,-1.5) {$\Inf_{-n}$};
\draw[->] (A) to node [left=3pt] {$\Ind$} (B);
\draw[->] (A) to node [right=3pt] {$\sim$} (C);
\draw[->] (B) to node [below] {$\sim$} (C);
\path (0,-.6) coordinate (basepoint);
\end{tikzpicture} \nonumber
\eeq
where the weak equivalence $\Bc\check{C}(\Inf_{n,\lambda})\to \Inf_{-n}$ follows from~\eqref{eq:Segalwe}, and homotopy commutativity follows from the argument in~\cite[Proposition~4.6]{HST}. By the 2-out-of-3 property, this shows that the map~\eqref{eq:index} is a weak equivalence. The second statement in the lemma then follows from Theorem~\ref{thm:HST}. 
%Homotopy commutativity: The homotopy sends all eigenspaces with eigenvalue $>\lambda$ to infinity. First use the homeomorphism from~\cite[Proposition 27]{HST}. Then apply a version of the homotopy in~\cite[Proposition 29]{HST},
%$$
%h_t(x)=\left\{\begin{array}{ll} \frac{x}{1-(t/\lambda)|x|} & x\in (-\lambda/t,\lambda/t)\\ \infty & {\rm else} \end{array}\right.
%$$ 
%(the one there is for $\lambda=1$). Or more easily:
%$$
%h_t(x)=\left\{\begin{array}{ll} x & x\in (-\lambda,\lambda)\\ x/(1-t) & {\rm else} \end{array}\right.
%$$ 
%Homotopy commutativity in more detail: Need to consider the map on the realization, so have to parameterize over a simplex as well. Have a pair of maps 
%$$
%(U_{\lambda_0}\bigcap \cdots \bigcap U_{\lambda_k}\times \Delta^k)\to \Inf_{-n}
%$$
%where the first sends $(D,\lambda_0,\dots,\lambda_k,t_1,\dots,t_k)$ to $D$ and the second maps to an unbounded operator whose domain is contained in $\mathcal{H}_{\lambda_k}$. There is a homotopy between these compositions that can be described in terms of configurations spaces: $\mathcal{H}_{\lambda_k}^\perp$ gets sent to the infinite eigenvalue, and then the homotopy on $\mathcal{H}^{<\lambda}$ depends on $(t_1,\dots,t_k)$. 
%The infinite eigenvalue homotopy is basically the same as in~\cite[Proposition 29.]{HST}, using that on $U_\lambda$ we can push all eigenspaces with eigenvalue larger the $\lambda$ away uniformly. The remaining homotopy is between maps on a finite dimensional vector space. We can use the linear homotopy here between any pair of endomorphisms. 
\ep

\begin{rmk}\label{rmk:HSTCheung}
The fact that $\Bc\M_n$ represents $\KO^n(-)$ follows directly from~\cite[Theorem~7.1]{HST}. An independent proof is given in \cite[\S1.3-1.4]{Cheung}, using a comparison between~$\Bc\M_n$ and Clifford linear Fredholm operators. The only new content in Lemma~\ref{lem:Kthy} is the observation that the index map implements this weak equivalence. 
\end{rmk} 

\begin{rmk}\label{rmk:ABScs}
We revisit the Atiyah--Bott--Shapiro description of compactly supported $\KO$-class $[V,Q]\in \KO_c^0(M)$ from~\eqref{eq:formaldiff}. Choose an open cover $U_0,U_1$ of $M$ with the property that the closure of $U_0$ is a compact subset of $M$ and $Q|_{U_1}$ is invertible. Next choose an embedding $V|_{U_0}\subset U_0\times\mathcal{H}_0$ as a summand of the trivial Hilbert bundle; such an embedding is unique up to a contractible space of choices by~\cite{Kuiper}. Next consider the constant map $U_1\to \M_0$ to the basepoint of $\M_0$, i.e., the zero module $\{0\}\subset \mathcal{H}_n$. On the overlap $U_0\bigcap U_1$, the odd operator~$Q$ determines a continuous family of morphisms in $\M_0$ between $V|_{U_0}\subset U_0\times\mathcal{H}_0$ and the zero object. All together, we obtain a continuous functor $(V,Q)\colon \Bc \check{C}(U_i)\to \M_0$ of topological categories. After realization, this determines a map 
$$
M\simeq \Bc \check{C}(U_i)\to \Bc \M_0
$$
representing the compactly supported class $[V,Q]\in \KO_c^0(M)$; see Proposition~\ref{prop:ClnKO} and Corollary~\ref{cor:cpt}. 
\end{rmk} 

\section{$\KO$-classes from $\Cl_n$-bundles}\label{sec:KO}

In this section we define categories of $\Cl_n$-bundles over a smooth manifold $M$ that provides a cocycle model for $\KO$-theory and recovers the Atiyah--Bott--Shapiro category when $M=\pt$. 
All super vector bundles are assumed to be finite-rank. 
A \emph{metric} on a super vector bundle is a metric on its underlying (ungraded) vector bundle such that the even and odd subspaces in each fiber are orthogonal. 

\subsection{$\Cl_n$-bundles}

We begin with a naive candidate for the category of $\Cl_n$-bundles. 

\begin{defn}\label{defn:naive1}
A \emph{basic $\Cl_{n}$-bundle} over $M$ is a metrized super vector bundle $V\to M$ with $\Cl_n$-action for which the fibers are self-adjoint left $\Cl_n$-modules. A \emph{stable isomorphism} between basic $\Cl_n$-bundles is a $\Cl_n$-equivariant inclusion $V\hookrightarrow V'$ over $M$ together with the data of a self-adjoint $\Cl_n\otimes \Cl_{-1}$-action on $V^\perp\subset V'$ extending the given $\Cl_n$-action. 
\end{defn}

Definition~\ref{defn:naive1} recovers a category equivalent to the (discrete) Atiyah--Bott--Shapiro category when~$M=\pt$. Furthermore, a $\Cl_n$-bundle in the above sense determines a map $M\to \Bc\M_n$; see the proof of Proposition~\ref{prop:ClnKO}. However, there are continuous maps $M\to \Bc\M_n$ that are not realized by a basic $\Cl_n$-bundle, and hence $\KO$-classes that do not admit a cocycle description using Definition~\ref{defn:naive1}. The easiest such example is a map $S^1\to \Bc\M_1$ whose class in $\widetilde{\KO}{}^1(S^1)\simeq \Z$ is nontrivial. The following definition repairs this defect. 
We use the notation from Example~\ref{ex:opencover0} for open subsets of a manifold~$M$ and their intersections. 

\begin{defn} \label{defn:Clnbundle}
A \emph{$\Cl_n$-bundle} over a manifold $M$ is data:
\begin{enumerate}
\item an open cover $\{U_i\}_{i\in I}$ of $M$;
\item for each nonempty finite subset~$\sigma\subset I$, a metrized super vector bundle $V_\sigma\to U_\sigma$ with a $\Cl_{n}$-action on $V_\sigma$ for which the fibers are self-adjoint left $\Cl_n$-modules;
\item for each inclusion of nonempty finite subsets $\sigma\subset \tau\subset I$, a $\Cl_n$-equivariant, injective map of metrized super vector bundles $g_{\sigma\tau}\colon V_\sigma|_{U_\tau} \hookrightarrow V_\tau$ over $U_\tau$, together with a $\Cl_{n}\otimes \Cl_{-1}$-action on the orthogonal complement of the image $V_\sigma^\perp\subset V_\tau$ that extends the $\Cl_n$-action.
\end{enumerate}
For $\sigma\subset \tau\subset \rho $, we require an equality of maps $g_{\sigma\tau}\circ g_{\tau\rho}=g_{\sigma\rho}\colon V_\sigma|_{U_\rho}\hookrightarrow V_\rho$, and we require the $\Cl_{-1}$-action on $V_\sigma^\perp\subset V_\tau \hookrightarrow V_\rho$ restricts from the $\Cl_{-1}$-action on~$V_\tau^\perp\subset V_\rho$. We use the notation $e_{\sigma\tau}\in \Gamma(U_{\tau},\End(V_\tau^\perp)_\odd)$ for the generator of the $\Cl_{-1}$-action on $V_\sigma^\perp\subset V_\tau$ and $V=\{V_\sigma,g_{\sigma\tau},e_{\sigma\tau}\}$ when referring to a $\Cl_n$-bundle over $M$. 
\end{defn}

\begin{ex}
Real vector bundles are examples of $\Cl_0$-modules. However, a $\Cl_0$-bundle need not be a vector bundle. 
For example, for a $\Cl_0$-bundle defined relative to an open cover $U_0,U_1\subset M$, Definition~\ref{defn:Clnbundle} requires data
\beq\label{eq:easyexample}
V_0\to U_0, \ \  V_1\to U_1, \qquad V_0|_{U_{01}}\hookrightarrow V_{01}\hookleftarrow V_1 |_{U_{01}}
\eeq
for real vector bundles $V_0$ and $V_1$, and a stable equivalences over $U_{01}=U_0\bigcap U_1$. Vector bundles $V_0$ and $V_1$ that are stably isomorphic over $U_{01}$ but have different ranks therefore determine a $\Cl_0$-bundle that is not a vector bundle.
\end{ex}

\begin{defn}
For a smooth map $\phi \colon N\to M$ and a $\Cl_n$-bundle $V=\{V_\sigma,g_{\sigma\tau},e_{\sigma\tau}\}$, the \emph{pullback} $\phi^*V:=\{\phi^*V_\sigma,\phi^*g_{\sigma\tau},\phi^*e_{\sigma\tau}\}$ is a $\Cl_n$-bundle on $N$ defined relative to the open cover $\{\phi^{-1}(U_i)\}$. In particular, $\Cl_n$-bundles can be restricted to submanifolds $N\subset M$. 
\end{defn}

We require the following versions of equivalence between $\Cl_n$-modules. Each of these generalizes a standard notion of equivalence between vector bundles; see Example~\ref{ex:stable}. 

\begin{defn}\label{defn:iso} An \emph{isomorphism} between $\Cl_n$-bundles $V=\{V_\sigma,g_{\sigma\tau},e_{\sigma\tau}\}$ and $V'=\{V'_{\alpha},g_{\alpha\beta}',e_{\alpha\beta}'\}$ defined relative to open covers $\{U_i\}_{i\in I}$ and $\{U_i'\}_{i\in I'}$, respectively, consists of $\Cl_n$-equivariant metrized vector bundles isomorphisms $\varphi_{\alpha\sigma}\colon V_\sigma|_{U_{\alpha}'}\stackrel{\sim}{\to} V_\alpha|_{U_\sigma}$ for $\sigma\subset I$ and $\alpha\subset I'$ with $U_\alpha'\bigcap U_\sigma$ nonempty. The isomorphisms $\varphi_{\alpha\sigma}$ are required to be compatible with the restrictions of $g_{\sigma\tau},g'_{\alpha\beta},e_{\sigma\tau},e_{\alpha\beta}'$: for $\sigma\subset \tau$ and $\alpha\subset \beta$, the above data are required to fit into a commuting diagram 
\beq
\begin{tikzpicture}[baseline=(basepoint)];
\node (A) at (0,0) {$V_\tau|_{U_\beta'}$};
\node (B) at (2,0) {$V_\beta'|_{U_\tau}$};
\node (C) at (0,-1.5) {$V_\sigma|_{U_\alpha'}$};
\node (D) at (2,-1.5) {$V_\alpha'|_{U_\sigma}$};
\draw[->] (A) to node [above] {$\varphi_{\beta\tau}$} (B);
\draw[->,left hook-latex] (C) to node [left] {$g_{\sigma\tau}$} (A);
\draw[->,left hook-latex] (D) to node [right] {$g_{\alpha\beta}'$} (B);
\draw[->] (C) to node [below] {$\varphi_{\alpha\sigma}$} (D);
\path (0,-.75) coordinate (basepoint);
\end{tikzpicture}\label{eq:compatdiagr}
\eeq
and the $\Cl_{-1}$-actions on orthogonal complements of the images are required to be compatible. We write $\varphi\colon V \stackrel{\sim}{\to} V'$ for an isomorphism of $\Cl_n$-bundles. 
\end{defn}

\begin{defn}\label{defn:stable}
Using the notation from Definition~\ref{defn:iso}, a \emph{stable equivalence} between $\Cl_n$-bundles $V$ and $V'$ is the data of $\Cl_n$-equivariant isometric inclusion $\varphi_{\alpha\sigma}\colon V_\sigma|_{U_\alpha'}\hookrightarrow V'_\alpha|_{U_\sigma}$ with a $\Cl_n\otimes \Cl_{-1}$-action on $V_\sigma^\perp\subset V'_\alpha|_{U_\sigma}$ for each $\sigma\subset I$ and $\alpha\subset I'$ with $U_\sigma\bigcap U'_\alpha\ne \emptyset$. These data are required to satisfy compatibility properties analogous to~\eqref{eq:compatdiagr}. 
\end{defn}

\begin{defn} A \emph{concordance} between $\Cl_n$-bundles $V$ and $V'$ on $M$ is a $\Cl_n$-bundle $\widetilde{V}$ on $M\times \R$ together with isomorphisms $i_0^*\widetilde{V}\cong V$ and $i_1^*\widetilde{V}\cong V'$, where $i_0,i_1\colon M\hookrightarrow M\times \R$ are the inclusions at~$0$ and~$1$ respectively. When such a concordance exists, we say that~$V$ and~$V'$ are \emph{concordant.} Concordance defines an equivalence relation on isomorphism classes of $\Cl_n$-bundles over~$M$ whose equivalence classes are called \emph{concordance classes}. 
\end{defn}

\begin{ex}\label{ex:homotopystable0}\label{ex:stable}
For vector bundles $V,V'\to M$ regarded as $\Cl_0$-modules, an isomorphism in the sense of Definition~\ref{defn:iso} is the same data as an isomorphism of vector bundles. If vector bundles $V,V'\to M$ are stably equivalent $\Cl_0$-bundles in the sense of Definition~\ref{defn:stable}, then $[V]=[V']\in \KO^0(M)$; indeed, $[V]-[V']=[V^\perp]\in \KO^0(M)$ and the existence of an odd isomorphism $e\colon V^\perp_\pm\to V_\mp^\perp$ implies $[V^\perp]=[V^\perp_\ev]-[V^\perp_\odd]=0$. In fact, a stable equivalence determines a concordance between $V$ and $V'$ as $\Cl_0$-modules, as constructed in the proof of the following. 
\end{ex}

\begin{lem} \label{lem:homotopystable}
If $\Cl_n$-bundles over $M$ are stable equivalent, then they are concordant. 
\end{lem} 
\bp
Let $V$ and $V'$ be stable equivalent $\Cl_n$-bundles over~$M$. By refining the cover, we can assume that $\Cl_n$-bundles $V$, $V'$, and the stable equivalence between them are defined subordinate to a fixed cover $\{U_i\}$ of~$M$. So let $V=\{V_\sigma,g_{\sigma\tau},e_{\sigma\tau}\}$ and $V'=\{V'_\sigma,g_{\sigma\tau}',e_{\sigma\tau}'\}$. Consider the open cover of $\R$ given by $\{(-\infty,1),(0,\infty)\}$. The concordance $\tilde{V}$ is defined via the basic $\Cl_n$-modules,
$$
V_\sigma\times (-\infty,1)\to U_\sigma\times(-\infty,1),\quad V_\sigma'\times (0,1)\to U_\sigma\times (0,1),\quad V'_\sigma\times (0,\infty)\to U_\sigma\times (0,\infty). 
$$
For inclusions $U_\sigma\times J\hookrightarrow U_\tau\times J$ with fixed $J=(-\infty,1)$, $(0,\infty)$ or $(0,1)$, we have the evident data $(g_{\sigma\tau},e_{\sigma\tau})$ or $(g_{\sigma\tau}',e_{\sigma\tau}')$ from $V$ or $V'$. It remains to specify the compatibility data (3) in Definition~\ref{defn:Clnbundle} for the inclusions 
$$
U_\sigma\times (-\infty,1)\hookleftarrow U_\sigma\times (0,1)\hookrightarrow U_\sigma\times (0,\infty). 
$$
For the right inclusion, the bundles are related by restriction, so the orthogonal complements are the zero bundle and there is no additional data to specify. For the left inclusion, we use the stable equivalence datum for $\varphi_\sigma\colon V_\sigma\times (0,1)\hookrightarrow V_\sigma'\times (0,1)$ to specify the compatibility data. This completes the construction of the concordance. 
\ep

\begin{rmk}
By composing concordances, a zig-zag of stable equivalences also determines a concordance between $\Cl_n$-bundles. 
\end{rmk}

\begin{defn}\label{defn:operations} Let $V=\{V_\sigma,g_{\sigma\tau},e_{\sigma\tau}\}$ be a $\Cl_n$-bundle on $M$ defined relative to an open cover $\{U_i\}$, and let $V'=\{V'_\alpha,g_{\alpha\beta}',e_{\alpha\beta}'\}$ be a $\Cl_m$-bundle on $M'$ defined relative to an open cover~$\{U_k'\}$. Define the following addition and multiplication operations on $\Cl_n$-bundles. 
\begin{enumerate}
\item If $n=m$, the \emph{external direct sum} is a $\Cl_n$-bundle $V\boxplus V':=\{V_\sigma\boxplus V_\alpha',g_{\sigma\tau}\boxplus g_{\alpha\beta}',e_{\sigma\tau}\boxplus e_{\alpha\beta}'\}$ on $M\times M'$ defined relative to the product open cover $\{U_i\times U_k'\}$, where $V_\sigma\boxplus V_\alpha'$ denotes the $\Z/2$-graded external direct sum of super vector bundles over $U_\sigma\times U'_\alpha$. 
\item The \emph{external tensor product} is a $\Cl_{n}\otimes \Cl_m$-bundle $V\boxtimes V'=\{V_\sigma\boxtimes V_\alpha',g_{\sigma\tau}\boxtimes g_{\alpha\beta}',e_{\sigma\tau}\boxtimes e_{\alpha\beta}'\}$ on $M\times M'$ defined relative to the open cover $\{U_i\times U_k'\}$, where $V_\sigma\boxtimes V_\alpha'$ denotes the $\Z/2$-graded external tensor product of super vector bundles over $U_\sigma\times U'_\alpha$. 
\end{enumerate}
Direct sum and tensor product operations for $\Cl_n$-bundles on $M$ are defined by pulling back the external direct sum and external tensor product along the diagonal~$M\hookrightarrow M\times M$. 
\end{defn}

\subsection{From $\Cl_n$-bundles to $\KO^n$-classes}

\begin{prop}\label{prop:ClnKO}
There is a surjective map from $\Cl_n$-bundles over $M$ to $\KO^n(M)$ that reduces to the Atiyah--Bott--Shapiro map~\eqref{eq:ABS} when $M=\pt$. A pair of $\Cl_n$-bundles determine the same class in $\KO^{n}(M)$ if and only if they are concordant. The external product and external sum of $\Cl_n$-bundles is compatible with the external product and external sum in~$\KO$. 
\end{prop}
\bp
Given a $\Cl_n$-bundle $V$ over $M$ defined relative to an open cover $\{U_i\}_{i\in I}$, consider the topological category $\check{C}(U_i)$ from Example~\ref{ex:opencover0}. 
The data in Definition~\ref{defn:Clnbundle} determines a map $\check{C}(U_i)_0 \to (\M_n)_0$ to the space of objects in the $n$th Atiyah--Bott--Shapiro category that is unique up to contractible choice, defined as follows. Choose $\Cl_n$-equivariant embeddings $\epsilon_\sigma\colon V_\sigma\subset \mathcal{H}_n\times U_\sigma$ over $U_\sigma$ for all $\sigma$. By~\cite{Kuiper}, the space of embeddings is contractible, so these choices are indeed unique up to contractible choice. Further specify homotopies between the embeddings
\beq\nonumber
\mathcal{H}_n\times U_\tau \stackrel{\epsilon_\sigma|_{U_\tau}}{\hookleftarrow}  V_\sigma|_{U_\tau} \stackrel{g_{\sigma\tau}}{\hookrightarrow} V_\tau \stackrel{\epsilon_\tau}{\hookrightarrow} \mathcal{H}_n\times U_\tau
\eeq
and higher homotopies between nested inclusions (again by~\cite{Kuiper}). 
This results in a map between realizations that is unique up to contractible choice
\beq
M\simeq |N_\bullet \check{C}(U_i)^\op|\to |N_\bullet\M_n|=\Bc\M_n.\label{eq:simpmap}
\eeq
By Lemma~\ref{lem:Kthy}, this determines a class $[V]\in \KO^n(M)$. 
By the nerve theorem~\cite{nerve}, any continuous map $M\to \Bc\M_n$ is homotopic to the realization of a map $N_\bullet \check{C}(U_i)^\op\to N_\bullet\M_n$ where $\{U_i\}$ is taken to be a good open cover; hence by Theorem~\ref{lem:Kthy} any class in $\KO^n(M)$ can be represented by a $\Cl_n$-bundle. 
Concordant $\Cl_n$-bundles determine homotopic maps to $\Bc\M_n$ by construction, and therefore represent the same $\KO$-class. Choosing a good open cover of $M\times \R$, homotopic maps can be lifted to concordant $\Cl_n$-bundles. 

The compatibility with external tensor product follows from Theorem~\ref{lem:Kthy} together with the fact that realization commutes with finite products; explicitly, the map~\eqref{eq:simpmap} associated with an external tensor product of $\Cl_n$- and $\Cl_m$-modules is homotopic to the composition
\beq
M\times M'&\simeq& |N_\bullet \check{C}(U_i\times U_j')^\op| \simeq |N_\bullet\check{C}(U_i)^\op\times N_\bullet\check{C}(U_i')^\op|\simeq |N_\bullet \check{C}(U_i)^\op|\times |N_\bullet \check{C}(U_i')^\op|\nonumber\\
&\to& \Bc\M_n\times \Bc\M_m \to \Bc\M_{n+m},\nonumber 
\eeq
where the first arrow is determined by the maps~\eqref{eq:simpmap} associated with the constituent $\Cl_n$- and $\Cl_m$-bundles, and the second arrow is the realization of~\eqref{eq:ABStensor}. Hence, this agrees with the product map $\KO^n(M)\times \KO^m(M')\to \KO^{n+m}(M\times M')$. 

The external sum is handled similarly, up to one subtlety coming from the fact that the sum operation on $\M_n$ is only partially defined (see Remark~\ref{rmk:productsum}). To ensure that the desired sums are defined, when choosing the maps to $\Bc\M_n$ associated with $\Cl_n$-bundles $V$ and $V'$, one chooses subspaces $V_\sigma\subset \mathcal{H}_n\times U_\sigma$ and $V_\alpha'\subset \mathcal{H}_n\times U_\alpha$ that are orthogonal to one another. Then the collection of subspaces $V_\sigma\boxplus V_\alpha'\subset \H_n\times U_\sigma\times U_\alpha'$ with the obvious compatibility data determines a map $M\times M'\to \Bc\M_n$ associated with the $\Cl_n$-bundle~$V\boxplus V'$. 
\ep

For a $\Cl_n$-bundle $V$ over $M$, we use the notation $[V]\in \KO^n(M)$ to denote the associated $\KO$-class from Proposition~\ref{prop:ClnKO}. 

\begin{cor} If $\Cl_n$-bundles $V$ and $V'$ over $M$ are stable equivalent, then $[V]=[V']\in \KO^n(M)$. \end{cor}

\bp This follows immediately from Lemma~\ref{lem:homotopystable}. \ep

\subsection{A representative of the suspension class}

\begin{defn} \label{defn:compactsupp}
A $\Cl_n$-bundle $V$ on $M$ has \emph{compact support} if on the complement of a compact subset $K\subset M$, $V$ has the zero $\Cl_n$-bundle as its underlying super vector bundle. 
Similarly, for a fiber bundle $\pi\colon M\to N$, a $\Cl_n$-bundle $V$ on $M$ has \emph{compact vertical support} if for each $y\in N$, there exists a neighborhood $U\subset N$ such that $V$ has compact support on~$\pi^{-1}(U)$. 
\end{defn}

\begin{cor}\label{cor:cpt}
There is a surjective map from compactly supported $\Cl_n$-bundles over~$M$ to $\KO^n_c(M)$. A pair of $\Cl_n$-bundles determine the same class in $\KO^{n}_c(M)$ if and only if they are concordant. The external product of $\Cl_n$-bundles is compatible with the external product in compactly supported~$\KO$. 
\end{cor}

\bp
This follows from the fact that the zero $\Cl_n$-bundle uniquely determines the constant map to basepoint of~$\Bc\M_n$ under~\eqref{eq:simpmap}. 
\ep

\begin{defn}
The \emph{suspension class} is a choice of generator $\KO^1_c(\R)\simeq \widetilde{\KO}{}^1(S^1)\simeq \Z$. 
\end{defn}

\begin{lem}\label{lem:cptsuspcocycle}
The suspension class is represented by the following $\Cl_1$-bundle. Consider
\beq\label{eq:nondiffsusp}
\begin{array}{c}
\{U_i\}=\{(-1,1), \ \R \setminus \{0\}\}\\
V_{\R\setminus \{0\}}={\bf 0}, \quad V_{(-1,1)}=\underline{\Cl}_1
\end{array}
\eeq
where $\Cl_1\simeq \R^{1|1}$ is a left module over itself, together with the stable isomorphism ${\bf 0 }\hookrightarrow \underline{\Cl}_1$ over $(-1,0)\coprod (0,1)$ from the $\Cl_1\otimes \Cl_{-1}$-action on $\R^{1|1}$ from Example~\ref{ex:Morita1}. 
\end{lem}

\bp
Let $\sigma_\R$ denote the $\Cl_1$-bundle~\eqref{eq:nondiffsusp}, and $[\sigma_\R]\in \KO^1_c(\R)$ the corresponding $\KO$-class using Corollary~\ref{cor:cpt}. The open cover $\{(-\infty,0), (-1,1), (0,\infty)\}$ is good, so (again by Corollary~\ref{cor:cpt}) the suspension class can be expressed in terms of basic $\Cl_n$-bundles on this cover together with stable isomorphisms on overlaps. The compact support condition requires that we take the zero module on $(-\infty,0)$ and $(0,\infty)$. On $(-1,1)$, a $\Cl_{1}$-module must be the trivial bundle with fiber the super vector space $\R^{k|k}\simeq (\R^{1|1})^{\oplus k}$ for some $k\in \N$; this is because $f\in \Cl_1$ gives an isomorphism between the even and odd subspaces of the bundle of $\Cl_1$-modules. On $(-1,0)$ and $(0,1)$, we require a $\Cl_1\otimes \Cl_{-1}\simeq \Cl_{1,1}$-module structure on $\R^{k|k}$. Up to parity reversal, $\R^{1|1}$ with $\Cl_{1,1}$-module structure from Example~\ref{ex:Morita1} is the unique irreducible $\Cl_{1,1}$-module. We conclude that any class in $\KO^1_c(\R)$ is of the form $[(\sigma_\R)^{\oplus k}]$ for some $k\in \Z$, where negative $k$ takes the parity reversed module. Hence, $[\sigma_\R]$ is a generator of $\KO^1_c(\R)\simeq \Z$, and therefore a representative of the suspension class. 
\ep

\section{The Pontryagin class of a $\Cl_n$-bundle as a \v{C}ech--de~Rham cocycle}\label{sec:superconn}

For an ordinary vector bundle $V\to M$, a choice of connection $\nabla$ on $V$ determines a closed differential form that refines the Pontryagin character $\Ch(V)\in \H(M;\R[u^{\pm 1}])$ of $[V]\in \KO^0(M)$, 
$$
\sTr(e^{-u\nabla^2})\in \Omega^\bullet(M;\R[u^{\pm 1}]),\qquad [\sTr(e^{-u\nabla^2})]=\Ch(V)\in \H(M;\R[u^{\pm 1}]).
$$
In this section we construct an analogous refinement for $\Cl_n$-bundles: after choosing connection data, a $\Cl_n$-bundle $V\to M$ determines a cocycle in the \v{C}ech--de~Rham complex of~$M$ that refines the Pontryagin character of~$[V]$; see Propositions~\ref{prop:CechdeRham} and~\ref{lem:yespontryagin}. 

\subsection{Review of Chern forms of superconnections} \label{sec:Quillen1}
For a real super vector bundle $V\to M$, let $V_\ev\subset V$ and $V_\odd\subset V$ denote its even and odd subbundles, respectively. Consider the space of $V$-valued differential forms $\Omega^\bullet(M;V)=\Gamma(M;\Lambda^\bullet(T^*M)\otimes V)$ as a super vector space.

\begin{defn}[{\cite{Quillensuperconn}}]
A \emph{superconnection} $\A$ on $V$ is an odd $\R$-linear map 
\beq
&&\A\colon \Omega^\bullet(M;V)\to \Omega^\bullet(M;V),\quad \A(\alpha\cdot v)=d\alpha \cdot v+(-1)^{|\alpha|} \alpha \cdot \A v\label{eq:ordsuperconn}
\eeq
satisfying the indicated Leibniz rule, where $v\in \Omega^\bullet(M;V)$ and $\alpha\in \Omega^\bullet(M)$. 
\end{defn}

On an $n$-manifold, one can express a superconnection as sum 
\beq
&&\A=\sum_{j= 0}^n \A_{(j)}, \label{eq:superconndegrees}
\eeq
where $\A_{(1)}=\nabla$ is an ordinary connection on $V$, $\A_{(2j)}\in  \Omega^{2j}(M;\End(V)_\odd)$ are valued in odd endomorphisms, and $\A_{(2j+1)}\in\Omega^{2j+1}(M;\End(V)_\ev)$ are valued in even endomorphisms. 

\begin{defn}[\cite{BGV}, page 117]\label{defn:sasuperconn}
Let $V\to M$ be a metrized real super vector bundle equipped with a superconnection~$\A$. The \emph{adjoint} superconnection is
\beq
\A^*=\sum_j(-1)^{j(j+1)/2} \A_{(j)}^\dagger\label{eq:superconnadjoint}
\eeq
where for $j\ne 1,$ $\A_{(j)}^\dagger$ is the adjoint of an endomorphism-valued form using the metric on~$V$, and $\A_{(1)}^\dagger=\nabla^\dagger$ is the adjoint of the ordinary connection characterized by 
$$
d\langle v,w\rangle=\langle \nabla v,w\rangle+\langle v,\nabla^\dagger w\rangle. 
$$
A superconnection is \emph{self-adjoint} if $\A^*=\A$. 
\end{defn}

\begin{rmk}\label{rmk:selfvsskew}
Self-adjointness is equivalent to $\A_{(j)}$ being self-adjoint in degrees~$j=4k,4k+3$, and skew-adjoint in degrees $j=4k+1, 4k+2$, $j\ne 1$. 
\end{rmk}

For $u^{1/2}$ a formal variable of degree $-1$, define 
\beq
&&\A(u):=u^{-1/2}\A_{(0)}+\A_{(1)}+u^{1/2}\A_{(2)}+\dots+ u^{j/2}\A_{(j)}+\dots\quad\label{eq:usuperconn} 
\eeq
Then the \emph{Pontryagin form of $\A$} is the differential form 
\beq\label{Eq:Chernform}
\Ch(\A):=\sTr(e^{-u\A(u)^2})\in \Omega^\bullet(M;\R[u,u^{-1}])
\eeq
of total degree~0, representing the Pontryagin character of~$V$. The absence of fractional powers of~$u$ in $\Ch(\A)$ follows from the fact that the super trace of an odd endomorphism is zero. Using Remark~\ref{rmk:selfvsskew}, self-adjointness further implies $\Ch(\A)$ takes values in differential forms in degrees 0 mod 4; see Lemma~\ref{lem:Chernbasics} below. The equalities of differential forms
\beq\label{eq:sumproduct}
&&\Ch(\A\oplus\A')=\Ch(\A)+\Ch(\A'),\quad \Ch(\A\otimes \A')=\Ch(\A)\cdot \Ch(\A')\in \Omega^\bullet(M;\R[u,u^{-1}])
\eeq
refine the direct sum and tensor product formulas for the Pontryagin character in cohomology.

\begin{defn}
Suppose we are given super vector bundles and superconnections $(V_0,\A_0)$ and $(V_1,\A_1)$ over $M$. A \emph{concordance} between the superconnections $\A_0$ and $\A_1$ is a superconnection $\tilde{\A}$ on a super vector bundle $\tilde{V}\to M\times \R$ with the data of isomorphisms $i_0^*\tilde{V}\simeq V_0$ and $i_1^*\tilde{V}\simeq V_1$ and the property that
$$
i_0^*\tilde{\A}=\A_0,\qquad i_1^*\tilde{\A}=\A_1
$$
for inclusions $i_j\colon M\times \{j\}\hookrightarrow M\times \R$, $j=0,1$. 
\end{defn}

\begin{defn}
The \emph{Chern--Simons form} of a concordance $\tilde{\A}$ of superconnections is 
\beq\label{eq:CSdefinition}
\CS(\tilde{\A}):=\int_{[0,1]} \Ch(\tilde{\A})\in \Omega^\bullet(M;\R[u,u^{-1}]).
\eeq
By Stokes' theorem, the Chern--Simons form satisfies
\beq\label{eq:CSStokes}
d\CS(\tilde{\A})=\Ch(\A_1)-\Ch(\A_0).
\eeq
\end{defn}

\begin{rmk}
Modulo exact forms $\CS(\tilde{\A})$ only depends on~$\A_1$ and~$\A_0$ and not on the choice of concordance $\tilde{\A}$. 
\end{rmk}

\begin{ex}
For superconnections $\A_0,\A_1$ on the same super vector bundle $V\to M$, affine interpolation
\beq\label{eq:affineinterp}
\tilde{\A}=dt\partial_t+(1-t)\A_0+t\A_1 \qquad  i_0^*\tilde{\A}=\A_0,\  \ i_1^*\tilde{\A}=\A_1
\eeq
gives a concordance on $p^*V\to M\times \R$ for $p\colon M\times \R\to M$ the projection. We use the notation
\beq\label{eq:specuialconmc}
\CS(\A_0,\A_1):=\CS(\tilde{\A})
\eeq
for concordances of the form~\eqref{eq:affineinterp}. 
\end{ex}

The existence of a concordance implies that $V_0\simeq V_1$, so often we shall implicitly assume that concordant superconnections $\A_0$ and $\A_1$ are defined on the same bundle. However, the concordance need not be of the form~\eqref{eq:affineinterp}. 

Using~\eqref{eq:sumproduct}, Chern--Simons forms have the following behavior under direct sums and tensor products
\beq\label{eq:CSproperties}
\begin{array}{ccl}
\CS(\tilde{\A}\oplus \tilde{\B})&=&\CS(\tilde{\A})+\CS(\tilde{\B})\\
\CS(\tilde{\A}\otimes \tilde{\B})&=&\CS(\tilde{\A})\Ch(\B)+\Ch(\A')\CS(\tilde{\B})\mod {\rm image}(d)\\
&=&\CS(\tilde{\A})\Ch(\B)+\Ch(\A)\CS(\tilde{\B})+d\CS(\tilde{\A})\CS(\tilde{\B}).\end{array}
\eeq

We require a generalization of~\eqref{eq:CSdefinition} to $(k+1)$-tuples of concordant superconnections. First define the affine smooth simplex $\Aff^k$ and the closed $k$-simplex $\Delta^k\subset \Aff^k$ as 
$$
\Aff^k=\{(t_0,t_1,\dots,t_k)\in \R^{k+1}\mid \sum t_i=1\}, \qquad \Delta^k=\{(t_0,t_1,\dots, t_k)\in \Aff^k\mid t_i\ge 0\}.
$$
Then given a super vector bundle $V\to M$ with $k+1$~superconnections $\A_0,\A_1,\A_2,\dots,\A_k$, define a superconnection $\tilde{\A}$ on $p^*V$ by
$$
\tilde{\A}:=\sum_{j=0}^k dt_j\partial_{t_j}+\sum_{j=0}^k t_jp^*\A_j
$$ 
where $p\colon M\times \Aff^k\to M$ is the projection. Then the differential forms 
\beq\label{eq:CSformsCdR2}
&&\CS_k(\tilde{\A})=\int_{\Delta^k} \Ch(\tilde{\A})\in \Omega^{\bullet}(M;\R[u^{\pm 1}])
\eeq
have total degree $-k$ and (by Stokes Theorem) satisfy 
\beq\label{eq:CSStokes2}
d\CS_k(\tilde{\A})=\int_{\partial\Delta^k} \Ch(\tilde{\A})=\sum_{j=0}^k (-1)^j\int_{\Delta^{k-1}} \Ch(\tilde{\A}_j)
\eeq
where $\tilde{\A}_j$ is the restriction of $\tilde{\A}$ to the $j$th face of $\Aff^k$. When $k=1$,~\eqref{eq:CSformsCdR2} reproduces the Chern--Simons form~\eqref{eq:CSdefinition}, where the property~\eqref{eq:CSStokes2} reduces to~\eqref{eq:CSStokes}.

\subsection{Clifford linear superconnections} \label{sec:naiveconn}
In~\cite[\S5]{Quillensuperconn}, Quillen uses superconnections to construct a differential form refinement of the Chern character for complex super vector bundles with $\cCl_1$-actions as an odd degree counterpart to~\eqref{Eq:Chernform}. The following is an elaboration on Quillen's ideas that treats real bundles with $\Cl_n$-actions. 

\begin{defn}
Let $V\to M$ be a real super vector bundle with fiberwise $\Cl_n$-action. A superconnection $\A$ on $V$ is \emph{$\Cl_n$-linear} if $\A$ graded commutes with the $\Cl_n$-action. 
\end{defn} 

\begin{defn}
Let $V\to M$ be a real, metrized super vector bundle with fiberwise self-adjoint $\Cl_n$-action. The \emph{Chern form} of a $\Cl_n$-linear self-adjoint superconnection $\A$ on $V$ is
\beq
&&\Ch_n(\A):=\sTr(u^{-n/2}\Gamma_n\circ \exp(-u\A(u)^2))\in \Omega^\bullet(M;\R[u^{1/2},u^{-1/2}]),\label{eq:Chernchar}
\eeq
using the super trace on $\Cl_n$-modules~\eqref{eq:Cliffordtrace} extended $\Omega^\bullet(M)$-linearly. 
\end{defn}

\begin{lem}\label{lem:Chernbasics}
For a self-adjoint $\Cl_n$-linear superconnection $\A$, the differential form $\Ch_n(\A)$ lies in the subalgebra $\Omega^\bullet(M;\R[\alpha,\alpha])\subset \Omega^\bullet(M;\R[u^{1/2},u^{-1/2}])$ where $\alpha=2u^{2}$. Furthermore, $\Ch_n(\A)$ is closed and of total degree $n$.
\end{lem}
\bp
The fact that $\Ch_n(\A)\in \Omega^\bullet(M;\R[u,u])\subset \Omega^\bullet(M;\R[u^{1/2},u^{-1/2}])$ follows from the super trace of an odd endomorphism being zero. In more detail,~\eqref{eq:superconndegrees} specifies whether a component of a superconnection takes values in even or odd endomorphisms; the odd endomorphisms correspond to the half-integer coefficients in~\eqref{eq:usuperconn}.  The endomorphism $\Gamma_n$ is even or odd depending on the parity of~$n$, and carries a factor of $u^{-n/2}$ from the definition of the Clifford super trace. Hence, the half-integer powers of~$u$ in the expansion of $\Ch_n(\A)$ are coefficients of odd endomorphisms and therefore sent to zero under the super trace. Self-adjoint superconnections give a further refinement to even powers of~$u$ using that the super trace of a skew-adjoint operator is zero, see Remark~\ref{rmk:selfvsskew}. This verifies that $\Ch_n(\A)\in \Omega^\bullet(M;\R[\alpha,\alpha])$. 

Since $u\A(u)^2$ has total degree zero and $u^{-n/2}\Gamma_n$ has total degree~$n$, $\Ch_n(\A)$ has total degree~$n$. Finally, to see that $\Ch_n(\A)$ is closed, one applies the same argument as in the case ordinary superconnections, e.g., \cite[Proposition~1.41]{BGV}. The key fact is 
\beq
\sTr_{\Cl_n}([\A,\omega])=d\sTr_{\Cl_n}(\omega). \label{eq:supertraceChern}
\eeq
This follows from a local computation, writing $\A=d+\gamma$ and using that the Clifford super trace vanishes on super commutators~\cite[Proposition~2]{Quillensuperconn}. Then the fact that $\Ch_n(\A)$ is closed follows from~\eqref{eq:supertraceChern} and the identity $[\A,\exp(-\A^2)]=0$. 
\ep

\begin{defn} A \emph{concordance} between self-adjoint $\Cl_n$-linear superconnections $\A_0,\A_1$ on $V_0,V_1\to M$ is a self-adjoint $\Cl_n$-linear superconnection $\tilde{\A}$ on $\tilde{V}\to M\times \R$
together with the data of isomorphisms
$$
i_0^*\widetilde{V}\simeq V_0,\qquad i_1^*\widetilde{V}\simeq V_1
$$
with the property that $i_0^*\tilde{\A}=\A_0$ and $i_1^*\tilde{\A}=\A_1$ for inclusions $i_j\colon M\times \{j\}\hookrightarrow M\times \R$, $j=0,1$.
\end{defn}

\begin{defn}
The \emph{Chern--Simons form} of a concordance between $\Cl_n$-linear self-adjoint superconnections is a differential form with total degree $n-1$, 
\beq\label{eq:CSsuperconn}
&&\CS(\tilde{\A}):=\int_{[0,1]} \Ch_n(\tilde{\A})\in \Omega^\bullet(M;\R[\alpha,\alpha^{-1}]).
\eeq
\end{defn}

The standard arguments show that the Chern--Simons form satisfies~\eqref{eq:CSStokes}, and only depends on~$\A_0$ and~$\A_1$ up to exact forms. We have the same formulas~\eqref{eq:CSproperties} for Chern--Simons forms involving direct sums and tensor products. 

\begin{lem}\label{lem:seCS}
Suppose we are given the data of a stable equivalence of  $\Cl_n$-bundles
\beq\label{eq:oddmeansequal}
V_0\hookrightarrow V_1,\qquad e\in \Gamma(M;\End(V^\perp)_\odd), \ e^2=\id
\eeq
together with $\Cl_n$-linear self-adjoint superconnections $\A_0$ on $V_0$ and $\A_1$ on $V_1$ (that need not be compatible with the inclusion $V_0\hookrightarrow V_1$). 
 Then there is a uniquely determined self-adjoint $\Cl_n$-linear superconnection $\A_0'$ on~$V_1$ with the equalities of differential forms of total degree~$n$
\beq\label{eq:chernformsequal2}
&&\Ch_n(\A_1)-\Ch_n(\A_0)=d\CS(\A_0',\A_1),\qquad  \Ch_n(\A_0)=\Ch_n(\A_0')\in \Omega^\bullet(M;\R[u^{\pm 1}])
\eeq
refining the equality in cohomology $\Ch(V_0)=\Ch(V_1)$ of Pontryagin characters.
In particular, a stable equivalence $V_0={\bf 0}\hookrightarrow V_1$ with a choice of super connection on $V_1$ determines a differential form~$\zeta$ satisfying $\Ch(\A_1)=d\zeta$. 
\end{lem}

\bp
Let $V_0^\perp\subset V_1$ denote the orthogonal complement to the inclusion $V_0\hookrightarrow V_1$. The existence of the $\Cl_n$-linear endomorphism $e\in \Gamma(\End(V_0^\perp)_\odd)$ with $e^2=\id_{V_0^\perp}$ implies that there is a direct sum decomposition of basic $\Cl_n$-bundles $V_0^\perp\simeq (V_0^\perp)_+\oplus (V_0^\perp)_-$ where $e\colon (V_0^\perp)_\pm \to (V_0^\perp)_\mp$ is a grading-reversing isomorphism that exchanges the summands. Let
$$
p^{+}\colon V_1\to(V_0^\perp)_{+}
$$
denote the orthogonal projection. Using the isomorphism of basic $\Cl_n$-bundles,
$$
V_1\simeq V_0\oplus V_0^\perp\simeq V_0\oplus (V_0^\perp)_{+}\oplus (V_0^\perp)_{-}
$$
define the self-adjoint $\Cl_n$-linear superconnection on $V_1$
\beq\label{eq:superconnprime}
\A_0'=\A_0\oplus (p^{+}\circ \nabla_1 \circ p^{+})\oplus e^*(p^{+}\circ \nabla_1 \circ p^{+})
\eeq
where $\nabla_1$ is the ordinary connection part of $\A_1$, see~\eqref{eq:superconndegrees}. Using~\eqref{eq:sumproduct} and properties of the Clifford supertrace, we compute
\beq
\Ch_n(\A_0')&=&\Ch_n(\A_0)+\Ch_n((p^{+}\circ \nabla_1 \circ p^{+}))+\Ch_n(e^*(p^{+}\circ \nabla_1 \circ p^{+}))\nonumber\\
&=&\Ch_n(\A_0)+\Ch_n((p^{+}\circ \nabla_1 \circ p^{+}))-\Ch_n((p^{+}\circ \nabla_1 \circ p^{+}))\nonumber\\
&=&\Ch_n(\A_0)\nonumber
\eeq
where the second line uses that $e$ is an odd automorphism. The result follows. \ep

\begin{ex}\label{ex:Bott} The modules in Example~\ref{eq:protoChern} can be regarded as basic $\Cl_n$-bundles over $M=\pt$. The Clifford super traces then compute the well-known Pontryagin characters of the associated classes in $\KO^\bullet(\pt)$. Specifically, 
\beq
&&\Ch_8({}_{\Cl_8}\R^{8|8})=u^{-4}=\beta^{-1}=\frac{1}{4}\alpha^{-2} \qquad \Ch_4({}_{\Cl_4}\R^{4|4})=2u^{-2}=\beta^{-1}\alpha=\frac{1}{4}\alpha^{-1},\label{eq:Cherncoeff}
\eeq
using the notation~\eqref{eq:gensrelation} for elements of $\KO^\bullet(\pt)$ and the identifications
$$
\Omega^\bullet(\pt;\R[\alpha,\alpha^{-1}])\simeq \R[\alpha,\alpha^{-1}]\subset \R[u,u^{-1}],\qquad 2u^2=\alpha.
$$
\end{ex}

\subsection{The Pontryagin form of a $\Cl_n$-bundle}\label{sec:CdR}

The Pontryagin form of a $\Cl_n$-bundle is locally determined by~\eqref{eq:Chernchar}, where the dependence of the differential form representative on the choice of superconnection is measured by~\eqref{eq:CSsuperconn}. Our present goal is to show that patching these local values together results in a cocycle in the \v{C}ech--de~Rham complex, $(\Omega^\bullet(U_*),d+\delta)$. Our conventions for the \v{C}ech--de~Rham double complex follow Bott and Tu~\cite[\S8]{BottTu}, and (in particular) the \v{C}ech differential $\delta$ requires an ordering on the indexing set $I$ of the cover. 

\begin{prop}\label{prop:CechdeRham}
Given a $\Cl_n$ bundle $V=\{V_\sigma,g_{\sigma\tau},e_{\sigma\tau}\}$ and an ordering on $I$, choices of self-adjoint superconnections $\A_\sigma$ on $V_\sigma\to U_\sigma$ for each~$\sigma \subset I$ determine a \v{C}ech--de~Rham cocycle denoted
\beq\label{eq:CechdeRham}
\Ch_n(\{U_i\},\A_i)\in \Omega^\bullet(U_*;\R[u^{\pm 1}]).
\eeq
\end{prop}

\begin{proof} 
Starting with the degree~$n$ Chern forms for each $i\in I$,
\beq\label{eq:ChernformsCdR}
\Ch_n(\A_i)\in \Omega^\bullet(U_i;\R[u^{\pm 1}])
\eeq
we construct an extension to a cocycle in the \v{C}ech--de~Rham complex. For $\sigma=\{i_0<i_1<\cdots <i_k\}$, the $\Cl_n$-bundle gives the data of stable isomorphisms $g_{i_j\sigma}\colon V_{i_j}|_{U_\sigma} \hookrightarrow V_\sigma$ over~$U_\sigma$. Using the chosen superconnections $\A_{i_j}$ on $V_{i_j}$ and $\A_\sigma$ on $V_\sigma\to U_\sigma$, Lemma~\ref{lem:seCS} produces superconnections~$\A_{i_j}'$ on $V_\sigma$ for each $i_j\in \sigma\subset I$, 
$$
\A_{i_0}',\A_{i_1}',\dots ,\A_{i_k}',\qquad \tilde{\A}_\sigma=\sum_{j=0}^k dt_j\partial_{t_j}+\sum_{j=0}^k t_j\A_j',\qquad i_0\le i_1\le \cdots \le i_k
$$ 
which determine the indicated superconnection $\tilde{\A}_\sigma$ on $p^*V_\sigma$ where $k+1=|\sigma|$ and $p\colon U_\sigma\times \Aff^k\to U_\sigma$ is the projection. Above, we order the superconnections $\A_{i_0}',\A_{i_1}',\dots ,\A_{i_k}'$ using the ordering on $I$. Following~\eqref{eq:CSformsCdR2}, define the higher Chern--Simons forms
\beq\label{eq:CSformsCdR}
&&\CS_k(\tilde{\A}_\sigma)=\int_{\Delta^k} \Ch_n(\tilde{\A}_\sigma)\in \Omega^{\bullet-k}(U_\sigma;\R[u^{\pm 1}])
\eeq
of total degree $n-k$. We note that changing the ordering of $I$ can change the sign of the form~\eqref{eq:CSformsCdR}.
The Chern forms~\eqref{eq:ChernformsCdR} and higher Chern--Simons forms~\eqref{eq:CSformsCdR} provide the necessary data for a degree~$n$ cocycle in the \v{C}ech--de~Rham complex; it remains to show that they satisfy the required properties. For $\delta$ the~\v{C}ech differential and $d$ the de~Rham differential in the \v{C}ech--de~Rham double complex, the first condition to check is
$$
\delta\{\Ch_n(\A_{i})\}=\delta\{\Ch_n(\widetilde{\A}_i)\}=\{d\CS(\A_i,\A_j)\}=\{d\CS(\tilde{\A}_{ij})\}\in \prod_{i<j} \Omega^{n}(U_{i}\bigcap U_{j})
$$
where the first and second equalities use~\eqref{eq:chernformsequal2}, and the third is the definition of~\eqref{eq:CSformsCdR}. The conditions on higher intersections follow similarly from 
\beq
d\int_{\Delta^k} \Ch_n(\A_{i_0,\dots,i_k})&=&\int_{\partial\Delta^k} \Ch_n(\A_{i_0,\dots,i_k})\nonumber\\
&=&\sum (-1)^j\int_{\Delta^{k-1}} \Ch_n(\A_{i_0,\dots,\widehat{i_j},\dots i_k})\nonumber\\
&=&\delta\int_{\Delta^{k-1}} \Ch_n(\A_{i_0,\dots,i_k}).\nonumber
\eeq
This shows we have produced a cocycle in the double complex, completing the proof.
\ep

\begin{prop} \label{lem:yespontryagin}
The \v{C}ech--de~Rham cocycle~\eqref{eq:CechdeRham} represents the Pontryagin character of the $\Cl_n$-bundle $V=\{V_\sigma,g_{\sigma\tau},e_{\sigma\tau}\}$. 
\end{prop}
\bp 
Since the space of superconnections is affine, $[\Ch_n(\{U_i\},\A_i)]$ only depends on~$[V]\in \KO^n(M)$. So by Proposition~\ref{prop:ClnKO},~$\Ch_n$ determines a map
\beq\label{eq:Pontrychar}
[\Ch_n]\colon \KO^n(M)\to \H^n(M;\R[\alpha,\alpha^{-1}]).
\eeq
The local description~\eqref{eq:Chernchar} along with the properties~\eqref{eq:Cliffordtraceplus} and~\eqref{eq:Cliffordtracetimes} of the Clifford supertrace imply that the map~\eqref{eq:Pontrychar} is additive and multiplicative. 

When $n=0$ (where $\Cl_0=\R$) the formula~\eqref{eq:Chernchar} coincides with the standard differential form refinement of the Pontryagin character for real vector bundles, verifying the lemma in this case. The statement for general~$n$ follows from the fact that~$\Ch_n$ is compatible with the suspension class and Bott element, as we now explain. Let $V=\{V_\sigma,g_{\sigma\tau},e_{\sigma\tau}\}$ be a $\Cl_j$-module, and choose $k\in \Z$ such that $8k-j\ge 0$. Then consider the diagram,
\beq
\begin{tikzpicture}[baseline=(basepoint)];
\node (A) at (0,0) {$\KO^j(M)$};
\node (B) at (5,0) {$\KO_c^{8k}(M\times \R^{8k-j})$};
\node (C) at (10,0) {$\KO_c^{0}(M\times \R^{8k-j})$};
\node (D) at (0,-1.5) {$\H^j(M;\R[u,u])$};
\node (E) at (5,-1.5) {$\H^{8k}_c(M\times \R^{8k-j};\R[u,u])$};
\node (F) at (10,-1.5) {$\H^{0}_c(M\times \R^{8k-j};\R[u,u])$};
\draw[->] (A) to node [above] {$\smallsmile [\sigma_{\R^{8k-j}}]$} (B);
\draw[->] (B) to node [above] {$\smallsmile [\beta^k]$} (C);
\draw[->] (A) to node [left] {$[\Ch_j]$} (D);
\draw[->] (B) to node [left] {$[\Ch_{8k}]$} (E);
\draw[->] (C) to node [left] {$[\Ch_0]$} (F);
\draw[->] (D) to node [below] {$\smallsmile [\Ch(\sigma_{\R^{8k-j}})]$} (E);
\draw[->] (E) to node [below] {$\smallsmile u^{4k}$} (F);
\path (0,-.75) coordinate (basepoint);
\end{tikzpicture}\nonumber
\eeq
The upper left horizontal arrow comes from the external product with the $\Cl_n$-bundle representing the suspension class $\sigma_{\R^{8k-j}}=(\sigma_\R)^{\boxtimes (8k-j)}$, where $\sigma_\R$ is defined in~\eqref{eq:nondiffsusp}. The image of this class under $\Ch_n$ is a \v{C}ech--de~Rham representative of the suspension class in ordinary cohomology, being compactly supported and nonzero in cohomology; compare Lemma~\ref{lem:suspensionch} below. This implies that the horizontal arrows on the left are both (suspension) isomorphisms. The middle horizontal arrows are the Bott isomorphisms, where Example~\ref{ex:Bott} computes the image of the Bott class under $\Ch_n$. The vertical arrow on the far right is the usual Pontryagin character in degree~0. Since all the horizontal arrows are isomorphisms, compatibility of the Pontryagin character in cohomology with suspension implies that the other vertical arrows are also given by the Pontryagin character. 
\ep

\section{Differential cocycles}\label{sec:dKOsec}
Our next goal is to refine the $\KO$-cocycles from Proposition~\ref{prop:ClnKO} to differential cocycles (see Definition~\ref{defn:superconngeneral}). The specifics are motivated by a generators and relations description of differential complex K-theory reviewed in~\S\ref{sec:FL} below. We then prove Theorem~\ref{thm1}, showing that $\dKO^\bullet(M)$ from Definition~\ref{defn:diffKO} is a model for differential $\KO$-theory.

\subsection{Review of differential complex K-theory}\label{sec:FL}
We review the geometric construction of differential complex K-theory from~\cite[Chapter~2]{Klonoff}. This is an extension of~\cite[Definition 2.14]{FreedLott} to super vector bundles with superconnection; see Remark~\ref{rmk:ordinaryconn}. For a smooth manifold $M$ consider triples $(V,\A,\phi)$ where $V\to M$ is a hermitian vector bundle, $\A$ is a self-adjoint superconnection on $V$, and $\phi\in \Omega^{\rm odd}(M)$. Define an equivalence relation 
\beq\label{eq:FLequiv}
&&(V,\A,\phi)\sim (V',\A',\phi')\iff  \begin{array}{l} {\rm there \ exists\ an \ isomorphism \ } f \colon V\xrightarrow{\sim} V' \\ {\rm with} \ \ \phi'=\CS(\A,f^*\A')+\phi\in \Omega^{\rm odd}(M)/\image(d).\end{array} 
\eeq
i.e., an isomorphism of hermitian vector bundles with a condition on the resulting Chern--Simons form.

\begin{defn} \label{defn:FL}
Define $\dK^0(M)$ as the quotient of the free abelian group generated by equivalence classes~\eqref{eq:FLequiv}, modulo the subgroup generated by
\beq\label{eq:FLequiv2}
&&(V,\A,\phi)+(V',\A',\phi')-(V\oplus V',\A\oplus \A',\phi+\phi')\quad{\rm and}\quad (W\oplus \Pi W, \nabla\oplus \Pi \nabla, 0)
\eeq
where $\Pi W$ is the parity reversal of a purely even vector bundle $W$, $\nabla$ is an ordinary connection on $W$, and $\Pi \nabla$ is the corresponding connection on $\Pi W$. 
A triple $(V,\nabla,\eta)$ is a \emph{differential cocycle} and $[V,\A,\phi]\in \dK^0(M)$ denotes its class in differential K-theory. 
\end{defn}

\begin{rmk}\label{rmk:ordinaryconn} In any given equivalence class~\eqref{eq:FLequiv} there is a representative involving an ordinary connection: every superconnection is concordant to an ordinary connection. 
\end{rmk}

We note that~\eqref{eq:FLequiv} implies there is a well-defined map 
$$
\dK^0(M)\to \Omega^{\rm even}(M),\qquad  [V,\nabla,\phi]\mapsto \Ch(V,\nabla,\phi):=\sTr(e^{-\nabla^2})+d\phi,
$$
valued in closed, even differential forms. For the generalization below, we emphasize that~$\phi$ is a coboundary mediating between the differential forms $\sTr(e^{-\nabla^2})$ and $\Ch(V,\nabla,\phi)$. 

\subsection{Differential refinements of $\Cl_n$-bundles}

\begin{defn} \label{defn:superconngeneral}
A \emph{differential refinement} of a $\Cl_n$-bundle $V=\{V_\sigma,g_{\sigma\tau},e_{\sigma\tau}\}$ defined relative to an open cover $\{U_i\}_{i\in I}$ is data
\begin{enumerate}
\item[(i)] for each $\sigma\subset I$, a self-adjoint $\Cl_n$-linear superconnection $\A_\sigma$ on $V_\sigma\to U_\sigma$, and
\item[(ii)] a closed global form $ \Ch(V,\A,\phi)\in \Omega^\bullet(M;\R[\alpha^{\pm 1}])$ of total degree~$n$;
\item[(iii)] a coboundary $\{\phi_\sigma\in \Omega^\bullet(U_\sigma;\R[\alpha^{\pm 1}])\}$ in the \v{C}ech--de~Rham complex mediating between the global form $\Ch(V,\A,\phi)$ and the cocycle $\Ch_n(\{U_i\},\A_i)\in \Omega^\bullet(U_*;\R[u^{\pm 1}])$ constructed in Proposition~\ref{prop:CechdeRham} (which depends on a choice of ordering on $I$). 
\end{enumerate}
A differential refinement of a $\Cl_n$-bundle is a \emph{differential cocycle} denoted by $(V,\A,\phi)=\{(V_\sigma,\A_\sigma,\phi_\sigma),(g_{\sigma\tau},e_{\sigma\tau})\}$.
\end{defn}

We unpack the data and property of the \v{C}ech--de~Rham coboundary in (iii). The data of this coboundary are differential forms $\phi_i\in \Omega^\bullet(U_i;\R[\alpha^{\pm1}])$  of total degree~$n-1$ for each~$i\in I$. These forms satisfy two properties. First, the locally defined quantities 
\beq
\Ch(\A_i)+d\phi_i\in \Omega^\bullet(U_i;\R[\alpha^{\pm1}]),\qquad \delta(\Ch(\A_i)+d\phi_i)=0\label{eq:determinesCh}
\eeq
determine the global form in (ii), i.e.,
$$
\Ch(V,\A,\phi)|_{U_i}-\Ch(\A_i)=d\phi_i\in \Omega^\bullet(U_i;\R[\alpha^{\pm 1}]),\qquad \Ch(V,\A,\phi)\in \Omega^\bullet(M;\R[\alpha^{\pm 1}]).
$$ 
Second, for all $i<j$, 
\beq
\CS(\tilde{\A}_{ij}) =\phi_j-\phi_i \in \Omega^\bullet(U_{ji};\R[\alpha,\alpha^{-1}]), \label{eq:compatibleChern}
\eeq
for the Chern--Simons form defined in~\eqref{eq:CSformsCdR}.

\begin{defn} A \emph{concordance} between differential cocycles $(V_0,\A_0,\phi_0),(V_1,\A_1,\phi_1)$ is a differential cocycle $(\tilde{V},\tilde{\A},\widetilde{\phi})$ over $M\times \R$ together with isomorphisms
$$
i_0^*\widetilde{V}\cong V_0, \ i_0^*(\tilde{\A},\widetilde{\phi})=(\A_0,\phi_0), \quad i_1^*\widetilde{V}\cong V_1, \ i_1^*(\tilde{\A},\widetilde{\phi})=(\A_1,\phi_1).
$$
satisfying the indicated compatibility with the superconnections and differential form data. 
\end{defn}

\begin{defn} \label{defn:CS}
Given a concordance between differential cocycles, define the \emph{total Chern--Simons form}
\beq
&&\CS((V_0,\A_0,\phi_0),(V_1,\A_1,\phi_1))=\int_{[0,1]} \Ch_n(\widetilde{V},\tilde{\A},\widetilde{\phi})\in \Omega^\bullet(M;\R[\alpha,\alpha^{-1}])\label{eq:CS1}
\eeq
of total degree $n-1$. Again by Stokes' theorem, the total Chern--Simons measures the difference between the Chern characters of concordant differential cocycles as in~\eqref{eq:CSStokes}. 
\end{defn}

For a smooth manifold $M$, consider the set of degree~$n$ differential cocycles modulo the following equivalence relation
\beq\label{eq:FLequiv3}
&&(V,\A,\phi)\sim (V',\A',\phi')\iff  \begin{array}{l} {\rm there \ exists\ a \ concordance \ (\tilde{V},\tilde{\A},\tilde{\phi})} \\ {\rm with} \ \ \CS((V,\A,\phi),(V,\A,\phi))\in \image(d).\end{array} 
\eeq

\begin{defn}\label{defn:diffKO} For a smooth manifold $M$, define $\dKO^n(M)$ as the quotient of the free abelian group generated by equivalence classes~\eqref{eq:FLequiv3} modulo the subgroup generated by 
\beq\label{eq:FLequiv4}
&&(V,\A,\phi)+(V',\A',\phi')-(V\oplus V',\A\oplus \A',\phi+\phi').
\eeq
Let $[V,\A,\phi]\in \dKO^n(M)$ denote the equivalence class of a differential cocycle in differential $\KO$-theory. 
\end{defn}

\begin{rmk} \label{ex:standard} 
Suppose $(V,\A,\phi)$ is a differential cocycle for which $V$ is a basic $\Cl_n$-bundle. Then by Stokes' theorem
\beq
\CS((V,\A,\phi),(V',\A',\phi'))&=&\int_{[0,1]} (\Ch_n(\widetilde{\A}) +d\widetilde{\phi})=\CS(\A,\A')+\phi'-\phi. \nonumber
\eeq
The condition~\eqref{eq:FLequiv3} asks that the right hand side be an exact form, which directly generalizes the equivalence relation~\eqref{eq:FLequiv} in differential complex K-theory.
\end{rmk}

\subsection{A ring structure on $\dKO^\bullet(M)$}
We define the following product structure on differential cocycles; compare~\cite[Equation~2.20]{FreedLott} or \cite[\S4.2]{Klonoff}.

\begin{defn}\label{eq:prodofCln}
Given a differential cocycle $(V,\A,\phi)$ on $M$ and a differential cocycle $(V',\A',\phi')$ on $M'$ define a differential cocycle on $V\boxtimes V'$ (see Definition~\ref{defn:operations}) as 
\beq
\{V_\sigma\boxtimes V_\alpha',\A_\sigma\boxtimes 1+1\boxtimes \A_\alpha', \phi_\sigma \Ch_n(\A_\alpha')+(-1)^{n}\Ch_n(\A_\sigma)\phi_\alpha'+\phi_\sigma d\phi_\alpha'\}\label{eq:superconnprodu}
\eeq
where $\A_\sigma\boxtimes 1+1\boxtimes \A_\alpha'$ is a superconnection on $V_\sigma\boxtimes V_\alpha'$ and $\phi_\sigma \Ch_n(\A_\alpha')+\phi_\alpha'\Ch_m(\A_\sigma)+\phi_\sigma d\phi_\alpha'\in \Omega^\bullet(U_\sigma\times U_\alpha';\R[\alpha,\alpha^{-1}])$. We use the shorthand $(V\boxtimes V',\A\boxtimes \A',\phi\boxtimes \phi')$ for the resulting differential cocycle. 
\end{defn}

\begin{lem}\label{lem:Chernchar} The assignment $(V,\A,\phi)\mapsto \Ch_n(V,\A,\phi)$ is compatible with sums and products:
\begin{enumerate}
\item For differential cocycles $(V,\A,\phi)$ and $(V,\A',\phi')$ of degree $n$ we have $\Ch_n(V\oplus V',\A\oplus \A',\phi+\phi')=\Ch_n(V,\A,\phi)+\Ch_n(V',\A',\phi')$.
\item For differential cocycles $(V,\A,\phi)$ and $(V,\A',\phi')$ of degree $n$ and $m$, we have $\Ch_{n+m}(V\boxtimes V',\A\boxtimes \A',\phi\boxtimes \phi')=\Ch_n(V,\A,\phi)\cdot\Ch_m(V',\A',\phi')$.
\end{enumerate}

\end{lem}

\bp
The statement can be checked locally, where it follows from the definition and properties~\eqref{eq:Cliffordtraceplus} and~\eqref{eq:Cliffordtracetimes} of the Clifford super trace. \ep

The following verifies that tensor products of differential cocycles endows $\dKO^\bullet(M)$ with a graded ring structure. 

\begin{lem} \label{lem:productlem}
If degree~$n$ differential cocycles $(V,\A,\phi)$ and $(V',\A',\phi')$ over $M$ are concordant, and degree~$m$ differential cocycles $(W,\B,\psi)$ and $(W',\B',\psi')$ over $N$ are concordant, then $(V\boxtimes W,\A\boxtimes \B,\phi\boxtimes\psi)$ is concordant to $(V'\boxtimes W',\A'\boxtimes \B',\phi'\boxtimes \psi')$ over $M\times N$. Furthermore, up to exact forms the total Chern--Simons form of the concordance is
\beq
\CS((V\boxtimes W,\A\boxtimes \B,\phi\boxtimes \psi),&&\!\!\!\!\!\!\!\!\!\!\!\!\!\!\!\!(V'\boxtimes W',\A'\boxtimes \B',\phi'\boxtimes \psi'))\nonumber\\
&=&\CS((V,\A,\phi),(V',\A',\phi'))\Ch_m(W,\B,\psi)\nonumber\\
&&+(-1)^{n}\Ch_n(V,\A,\phi)\CS((W,\B,\psi),(W',\B',\psi'))\nonumber\\
&&+\CS((V,\A,\phi),(V',\A',\phi))d\CS((W,\B,\psi),(W',\B',\psi')). \nonumber
\eeq
\end{lem} 

\bp
Let $\tilde{W}$ be the concordance over $M\times \R$ and $\tilde{V}$ be the concordance over $N\times \R$. Then we pullback $\tilde{W}\boxtimes \tilde{V}$ along the map 
$$
M\times N\times \R\stackrel{\id\times \Delta}{\longrightarrow} M\times N\times \R\times \R\simeq (M\times \R)\times (N\times \R)
$$
where $\Delta\colon \R\to \R\times \R$ is the diagonal map. This pullback gives the desired concordance between $(W\boxtimes V,\A\boxtimes \B,\phi\boxtimes \psi)$ and $(V'\boxtimes W',\A'\boxtimes \B',\phi'\boxtimes\psi')$. By Lemma~\ref{lem:Chernchar} part (2) and a short computation using~\eqref{eq:CSproperties}, the Chern--Simons form is as claimed.
\ep

\subsection{Differential refinements of stable equivalences}

\begin{defn}\label{defn:naiveconn}
Let $V_0,V_1\to M$ be basic $\Cl_n$-modules. 
A \emph{differential stable equivalence} between differential cocycles $(V_0,\A_0,\phi_0)$ and $(V_1,\A_1,\phi_1)$ is an inclusion $g\colon V_0\hookrightarrow V_1$ with the data of a $\Cl_n\otimes \Cl_{-1}$-action on $V_0^\perp\subset V_1$ extending the given $\Cl_n$-action such that 
\beq\label{eq:CSequals}
\CS(\A_0',\A_1)=\phi_0-\phi_1\in \Omega^\bullet(M;\R[u,u])^{n-1}/\image(d)
\eeq
for $\CS(\A_0',\A_1)$ the differential form constructed in Lemma~\ref{lem:seCS}. 
\end{defn}

\begin{lem}\label{lem:stableimpliesconcordant}
Let $V_0,V_1\to M$ be $\Cl_n$-modules, and suppose that differential cocycles $(V_0,\A_0,\phi_0)$ and $(V_1,\A_1,\phi_1)$ are differential stable equivalent in the sense of Definition~\ref{defn:naiveconn}. Then $(V_0,\A_0,\phi_0)$ and $(V_1,\A_1,\phi_1)$ are concordant with vanishing total Chern--Simons form. 
\end{lem}
\bp
Lemma~\ref{lem:homotopystable} constructs a concordance between the $\Cl_n$ bundles $V$ and $V'$, so it remains to construct a differential refinement of this concordance with the claimed properties. The superconnection data is specified by
$$
(p_-^*V_0,p_-^*\A_0,p_-^*\phi_0)\to M\times (-\infty,1),\quad (p^*V_1,p^*\A_1,p^*\phi_1)\to M\times (0,1),$$$$ (p_+^*V_1,p_+^*\A_1,p_+^*\phi_1)\to M\times (0,\infty)
$$
where $p\colon M\times (0,1)\to M$, $p_-\colon M\times (-\infty,1)\to M$ and $p_+\colon M\times (0,\infty)\to M$. The \v{C}ech--de~Rham coboundary is given by $\phi_0$ and $\phi_1$: the condition~\eqref{eq:compatibleChern} follows from the property~\eqref{eq:CSequals} of a stable equivalence, and the condition~\eqref{eq:determinesCh} follows from the equality
\beq\label{eq:diffconcodcond2}
\Ch(\A_0)+d\phi_0=\Ch(\A_1)+d\phi_1.
\eeq
Finally,~\eqref{eq:diffconcodcond2} also implies that the total Chern form of this differential cocycle over $M\times \R$ is constant in $\R$, being equal to the pullback of~\eqref{eq:diffconcodcond2} along the projection $M\times \R\to M$. Having no $dt$ component, the integral defining the Chern--Simons form vanishes identically. 
\ep

By Lemma~\ref{lem:seCS}, we obtain the following corollary that produces the claimed generalization of the second relation in~\eqref{eq:FLequiv2}. 
\begin{cor}\label{cor:trivsuper}
Let $(W,\A,\phi)$ be a differential cocycle over $M$. If $W$ has an extension to a $\Cl_n\otimes \Cl_{-1}$-bundle over $M$, then there exists a concordance of differential cocycles with vanishing total Chern--Simons form between $(W,\A,\phi)$ and $({\bf 0},0,\zeta)$, where ${\bf 0}$ is the zero bundle over~$M$ with zero connection~$0$, and $\zeta \in \Omega^\bullet(M;\R[u,u])^{n-1}/\image(d)$ satisfies
\beq
\Ch_n(W,\A,\phi)=d\zeta.\label{eq:trivsuper}
\eeq
\end{cor}

\subsection{A differential refinement of the $\KO$-suspension class}\label{sec:supsension}

The suspension isomorphism in de~Rham cohomology has a cocycle-level description as external multiplication with a differential form $U_{\R^n}$ 
\beq
&&\Omega^\bullet(M)\to \Omega_c^{\bullet+n} (M\times \R^n), \quad \alpha\mapsto \alpha\wedge U_{\R^n}, \quad 
U_{\R^n}\in \Omega^n_c(\R^n), \quad (2\pi)^{-n/2} \int_{\R^n} U_{\R^n}=1\label{eq:deRhamsusp0}\nonumber
\eeq
where $\Omega^n_c(\R^n)\subset \Omega^n(\R^n)$ denotes the subspace of rapidly decreasing forms; these forms give a model for compactly supported cohomology~\cite[\S4]{MathaiQuillen}. Any rapidly decreasing form with total integral~1 gives a suspension class (note our normalization above for this integral). A standard choice for $U_{\R^n}$ is the Gaussian
\beq
U_{\R^n}:=2^{n/2} e^{-|x|^2}d{\rm vol}.\label{eq:Gaussian}
\eeq
The suspension isomorphism $\KO^k(M)\to \KO^{k+n}_c(M\times \R^n)$ is similarly determined by external multiplication with a generator of $\KO^n_c(\R^n)\simeq \widetilde{\KO}{}^n(S^n)\simeq \Z$; see Lemma~\ref{lem:cptsuspcocycle}. We use Definition~\ref{defn:compactsupp} of a compactly supported $\Cl_n$-bundle in the following.

\begin{defn}
Define a \emph{rapidly decreasing degree~$k$ differential cocycle} over $\R^n$ as a compactly supported $\Cl_n$-bundle $V$ and a differential refinement $(V,\A,\phi)$ for which $\Ch_k(V,\A,\phi)\in \Omega(\R^n;\R[\alpha,\alpha^{-1}])$ and $\{\phi_\sigma\}$ are rapidly decreasing differential forms on $\R^n$. Consider the set of equivalence classes of rapidly decreasing differential cocycles modulo the relation~\eqref{eq:FLequiv3}. Define $\dKO^k_c(\R^n)$ as the quotient of the free abelian group on this set by the subgroup generated by~\eqref{eq:FLequiv4}. 
\end{defn}

By the above definition, there are maps
\beq
&&\KO^n_c(\R^n) \xleftarrow{[-]} \dKO^n_c(\R^n)\xrightarrow{\Ch_n} \Omega_c(\R^n;\R[\alpha,\alpha^{-1}])^n,\label{eq:curvunderlying}
\eeq
and so one can ask for a differential cocycle $\hat\sigma_{\R^n}\in \dKO^n_c(\R^n)$ that simultaneously refines the suspension class in $\KO$-theory and the differential form representative~\eqref{eq:Gaussian} of the suspension class in de~Rham cohomology. 

\begin{defn} \label{defn:suspcocycle}
Let $(\underline{\R}^{1|1},\A)$ denote the basic $\Cl_1$-bundle on $\R$ with super vector bundle  $\R^{1|1}\times \R\to \R$ and $\Cl_1$-action determined by the equation for $f$ in~\eqref{eq:Cl11bimod}. Equip this basic $\Cl_1$-bundle with the superconnection
\beq
\A=\A_{(0)}+\A_{(1)}=\left[\begin{array}{cc} 0 & x \\ x & 0\end{array}\right]+d=x e+ d,\label{eq:suspensionsuperconn}
\eeq
where $x\in \R$ is the standard coordinate, and we identify the odd endomorphism with the left action by $x e \in \Cl_{-1}$ from Example~\ref{ex:Morita1}. For $n\in \N$, we may also consider  the $\Cl_n$-bundle with superconnection on $\R^n$ given by the $n$th external product $(\underline{\R}^{1|1},\A)^{\boxtimes n}$. 
\end{defn} 

The following is an elaboration on computations in~\cite[\S2]{MathaiQuillen}.

\begin{lem}\label{lem:suspensionch}
We have the equality of differential forms 
$$
\Ch_n(\A^{\boxtimes n})=U_{\R^n} \in \Omega^{n}_c(\R^n)\label{eq:deRhamsusp}
$$
for $U_{\R^n}$ defined in~\eqref{eq:Gaussian}.
\end{lem}

\bp It suffices to prove the lemma in the case $n=1$. 
We compute
\beq
\Ch_1(\A)&=&\sTr(\Gamma_1 u^{-1/2} \exp(-u\A(u)^2))\nonumber\\
&=&\frac{1}{\sqrt{2}} \sTr\left(u^{-1/2}\left[\begin{array}{cc} 0 & 1 \\ -1 & 0\end{array}\right]\exp(-u\A(u)^2)\right)\nonumber\\
&=&\frac{2}{\sqrt{2}}e^{-x^2}dx=\sqrt{2} e^{-x^2}dx\in \Omega^1_c(\R;\R[u,u^{-1}]) \nonumber
\eeq
where we remind that the factor of $2^{-1/2}$ in the second line is from the definition of $\Gamma_1$ in~\eqref{eq:Gammadef}, and we used that 
$$
u \A(u)^2=\left[\begin{array}{cc} x^2 & 0 \\ 0 & x^2 \end{array}\right]+u^{1/2}\left[\begin{array}{cc} 0 & dx \\ dx & 0 \end{array}\right].
$$
This shows $\Ch_1(\A)=U_{\R}$ and the statement is proved. 
\ep

The $\Cl_n$-bundle underlying $(\underline{\R}^{1|1},\A)^{\boxtimes n}$ is not compactly supported in the sense of Definition~\ref{defn:compactsupp}, and hence cannot be used to define a differential refinement of the suspension class. However, $(\underline{\R}^{1|1},\A,0)^{\boxtimes n}$ is concordant to a compactly supported differential cocycle defined as follows; compare Lemma~\ref{lem:cptsuspcocycle}.

\begin{defn}\label{defn:cptsuspcocycle}
Define the \emph{suspension cocycle} $\hat\sigma_\R=(V,\A,\phi)\in \dKO^1_c(\R)$ as the equivalence class associated with the following differential cocycle on $\R$:
$$
\{U_i\}=\{(-1,1),\R\setminus \{0\}\}
$$
$$
V_{\R\setminus \{0\}}={\bf 0}, \ \A_{\R\setminus \{0\}}=0, \quad V_{(-1,1)}=\underline{\R}^{1|1}, \A_{(-1,1)}=d+xe,
$$
$$ \phi_{\R\setminus \{0\}}=(\sqrt{2}\int_{-\infty}^x e^{-|t|^2}dt\coprod \sqrt{2}\int_\infty^x e^{-|t|^2}dt), \quad \phi_{(-1,1)}=0. 
$$
The stable isomorphism data over $(-1,0)$ and $(0,1)$ comes from the $\Cl_1\otimes \Cl_{-1}$-action on $\R^{1|1}$ from Example~\ref{ex:Morita1}. 
For $n\in \N$, define the \emph{$n$th suspension cocycle} as external tensor product~$\hat\sigma_{\R^n}:=(\hat\sigma_\R)^{\boxtimes n}$. 
\end{defn}

\begin{lem}\label{lem:suspconc} There is a concordance between $\hat\sigma_{\R^n}$ and $(\underline{\R}^{1|1},\A,0)^{\boxtimes n}$ whose Chern--Simons form vanishes. 
\end{lem} 

\bp
On $\R\setminus \{0\}$, consider the canonical inclusion ${\bf 0}\hookrightarrow \underline{\R}^{1|1}$ with $\Cl_{-1}$-action on ${\bf 0}^\perp=\underline{\R}^{1|1}$ from Example~\ref{ex:Morita1}, and on $(1,-1)$ consider the identity inclusion $\underline{\R}^{1|1}\to \R^{1|1}$. Corollary~\ref{cor:trivsuper} verifies that this concordance of $\Cl_1$-bundles is a concordance between differential cocycles.
\ep

\begin{lem} \label{lem:suspension}
The $\KO$-class $[\hat\sigma_{\R^n}]\in \KO^n_c(\R^n)\simeq \Z$ underlying the $n$th suspension cocycle is the $n$th suspension class, i.e., a generator. 
\end{lem}
\bp
It suffices to prove the $n=1$ case, which follows from Lemma~\ref{lem:cptsuspcocycle}.
\ep

\begin{defn} The \emph{$n$th differential suspension map} is 
\beq
\boxtimes \hat\sigma_{\R^n}\colon \dKO^k(M)\to \dKO^{k+n}_{cv}(M\times \R^n),\qquad (W,\A,\phi)\mapsto (W,\A,\phi)\boxtimes \hat\sigma_{\R^n} \label{eq:susp}
\eeq
where the subscript $cv$ stands for compact vertical support relative to $M\times \R^n\to M$. 
\end{defn}

\subsection{A model for differential $\KO$-theory from $\Cl_n$-bundles}

\begin{lem} \label{lem:square} For each $n\in \Z$, the following square commutes: 
\beq
&&\begin{tikzpicture}[baseline=(basepoint)];
\node (A) at (0,0) {$\dKO^n(M)$};
\node (B) at (5,0) {$\KO^n(M)$};
\node (C) at (0,-1.5) {$\Omega_\cl(M;\R[\alpha,\alpha^{-1}])^n$};
\node (D) at (5,-1.5) {$\H^n(M;\R[\alpha,\alpha^{-1}]).$};
\draw[->] (A) to  (B);
\draw[->] (B) to node [right] {$[\Ch_n]$} (D);
\draw[->] (A) to node [left] {$\Ch_n$} (C);
\draw[->] (C) to (D);
\path (0,-.75) coordinate (basepoint);
\end{tikzpicture}\label{diag:square}
\eeq
where the upper horizontal arrow sends a differential cocycle in $\dKO^n(M)$ to its underlying class in $\KO^n(M)$ (via Proposition~\ref{prop:ClnKO}), the lower horizontal arrow sends a closed form to its de~Rham cohomology class, and the arrow $\Ch_n$ extracts the datum of the closed differential form $\Ch_n(V,\A,\phi)\in \Omega_\cl(M;\R[\alpha,\alpha^{-1}])^n$ of total degree~$n$.
\end{lem}
\bp 
This follows from Proposition~\ref{lem:yespontryagin} and the definition of a differential cocycle. 
\ep

\begin{lem}\label{lem:exact}
There is an exact sequence 
\beq
&&\KO^{n-1}(M)\stackrel{[\Ch]}{\to} \Omega^{n-1}(M;\R[\alpha,\alpha^{-1}])/\image(d) \stackrel{a}{\to} \dKO^n(M)\to \KO^n(M)\to 0
\eeq
where the map $a(\phi)=[{\bf 0},0,\phi]\in \dKO^n(M)$ sends a differential form $\phi$ of total degree $n-1$ to the zero $\Cl_n$-bundle with differential refinement determined by $\phi$.
\end{lem}
\bp
We work from right to left. Contractibility of spaces of superconnections and differential forms implies that the map $\dKO^n(M)\to \KO^n(M)$ is surjective. 
Exactness at $\dKO^n(M)$ follows from Corollary~\ref{cor:trivsuper}. Next we check exactness at $\Omega^{n-1}(M;\R[\alpha,\alpha^{-1}])$. By Proposition~\ref{prop:ClnKO}, it suffices to show that for all degree~$n-1$ differential cocycles $(V,\A,\phi)$ the degree~$n$ differential cocycle $({\bf 0},0,\Ch_{n-1}(V,\A,\phi))$ is concordant to zero. In turn, this requires we find a concordance from the zero $\Cl_n$-bundle and a differential refinement of this concordance whose Chern--Simons form is $\Ch_{n-1}(V,\A,\phi)$. Apply the differential suspension map~\eqref{eq:susp} to $(V,\A,\phi)$. This yields a compactly supported degree~$n$ differential cocycle over $M\times \R$, which (after a reparameterization of $\R$) determines a concordance from the zero $\Cl_n$-bundle to itself. By Lemma~\ref{lem:Chernchar} and the fact that $\Ch(\sigma_\R)=\U_\R$ from~\eqref{eq:Gaussian} integrates to~1, the Chern--Simons form of this concordance is precisely $\Ch_{n-1}(V,\A,\phi)$. The lemma is proved. 
\ep

\begin{prop}\label{prop:diffsusp} The $k$th differential suspension map~\eqref{eq:susp} is an isomorphism. 
\end{prop}
\bp
 It suffices to prove the statement for $k=1$. Consider the diagram 
\beq
\begin{tikzpicture}[baseline=(basepoint)];
\node (AA) at (-4.5,0) {$\KO^{n-1}(M)$};
\node (BB) at (-4.5,-1.5) {$\KO^{n}_c(M\times \R)$};
\node (A) at (0,0) {$\Omega^{n-1}(M;\R[\alpha,\alpha^{-1}])/\image(d)$};
\node (B) at (5,0) {$\dKO^n(M)$};
\node (C) at (8,0) {$\KO^n(M)$};
\node (D) at (0,-1.5) {$\Omega^{n}_c(M\times \R;\R[\alpha,\alpha^{-1}])/\image(d)$};
\node (E) at (5,-1.5) {$\dKO_c^{n+1}(M\times \R)$};
\node (F) at (8,-1.5) {$\KO_c^{n+1}(M\times \R),$};
\draw[->] (AA) to node [left] {$[\sigma_\R]$} (BB);
\draw[->] (AA) to node [above] {$[\Ch_{n-1}]$} (A);
\draw[->] (BB) to node [above] {$[\Ch_n]$} (D);
\draw[->] (A) to (B);
\draw[->] (B) to (C);
\draw[->] (A) to node [left] {$(-)\wedge U_\R$} (D);
\draw[->, bend right] (D) to node [right] {$(2\pi)^{-1/2}\int_\R$} (A);
\draw[->,dashed] (B) to node [left] {$\sigma_\R$} (E);
\draw[->] (C) to node [left] {$[\sigma_\R]$} (F);
\draw[->] (D) to  (E);
\draw[->] (E) to  (F);
\path (0,-.75) coordinate (basepoint);
\end{tikzpicture}\nonumber
\eeq
where the vertical arrows on the far left and far right are the suspension isomorphism in (non-differential) $\KO$-theory. The arrow given by the external tensor product with the rapidly decreasing 1-form $U_\R=\sqrt{2} e^{-x^2}dx$ has an inverse given by integration over $\R$. Hence, three of the vertical arrows are isomorphisms. By exactness of the rows and surjectivity of the horizontal maps on the far right, we conclude the dashed vertical arrow is an isomorphism by the five lemma. 
\ep

\begin{proof}[Proof of Theorem~\ref{thm1}] 
The result follows from the uniqueness of multiplicative differential extensions for rationally even cohomology theories~\cite[Theorem 1.7]{BunkeSchick}, using that $\dKO^\bullet(-)$ is multiplicative together with Lemmas~\ref{lem:square} and~\ref{lem:exact}. 
\ep

\begin{rmk}
One can also prove Theorem~\ref{thm1} directly without appealing to Bunke and Schick's unicity theorem. To sketch the idea, for $n=0$ Theorem~\ref{thm1} follows from the same arguments for the generators and relations presentation for complex K-theory in degree~0, e.g., see~\cite[\S2.2]{FreedLott}. 
With the $n=0$ case in hand, the general statement follows the differential suspension isomorphism using the same argument as~\cite[Proposition 9.21]{FreedLott}. Yet another proof comes from constructing an isomorphism between $\dKO^\bullet$ and the differential refinement constructed by Hopkins and Singer~\cite[\S4]{HopSing} as a homotopy pullback; for this comparison, one uses the point-set model for the $\KO$-spectrum from~\S\ref{sec:ABS}, together with the description of the Pontryagin character from~\S\ref{sec:superconn}. 
\end{rmk}

\section{The families differential index}

\subsection{Review of the (families) analytic index} \label{sec:Cliffordindex}
Let $X$ be an $n$-dimensional Riemannian spin manifold with principal $\Spin(n)$-bundle~$P\to X$. Define the vector bundle~\cite[II.7]{LM}
\beq
\bS_X:=P\times_{\Spin(n)}\Cl_n\simeq P\times_{\Spin(n)}(\Cl_n^\ev\oplus \Cl_n^\odd)\label{eq:Riemspinbundle}
\eeq
with $\Z/2$-grading inherited from the grading on the Clifford algebra as indicated. The fibers of $\bS_X$ are invertible $\Cl(TX)-\Cl_n$-bimodules. The Dirac operator is the composition
\beq
&&\slashed{D}\colon \Gamma(\bS_X)\stackrel{\nabla}{\to} \Gamma(T^*X)\otimes_{C^\infty(X)} \Gamma (\bS_X)\stackrel{\sim}{\to} \Gamma(TX)\otimes_{C^\infty(X)} \Gamma (\bS_X)\stackrel{\cl}{\to} \Gamma(\bS_X),\label{eq:defnDirac}
\eeq
where $\nabla$ is the Levi-Civita connection, the middle isomorphism uses the metric, and the final arrow on the right is (left) Clifford multiplication. Hence, $\slashed{D}$ commutes with the right $\Cl_n$-action on~$\Gamma(\bS_X)$. Given a real vector bundle $V\to X$ with metric and compatible connection~$\nabla^V$, the $V$-twisted $\Cl_n$-linear Dirac operator $\slashed{D}\otimes V$ is defined similarly on~$\Gamma(\bS_X\otimes V)$. After trading the right $\Cl_n$-action in~\eqref{eq:Riemspinbundle} for a left $\Cl_{-n}$-action (using $\Cl_n^\op=\Cl_{-n}$), the Dirac operator $\slashed{D}\otimes V$ is a $\Cl_{-n}$-linear unbounded operator on the $L^2$-completion of sections $L^2(\bS_X\otimes V)$, which we identify with a stable $\Cl_{-n}$-module. The space of identifications $L^2(\bS_X\otimes V)\simeq \mathcal{H}_{-n}$ is contractible by~\cite{Kuiper}. The domain of~$\slashed{D}\otimes V$ is the dense subspace of smooth spinors, $\Gamma(\bS_X\otimes V)\subset L^2(\bS_X\otimes V)$. When the spin manifold is compact, $\slashed{D}\otimes V$ has compact resolvent and so determines a point in the space $\Inf_{-n}$ of unbounded operators defined in~\S\ref{sec:unbounded}. Under Theorem~\ref{thm:HST} we then obtain a class, 
\beq
[\slashed{D}\otimes V]\in \KO^{-n}(\pt).\label{eq:analyticindexreview}
\eeq
By Lemma~\ref{lem:Kthy}, this class equals the one computed by the Clifford index~\eqref{eq:KOindex}. 

The families generalization of~\eqref{eq:analyticindexreview} starts with the following geometric data. 

\begin{defn}[{\cite[\S10.1]{BGV}}] \label{defn:spinfamily} A \emph{family of Riemannian spin manifolds} is
\begin{enumerate}
\item a smooth fibration of manifolds $\pi\colon X\to B$ whose fibers are compact of constant dimension~$n$ and have chosen spin structures;
\item a metric $g^{X/B}$ on the vertical tangent bundle $T(X/B):={\rm ker}(d\pi)\subset TX$; and
\item a projection $p\colon TX\to T(X/B)$. 
\end{enumerate}
\end{defn}

Letting $P\to X$ denote the principal spin bundle for $T(X/B)\to X$, we define the families spinor bundle $\bS:=P\times_{\Spin(n)}\Cl_n$ over~$X$ analogously to~\eqref{eq:Riemspinbundle}. Given a real vector bundle $V\to X$ with metric and compatible connection, we may form the $V$-twisted spinor bundle $\bS\otimes V$. The associated $\Cl_{-n}$-linear Dirac operator $\slashed{D}\otimes V$ is defined in the same manner as before. Restricting $\bS\otimes V$ and $\slashed{D}\otimes V$ to a fiber $X_b:=\pi^{-1}(\{b\})$ at $b\in B$, we recover the spinor bundle $\bS_{X_b}\otimes V|_{X_b}$ and Dirac operator $\slashed{D}\otimes V|_{X_{b}}$ for the Riemannian spin manifold~$X_b$. 

The direct image $\pi_*(\bS\otimes V)$ determines a smooth $\Z/2$-graded Fr\'echet vector bundle on~$B$~\cite[\S9.2]{BGV}. The fibers of this bundle have the structure of left $\Cl_{-n}$-bundles, and the $\Cl_{-n}$-linear Dirac operator determines a family of Clifford linear odd operators. The fiberwise $L^2$-completion of $\pi_*(\bS\otimes V)$ determines a continuous (but not smooth) $\Z/2$-graded Hilbert bundle over~$B$. By Kuiper's theorem~\cite{Kuiper}, we may choose an identification with~$\mathcal{H}_{-n}\times B\to B$, the trivial bundle. The family of Dirac operators then determines a continuous map 
\beq
B\to \Inf_{n}\qquad b\mapsto (\slashed{D}\otimes V)_b\in \Inf_{n}\label{eq:maptoInf}
\eeq
for the spaces $\Inf_{n}$ constructed in~\cite{HST} and reviewed in~\S\ref{sec:unbounded}. 

\begin{defn} Define the \emph{families analytic index} of $\pi\colon X\to B$ as the $\KO$-class
\beq
[\slashed{D}\otimes V]\in [B,\Inf_{n}]\simeq \KO^{-n}(B)\label{eq:familiesindex}
\eeq
determined by the map~\eqref{eq:maptoInf}. 
\end{defn} 

\begin{rmk}
We learned the above description of the families analytic index from Stolz and Teichner, see~\cite[Remark 3.2.21]{ST04} and~\cite[\S3 and \S6]{HST}. It is equivalent to the more standard Clifford linear Fredholm index from~\cite[III.16]{LM} using the comparison results~\cite[Proposition 3.11 and Theorem 8.2]{HST}, stated as Theorem~\ref{thm:HST} above. 
\end{rmk}

\subsection{The families index as a $\Cl_{-n}$-bundle}\label{sec:familiesindex2}

We refine the families analytic index~\eqref{eq:familiesindex} to a $\Cl_{-n}$-module that is compatible with~\eqref{eq:familiesindex} via Proposition~\ref{prop:ClnKO}. 

\begin{prop}\label{thm:Cliffordindex}
Given a family of $n$-dimensional Riemannian spin manifolds~$\pi\colon X\to B$ and a real vector bundle $V\to X$ with metric and compatible connection, there is a $\Cl_{-n}$-bundle on $B$ that is unique up to concordance representing the $\KO$-class $[\slashed{D}\otimes V]\in \KO^{-n}(B).$ 
\end{prop}
\bp
The existence of a continuous $\Cl_{-n}$-bundle follows formally from the construction~\eqref{eq:maptoInf} and Lemma~\ref{lem:Kthy}. A smooth refinement critically uses properties of Dirac operators. 
For the family $\slashed{D}\otimes V$ of Dirac operators acting on the (continuous) Hilbert bundle $\mathcal{H}\to B$, 
define an ordered open cover via subsets $U_\lambda\subset B$ for $\lambda\in \R_{>0}$
\beq
U_\lambda=:\{b\in B\mid \lambda\notin {\rm Spec}(\slashed{D}_b^2)\}\subset B.\label{eq:Ulambda2}
\eeq
For a fixed $\lambda\in \R_{>0}$ and $b\in U_\lambda$, define
\beq
 \mathcal{H}^{<\lambda}_{b}:=\bigoplus_{\nu<\lambda} \mathcal{H}_{b}^\nu \subset \mathcal{H}_b\label{eq:bundleofClnmodules}
 \eeq
as the direct sum of $\nu$-eigenspaces $\mathcal{H}_{b}^\nu$ of $\slashed{D}^2_b$ with eigenvalue $\nu<\lambda$. Let $\mathcal{H}^{<\lambda}$ denote the $\Cl_{-n}$-invariant subset of~$\mathcal{H}|_{U_\lambda}\simeq U_\lambda\times \mathcal{H}_{-n}$ whose fiber at $b$ is~\eqref{eq:bundleofClnmodules}. Ellipticity implies that the subsets~\eqref{eq:Ulambda2} are smooth submanifolds and the subsets~\eqref{eq:bundleofClnmodules} are smooth, finite-rank bundles of $\Cl_{-n}$-modules on $U_\lambda$; see~\cite[Proposition~9.10]{BGV} or~\cite[III, Theorem~5.8]{LM}.

We shall define a $\Cl_{-n}$-bundle on $B$ subordinate to the open cover~$\{U_\lambda\}_{\lambda\in \Lambda}$ for any choice of subset $\Lambda\subset \R_{>0}$ that determines a locally finite open cover of~$B$. For $\sigma\subset \Lambda\subset \R_{>0}$ a finite subset determining a nonempty intersection $U_\sigma$, let $\lambda\in \sigma\subset \R_{>0}$ be the largest element and assign the bundle $\mathcal{H}^{<\lambda}|_{U_\sigma}\to U_\sigma$ from restricting along the inclusion $U_\sigma\subset U_\lambda$. Over $U_\lambda\bigcap U_\mu$ for $\lambda<\mu$, there is a canonical inclusion of bundles of $\Cl_{-n}$-modules 
\beq
g_{\lambda\mu}\colon \mathcal{H}^{<\lambda}|_{U_\mu}\hookrightarrow \mathcal{H}^{<\mu}|_{U_\lambda}.\label{eq:compatibilitydata}
\eeq
Indeed, there is an orthogonal decomposition $\mathcal{H}^{<\mu}\simeq \mathcal{H}^{<\lambda}\oplus \mathcal{H}^{(\lambda,\mu)}$ over $U_\lambda\bigcap U_\mu$, where $\mathcal{H}^{(\lambda,\mu)}$ is the finite-rank smooth vector bundle whose fiber at $b\in U_\lambda\bigcap U_\mu$ is the direct sum of eigenspaces of $\slashed{D}^2_b$ with eigenvalues in $(\lambda,\mu)$. Hence, $(\mathcal{H}^{<\lambda})^\perp=\mathcal{H}^{(\lambda,\mu)}$, for the orthogonal complement taken relative to the inclusion~\eqref{eq:compatibilitydata}. Define a family of odd endomorphisms of $\mathcal{H}^{(\lambda,\mu)}$ that at $b\in U_\lambda\bigcap U_\mu$ is given by a direct sum of linear maps on each eigenspace
\beq
e_{\lambda\mu}(b):=\bigoplus_{\nu\in (\lambda,\mu)} \frac{1}{\sqrt{\nu}}\slashed{D}_b|_{\mathcal{H}^\nu_{b}},\qquad e_{\lambda\mu}^2=\id_{\mathcal{H}^{(\lambda,\mu)}} \label{eq:CLminus1action}
\eeq
where sum is indexed by eigenvalues $\nu\in (\lambda,\mu)$ of $\slashed{D}_b^2$, and $e_{\lambda\mu}$ generates a self-adjoint $\Cl_{-1}$-action on~$\mathcal{H}^{(\lambda,\mu)}$. For arbitrary $\sigma\subset \tau\subset \Lambda$, the above data restrict to give the data of a $\Cl_{-n}$-bundle associated to the inclusion $U_\tau\subset U_\sigma$. This completes the construction of a~$\Cl_{-n}$-bundle associated to a family of spin manifolds. 

A priori the above construction depends on a choice of subcover $\{U_\lambda\}_{\lambda\in \Lambda}$ for $\Lambda\subset \R_{>0}$. For a different choice of subcover $\{U_\lambda\}_{\lambda\in \Lambda'}$, we can compare the resulting $\Cl_{-n}$-bundles on the mutual refinement $\{U_\lambda\}_{\lambda\in (\Lambda\bigcup \Lambda')}$. For $\lambda<\lambda'$, we obtain inclusions~\eqref{eq:compatibilitydata} of $\Cl_{-n}$-bundles with commuting $\Cl_{-1}$-actions on the orthogonal complement to the inclusion. This gives the data of a stable equivalence of $\Cl_{-n}$-bundles. By Lemma~\ref{lem:homotopystable}, the $\Cl_{-n}$-bundles are therefore concordant. 

The $\Cl_{-n}$-bundle constructed above determines a continuous functor $\check{{\rm C}}(U_\lambda)^\op\to \M_{-n}$ that (by comparing with Lemma~\ref{lem:Kthy}) factors the map $B\to \Inf_{n}$ in~\eqref{eq:familiesindex}. Hence, the associated $\KO$-class agrees with the families analytic index.
\ep

\subsection{A differential refinement of families index }\label{sec:familiesindex}

Following Bismut~\cite{Bismutindex}, the geometric data in Definition~\ref{defn:spinfamily} allows for the construction of a superconnection on the infinite rank bundle $\pi_*(\bS\otimes V)$ whose degree zero part is the family of Dirac operators $\slashed{D}\otimes V$. To review this, let $\Omega^\bullet(B;\pi_*(\bS\otimes V))$ denote the real Fr\'echet space of differential forms with values in the $\Z/2$-graded Fr\'echet bundle $\pi_*(\bS\otimes V)$. We regard $\Omega^\bullet(B;\pi_*(\bS\otimes V))$ as a left $\Cl_{-n}$-module. Following~\cite[Definition~9.12]{BGV}, a \emph{superconnection on $\pi_*(\bS\otimes V)$ adapted to $\slashed{D}\otimes V$} is a differential operator $\A\colon \Omega^\bullet(B;\pi_*(\bS\otimes V))\to \Omega^{\bullet}(B;\pi_*(\bS\otimes V))$ of odd degree satisfying a Leibniz rule and taking the form
$$
\A=\slashed{D}\otimes V+\sum_{i=1}^{{\rm dim}(B)} \A_i
$$ 
i.e., the degree zero part of the superconnection is the family of Dirac operators. A superconnection is \emph{$\Cl_{-n}$-linear} if it graded commutes with the $\Cl_{-n}$-action on $\Omega^\bullet(B;\pi_*(\bS\otimes V))$. 

The Fr\'echet bundle $\pi_*(\bS\otimes V)$ carries a unitary connection $\nabla^{\pi_*(\bS\otimes V)}$ constructed from the Levi-Civita connection on $\bS$, the connection on $V$, and the mean curvature of the fibers of $\pi\colon X\to B$~\cite[Proposition~9.13]{BGV}. The \emph{Bismut superconnection} is defined as~\cite[Proposition~10.15]{BGV}
\beq
\B:=\slashed{D}\otimes V+\nabla^{\pi_*(\bS\otimes V)}-\frac{c(T)}{4}\label{eq:Bismutsuper}
\eeq
where $c(T)$ is Clifford multiplication by the curvature 2-form $T$ of the bundle $\pi\colon X\to B$. The Bismut superconnection is $\Cl_{-n}$-linear: each of the three terms on the right of~\eqref{eq:Bismutsuper} independently commutes with the $\Cl_{-n}$-action, where it is important to note that the operator $c(T)$ uses the left $\Cl(TX)$-action, which commutes with the left $\Cl_{-n}$-action. 

\begin{proof}[Proof of Theorem~\ref{thm2}] 

It remains to construct a differential refinement of the $\Cl_{-n}$-bundle $\{\mathcal{H}^{<\lambda},g_{\lambda\mu},e_{\lambda\mu}\}$ from Proposition~\ref{thm:Cliffordindex}.
Following Definition~\ref{defn:superconngeneral}, we will supply data:
\begin{enumerate}
\item[(i)] ordinary $\Cl_{-n}$-linear connections $\nabla^{<\lambda}$ on each $\mathcal{H}^{<\lambda}$; 
\item[(ii)] a global differential form of total degree $-n$; and
\item[(iii)] differential forms $\eta_\lambda\in \Omega^\bullet(B;\C[u,u^{-1}])$ determining a \v{C}ech--de~Rham coboundary between the Pontryagin character associated with (i) and the global form from (ii). 
\end{enumerate}
%Parts (i) and (iii) analysis of the $\Cl_{-n}$-linear Bismut superconnection~\eqref{eq:Bismutsuper}. 

For (i), define connections 
\beq\label{eq:cutoffconnection}
\nabla^{<\lambda}:=p_\lambda\circ \nabla^{\pi_*(\bS\otimes V)} \circ p_\lambda
\eeq 
on each  $\mathcal{H}^{<\lambda}\to U_\lambda$, where $p_\lambda\colon\mathcal{H}|_{U_\lambda}\to \mathcal{H}^{<\lambda}$ is the fiberwise orthogonal projection. By~\cite[Proposition~9.10]{BGV}, $p_\lambda$ is a smoothing operator, and hence the resulting operator $\nabla^{<\lambda}$ on the finite-rank vector bundle $\mathcal{H}^{<\lambda}\to U_\lambda$ is smooth. Finally, since $p_\lambda$ and $\nabla^{\pi_*(\bS\otimes V)}$ are $\Cl_{-n}$-linear, the resulting connection $\nabla^{<\lambda}$ on $\mathcal{H}^{<\lambda}$ is $\Cl_{-n}$-linear. 

For (ii), we take the differential form
\beq\label{eq:indexdensity}
u^{-n/2}\int_{X/B} \Ch(\nabla^V)\wedge \hat{A}(X/B) \in \Omega^\bullet_\cl(B;\R[\alpha,\alpha^{-1}])
\eeq
where 
\beq\label{eq:Ahatfam}
\hat{A}(X/B):=\det{}^{1/2}\left(\frac{u\cdot (\nabla^{X/B})^2}{\sinh(u\cdot (\nabla^{X/B})^2)}\right)\in \Omega(X;\R[\alpha,\alpha^{-1}])^0
\eeq
is the $\hat{A}$-form of the vertical tangent bundle of $\pi\colon X\to B$ constructed from the Levi-Civita connection $\nabla^{X/B}$ for the family. The integral in~\eqref{eq:indexdensity} is over the fibers of $\pi\colon X\to B$. Since $\Ch(\nabla^V)$ and $\hat{A}(X/B)$ have total degree zero, \eqref{eq:indexdensity} has total degree~$-n$. 

The construction of the data (iii) is a bit more involved, but essentially follows the same arc as Bismut's proof of the local index theorem. Indeed, we shall construct differential forms $\eta_\lambda\in \Omega^\bullet(U_\lambda;\R[\alpha,\alpha^{-1}])$ of total degree~$n-1$ satisfying 
\beq
&&\Ch_{-n}(\nabla^{<\lambda})+d\eta_\lambda=\left(u^{-n/2}\int_{X/B} \Ch(\nabla^V)\wedge \hat{A}(X/B)\right)|_{U_\lambda} \in \Omega^\bullet_\cl(U_\lambda;\R[\alpha,\alpha^{-1}]).\label{eq:etaform}
\eeq
Using that the Pontryagin character is the Chern character of the complexification, the $\eta_\lambda$ come from applying Morita equivalences to reduce the problem to Bismut's known construction for the families index over~$\C$, see~\cite[Chapter~10]{BGV} for $n$ even and~\cite[\S{II.f}]{BismutFreed2} for $n$ odd. Our particular method of constructing $\eta_\lambda$ is an adaptation of~\cite[\S7.2]{FreedLott}. 

Recall the invertible bimodule $\C^{1|1}$ implementing the Morita equivalence between $\cCl_2$ and $\C$ from Example~\ref{ex:cMorita2}. Taking tensor powers, we obtain invertible bimodules
$$
M(n):=\left\{\begin{array}{ll} {\rm invertible} \ \C-\cCl_{-n} \ {\rm bimodule} \ {\rm for} \  n \ {\rm even} \\ {\rm invertible} \ \cCl_{-1}-\cCl_{-n} \ {\rm bimodule} \ {\rm for} \ n \ {\rm odd}.\end{array}\right.
$$
Consider the Fr\'echet vector bundles on $B$,
\beq
&&M(n)\otimes_{\cCl_{-n}} \pi_*(\bS\otimes V_\C )\simeq  \pi_*(M(n)\otimes_{\cCl_{-n}} \bS\otimes V_\C )\simeq \pi_*(\slashed{S}\otimes_{\R} V)\label{eq:slashedspin}
\eeq
as a bundle of $\cCl_0=\C$-modules when $n$ is even and $\cCl_{-1}$-modules when $n$ is odd. The last isomorphism in~\eqref{eq:slashedspin} translates between Clifford linear Dirac operators and Dirac operators built out of spinor representations (see~\cite[\S{II.7}]{LM}): when $n$ is even, we identify~\eqref{eq:slashedspin} with the direct image of the standard $V\otimes \C$-twisted complex spinor bundle 
$$
\slashed{S}:=P\times_{\Spin(n)} \Delta \simeq M(n)\otimes_{\cCl_{-n}} (\bS\otimes \C)
$$ 
where $\slashed{S}$ is constructed as the associated bundle for the $\Z/2$-graded complex spin representation $\Delta$ (compare~\ref{eq:Riemspinbundle}).  Similarly, when $n$ is odd we identify~\eqref{eq:slashedspin} with the $\cCl_{-1}$-module of $V_\C$-twisted complex spinors, again denoted $\slashed{S}\simeq M(n)\otimes_{\cCl_{-n}} \bS_\C$ above. Under these identifications, the operator 
\beq
D:=\id_{M(n)}\otimes_{\cCl_{-n}}  (\slashed{D}\otimes V_\C)\label{eq:usualDirac}
\eeq
can be identified with the standard $V_\C$-twisted Dirac operator acting on $\slashed{S}\otimes_\R V$. Furthermore, the Clifford linear Bismut superconnection~\eqref{eq:Bismutsuper} becomes the standard Bismut superconnection,
$$
\A_{\rm Bis}=M(n)\otimes_{\cCl_{-n}}\otimes \B,
$$ 
e.g., as defined in~\cite[Proposition~9.10]{BGV} for $n$ even, and~\cite[\S{II.f}]{BismutFreed2}  for $n$ odd. 

For each $\lambda\in \R_{>0}$ define the finite-rank complex vector bundle over~$U_\lambda$ 
\beq
M(n)\otimes_{\cCl_{-n}} (\mathcal{H}^{<\lambda}_\C)=:\left\{\begin{array}{ll}  \mathcal{H}^{<\lambda}_0 & n\ {\rm even} \\ \mathcal{H}^{<\lambda}_1 & n \ {\rm odd} \end{array}\right. 
\nonumber
\eeq
whose fibers are $\C=\cCl_0$-modules when $n$ is even and $\cCl_{-1}$-modules when $n$ is odd.  For $n$ even, the $\Cl_{-n}$-linear connection on $\mathcal{H}^{<\lambda}$ determines a $\cCl_0$-linear connection $\nabla^{\mathcal{H}^{<\lambda}_0}$ on $\mathcal{H}^{<\lambda}_0$, and for $n$ odd we obtain a $\cCl_{-1}$-linear connection $\nabla^{\mathcal{H}^{<\lambda}_1}$ on $\mathcal{H}^{<\lambda}_1$. Furthermore, using~\eqref{eq:Botttrace} to compute the Clifford super dimension of the modules $M(n)$, the Clifford super trace in~\eqref{eq:etaform} can computed by the complex Clifford super trace,
$$
\Ch_{-n}(\nabla^{<\lambda})=
\left\{\begin{array}{ll} u^{-n/2}\Ch_0(\nabla^{\mathcal{H}^{<\lambda}_0}) & n \ {\rm even} \\ u^{-(n-1)/2}\Ch_{-1}(\nabla^{\mathcal{H}^\lambda_1}) & n \ {\rm odd} \end{array}\right.
$$
We will use the notation
$$
\mathcal{H}_\C:=(M(n)\otimes_{\cCl_{-n}} \mathcal{H}\otimes\C)= L^2(\pi_*(\slashed{S}\otimes_\R V))
$$
to denote the continuous complex Hilbert bundle that is related to the real Hilbert bundle $\mathcal{H}=L^2(\pi_*(\mathbb{S}\otimes V_\C))$ by a Morita equivalence. 

Next we apply an argument from~\cite[\S7.2]{FreedLott} to construct $\eta_\lambda$. For $l\in \{0,1\}$, let $\Pi \mathcal{H}^{<\lambda}_l$ denote the parity reversal of the super vector bundle $\mathcal{H}^{<\lambda}_l$ and consider the (continuous) Hilbert bundle $\mathcal{H}_\C \oplus \Pi \mathcal{H}^{<\lambda}_l$ on $U_\lambda$ with the $\R$-family of self-adjoint unbounded operators
\beq
D^\lambda(\alpha)=\left[\begin{array}{cc} D & \alpha i_\lambda \\ \alpha p_\lambda & 0 \end{array}\right],\quad \alpha \in \R
\eeq
where $D$ is as in~\eqref{eq:usualDirac}, $p_\lambda \colon \mathcal{H}_\C \to \mathcal{H}^{<\lambda}_l$ is the spectral projection and $i_\lambda\colon \mathcal{H}^{<\lambda}_l\to \mathcal{H}_\C$ is the inclusion. Fix a real-valued function $\alpha(s)\in C^\infty(\R)$ that vanishes near $0\in \R$ and has $\alpha(s)=1$ for $s\ge 1$. Then define the 1-parameter family of superconnections
$$
\A_{{\rm Bis},\lambda}(s)=s D^\lambda(\alpha(s))+\nabla^{\pi_*(\slashed{S}\otimes V_\C)}\oplus \nabla^{\mathcal{H}^{<\lambda}_l}-s^{-1}\frac{c(T)}{4}\oplus 0_{\mathcal{H}^{<\lambda}_l}.
$$
For each $s$, this is a superconnection adopted to a finite-rank deformation of the family of Dirac operators~$D$. Then for $s$ near zero $\A_{{\rm Bis},\lambda}(s)=\A_{\rm Bis}\oplus \nabla^{\Pi \mathcal{H}^{<\lambda}_0}$ and
$$
\lim_{s\to 0} \Ch_{-l}(\A_{{\rm Bis},\lambda}(s))=u^{-n/2}\int_{X/B}\Ch(V)\wedge \hat{A}(X/B)-\Ch_{-l}(\nabla^{\mathcal{H}^{<\lambda}_l}),
$$
by standard properties of the Bismut superconnection; for $l=0$ see~\cite[Theorem~10.23]{BGV} and for $l=1$ see~\cite[Theorem~2.10]{BismutFreed2}. For $\alpha(s)>0$, $D_\lambda(\alpha(s))$ is a family of invertible operators by~\cite[Lemma~7.18]{FreedLott}. Then by~\cite[Theorem 9.26]{BGV}, 
$$
\lim_{s\to \infty} \Ch(\A_{{\rm Bis},\lambda}(s))=0.
$$
Therefore, the differential form (compare~\cite[Equation 7.23]{FreedLott})
$$
\eta_\lambda:=u^{-n/2} \int_0^\infty \sTr(u\frac{d\A_{{\rm Bis},\lambda}(s)}{ds} e^{u \A_{{\rm Bis},\lambda}(s)^2})ds\in \Omega^\bullet(U_\lambda;\R[\alpha,\alpha^{-1}])
$$
is well-defined and satisfies~\eqref{eq:etaform}. 
The $\{\nabla^{<\lambda},\eta_\lambda\}$ provide the data for a differential refinement of the $\Cl_{-n}$-bundle $\{\mathcal{H}^{<\lambda},g_{\lambda\mu}\}$, where to $U_\sigma$ we assign the restriction of $\{\mathcal{H}^{<\lambda},\nabla^{<\lambda},\eta_\lambda\}$ for $\lambda\in \sigma\subset \Lambda\subset \R_{>0}$ the largest element. It remains to verify the compatibility conditions~\eqref{eq:determinesCh} and~\eqref{eq:compatibleChern}, so that the $\eta_\lambda$ determine the desired \v{C}ech--de~Rham coboundary. The description~\eqref{eq:etaform} identifies the global form and verifies condition~\eqref{eq:determinesCh}. The integrals for $\eta_\lambda$ and $\eta_\mu$ differ by a super trace over a finite-dimensional subspace on which the restriction of the Dirac operator is invertible, coming down to the exact sequence of metrized super vector bundles
$$
0\to \mathcal{H}^{<\lambda}\hookrightarrow \mathcal{H}^{<\mu}\to \mathcal{H}^{(\lambda,\mu)}\to 0,
$$
and the Chern--Simons form of the associated connections. As the first arrow is a stable equivalence, the discrepancy between the integrals $\eta_\lambda$ and $\eta_\mu$ is measured by the same construction as in Lemma~\ref{lem:seCS}, and~\eqref{eq:compatibleChern}  is satisfied by construction.
This completes the construction of the differential index
$$
\widehat{\pi}_!(V,\nabla):=\{\mathcal{H}^{<\lambda}\to U_\lambda,\nabla^{<\lambda},\eta_\lambda,g_{\lambda\mu},e_{\lambda\mu}\} \in \dKO^{-n}(B)
$$
for  $[V,\nabla]\in \dKO^0(X)$. 

Just as in Proposition~\ref{thm:Cliffordindex}, the only choice involved in constructing $\pi_!(V,\nabla)$ is the choice of cover $\{U_\lambda\}$ of $B$ for $\Lambda\subset \R_{>0}$. A concordance between $\Cl_{-n}$-bundles for different choices was constructed in Proposition~\ref{thm:Cliffordindex} via a stable equivalence of $\Cl_{-n}$-bundles. By Corollary~\ref{cor:trivsuper}, the Chern--Simons form for this stable equivalence is exact, and hence the differential class $\pi_!(V,\nabla)$ does not depend on the choice of cover. 

Finally, Proposition~\ref{thm:Cliffordindex} shows that the differential index refines the families analytic index. This completes the proof. 
\ep

\bibliographystyle{amsalpha}
\bibliography{references}

\end{document}